\documentclass[12pt]{article}

\usepackage{amsmath}
\usepackage{amssymb}
\usepackage{amsfonts}
\usepackage{a4}
\usepackage{graphicx}
\usepackage[abs]{overpic}
\usepackage{caption}
\usepackage{wrapfig}
\usepackage{float}
\usepackage{graphicx}

\newtheorem{theorem}{Theorem}[section]
\newtheorem{lemma}[theorem]{Lemma}

\newtheorem{corollary}[theorem]{Corollary}

\newtheorem{rem}[theorem]{Remark}
\newtheorem{example}[theorem]{Example}
\newcommand{\Proof}{\par\noindent{\em Proof. }}
\newcommand{\eop}{\nopagebreak\hspace*{\fill}$\Box$\smallskip}

\newcommand{\N}{\Bbb N}
\newcommand{\Z}{\Bbb Z}
\newcommand{\R}{\Bbb R}

\def\Id{\mathbf{Id}}
\def\id{\mathbf{id}}

\def\eps{\varepsilon}

\def\e{\mathbf{e}}

\def\dist{\operatorname{dist}}

\def\XXint#1#2#3{{\setbox0=\hbox{$#1{#2#3}{\int}$}
     \vcenter{\hbox{$#2#3$}}\kern-.5\wd0}}

\usepackage{color}

\numberwithin{equation}{section}

\begin{document}

\begin{center}
\begin{Large}
{\bf {A quantitative geometric rigidity result in SBD}}
\end{Large}
\end{center}

\begin{center}
\begin{large}
Manuel Friedrich\footnote{Faculty of Mathematics, University of Vienna, 
Oskar-Morgenster-Platz 1, 1090 Vienna, Austria. {\tt manuel.friedrich@math.uni-augsburg.de}}
and Bernd Schmidt\footnote{Universit{\"a}t Augsburg, Institut f{\"u}r Mathematik, 
Universit{\"a}tsstr.\ 14, 86159 Augsburg, Germany. {\tt bernd.schmidt@math.uni-augsburg.de}}\\
\end{large}
\end{center}

\begin{center}
\today
\end{center}
\bigskip

\begin{abstract}
We present a quantitative geometric rigidity estimate for special functions of bounded deformation in a planar setting generalizing a result by Friesecke, James and M\"uller 
for Sobolev functions obtained in nonlinear elasticity theory and a qualitative piecewise rigidity result by Chambolle, Giacomini and Ponsiglione for brittle materials which do not store elastic energy. We show that for each deformation there is an associated triple consisting of a partition of the domain, a corresponding piecewise rigid motion being constant on each connected component of the cracked body and a displacement field measuring the distance of the deformation from the piecewise rigid motion. We also present a related estimate in the geometrically linear setting which can be interpreted as a `piecewise Korn-Poincar\'e inequality'.
\end{abstract}
\bigskip

\begin{small}


\end{small}

\tableofcontents

\section{Introduction}

It is a subtle problem in mathematical analysis to infer global properties of a function $u$ from conditions on its derivative $\nabla u$ given in terms of partial differential relations such as $\nabla u \in K$ or approximate relations such as $\dist(\nabla u, K) \ll 1$, where $K$ denotes a specific set of matrices. In particular, constraining $\nabla u$ to be in, or close to, the set $K = SO(d)$ of rigid motions, one is led to the question to what extend such a pointwise (approximate) isometry constraint has the global consequence of rendering $u$ itself (approximately) rigid. As will be detailed below, notably the last decades have witnessed a tremendous progress in establishing such geometric rigidity results; classical theorems for smooth functions have been extended to Sobolev functions and even sharp rigidity estimates have been derived for such functions. In this article we address the problem of deriving a quantitative rigidity estimate beyond the setting in Sobolev spaces, specifically, allowing for functions with jump discontinuities. As such a lack of regularity impedes a direct extension, our main rigidity result has to be formulated in a considerably more complex way. Moreover, major challenges arise in our framework from the fact that the distributional derivative of the mappings under consideration is barely a measure and from the necessity to gain control over both bulk and surface contributions. 

Our main motivation comes from variational fracture mechanics. Since the pioneering work of Griffith \cite{Griffith:1921} the propagation of crack is viewed as the result of a competition between the surface energy and the reduction of bulk energy during an infinitesimal increase of the cracked region.   Based on this idea Francfort and Marigo \cite{Francfort-Marigo:1998} have introduced an energy functional comprising elastic bulk and surface contributions in order to tackle problems in fracture mechanics with variational methods, where the displacements and crack paths are determined from an energy minimization principle.

To simplify the mathematical description, problems in this context are often studied in the case of anti-planar shear (see e.g. \cite{DalMaso-Toader:02, Francfort-Larsen:2003}) or in the realm of linearized elasticity (see e.g. \cite{Bellettini-Coscia-DalMaso:98, Chambolle:2003, Chambolle:2004, Focardi-Iurlano:13, Iurlano:13, SchmidtFraternaliOrtiz:2009})  since  such models are usually significantly easier to treat as their nonlinear counterparts. In fact, in the regime of finite elasticity the energy density of the elastic contributions is genuinely geometrically nonlinear due to frame indifference rendering the problem highly non-convex. Consequently, in contrast to linear models already the fundamental question if minimizing configurations for given boundary data exist at all is a challenging problem.

To gain a deeper understanding of nonlinear models in fracture mechanics it is therefore desirable to identify an effective linear theory and in this way to rigorously show  that in the small displacement regime the neglection of effects arising from the non-linearities is a good approximation of the problem. Indeed, for elastic bodies not exhibiting cracks the passage from nonlinear to linearized models is by now well understood via $\Gamma$-convergence (cf. \cite{DalMasoNegriPercivale:02, Schmidt:08}). It turns out that a fundamental issue in this context is the derivation of suitable rigidity estimates  which, based on the deformation of a material, allow to control an  associated  infinitesimal displacement field measuring the distance from a rigid motion and being the essential quantity  on which the linearized elastic energy depends.

Rigidity estimates  have a long history going back to the fundamental result of Liouville which states that a smooth function has to be an affine mapping if its gradient is a rotation everywhere. Various generalizations of this classical qualitative theorem in the realm of nonlinear elasticity theory have appeared over the last decades (see e.g. \cite{John:1961, Reshetnyak:1961}). For brittle materials the problem is more subtle as additional difficulties arise from the fact that the body might be disconnected by the jump set into various components.  Chambolle, Giacomini and Ponsiglione \cite{Chambolle-Giacomini-Ponsiglione:2007} recently showed that also in this setting a Liouville-type result holds and that the body behaves piecewise rigidly. In fact, under the constraint that the material does not store elastic energy the only possibility that global rigidity can fail is that the body is divided into various parts each of which subject to a different rigid motion.

However, the above mentioned results fall short of being useful for the investigation of variational models due to the  restrictive constraint on the deformation gradient. The fundamental step towards quantitative results was a geometric rigidity estimate by Friesecke, James and M\"uller \cite{FrieseckeJamesMueller:02} which states that, loosely speaking, if the deformation gradient of an $H^1$-function is close to the set of rotations (e.g. in an $L^2$ sense), then it is in fact close to one single rotation. This result provides the essential relation between the deformation and a corresponding displacement field and allows to establish a compactness result  for a sequence of displacements with uniformly bounded elastic energy.

Whereas this estimate in elasticity theory was generalized to various settings including \cite{Conti-Dolzmann-Muller:14, Muller-Scardia-Zeppieri:14}, to the best of our knowledge a corresponding general estimate for brittle materials has not yet been established. The farthest reaching result in this direction seems to be a recent contribution by Negri and Toader \cite{NegriToader:2013} where rigidity estimates are provided in the context of quasistatic evolution for a restricted class of admissible cracks. In particular, in their model the different components of the jump set are supposed to have a least positive distance rendering the problem considerably easier. In fact, one can essentially still employ the result in \cite{FrieseckeJamesMueller:02} and the specimen cannot be separated into different parts  effectively leading to a simple relation between the deformation and the displacement field.

The goal of the present work is the derivation of a new kind of quantitative geometric rigidity estimate in the framework of geometric measure theory without  any a priori assumptions on the deformation and the crack geometry, i.e we treat a full free discontinuity problem in the language of Ambrosio and De Giorgi \cite{DeGiorgi-Ambrosio:1988}. We call this estimate for brittle materials, which we establish in a planar setting, an SBD-rigidity result as it is formulated in terms of \emph{special functions of bounded deformation} (see \cite{Ambrosio-Coscia-Dal Maso:1997, Bellettini-Coscia-DalMaso:98}). The result may be seen as a suitable combination of the aforementioned estimate for elastic materials \cite{FrieseckeJamesMueller:02} and the qualitative result in \cite{Chambolle-Giacomini-Ponsiglione:2007},   being tailor-made for general Griffith models where both energy forms are coexistent. 

The rigidity result provides the relation between the deformation of a brittle material and the associated displacements. Whereas in elasticity theory there is a simple connection between these two objects, in the present context the description is rather complicated since the deformation is related to a triple consisting of a partition of the domain, a corresponding piecewise rigid motion being constant on each connected component of the cracked body and a displacement field which is defined separately on each piece of the specimen. The result in the present work proves to be the fundamental ingredient to identify an effective linearized theory. For a detailed analysis of compactness results and the derivation  of linearized Griffith models from nonlinear energies via $\Gamma$-convergence in a small strain limit we refer to the subsequent paper \cite{Friedrich:15-2}.

One essential point in the analysis is the derivation of an inequality for the symmetric part of the gradient. We also see that in general it is not possible to gain control over the full gradient which is not surprising as there is no analogue of Korn's inequality for SBV functions. Consequently, the result is naturally an SBD estimate. In addition, we provide an $L^2$-bound for the configurations measuring the distance of the deformation itself from a piecewise rigid motion. In contrast to the setting in elasticity theory this is highly nontrivial as Poincar\'e's inequality cannot be applied due to the possibly present complicated crack geometry. Consequently, our findings are not only interesting  in the realm of finite elasticity, but also in a geometrically linear setting and can be interpreted as a `piecewise Korn-Poincar\'e inequality'. Moreover, we remark that our main estimate can only be established under the additional condition that we admit an arbitrarily small modification of the deformation.

The derivation of the main result is very involved as among other things one has to face the problems that (1) the body might be disconnected by the jump set, (2) the body might be still connected but only in a small region where the elastic energy is possibly large, (3) the crack geometry might become extremely complex due to relaxation of the elastic energy by oscillating crack paths and infinite
crack patterns occurring on different scales.  The common difficulty  of all these phenomena is the possible high irregularity of the jump set. Even if one can assume that the domain can be decomposed into different sets with Lipschitz boundary (e.g. by a density argument), there are no uniform bounds on the constants of several necessary inequalities such as the Poincar\'e and Korn inequality and the rigidity estimate \cite{FrieseckeJamesMueller:02}.    

To avoid further complicacies of technical nature concerning the topological structure of cracks in higher dimensions and to concentrate on the essential difficulties arising from the frame indifference of the energy density, we will tackle the problem in a planar setting  with isotropic crack energies. However, we believe that our results can be extended to  anisotropic surface terms and that the proof provides the principal techniques being necessary to establish the result in arbitrary space dimension. In fact, many arguments are valid also in dimension $d \ge 3$ and we hope  that our methods, in particular the modification scheme for deformations and jump sets, may also contribute to solve related problems in the future.  One of the essential reasons why we restrict ourselves to the two-dimensional framework is the usage of a Korn-Poincar\'e-type inequality (see \cite{Friedrich:15-1}) which was only established in a planar setting due to a lot of technical difficulties concerning the jump set geometry.

The paper is organized as follows. In Section \ref{rig-sec: main} we present the main results about geometric rigidity in SBD and also state a corresponding estimate in the geometrically linear setting which is interesting on its own and considerably simpler to prove  than its nonlinear counterpart.  As the proof  is very long and technical, we give an  overview and highlight the principal strategies for the convenience of the reader in Section \ref{rig-sec: sub, proofoverview}.

Section \ref{rig-sec: pre} is devoted to some preliminaries. We first recall the definition of special functions of bounded variation and discuss basic properties. Then we recall a (local) Korn-Poincar\'e-type inequality in SBD (see \cite{Friedrich:15-1} and Section \ref{rig-sec: subsub,  mainproof}) which measures the distance of the displacement field from an infinitesimal rigid motion in terms of the elastic energy. It turns out that this inequality is one of the key ingredients to derive our main result which can be compared with the fact that in elasticity theory the linearized rigidity estimate, called Korn's inequality (see \cite{ContiFaracoMaggi:2005}), is one of the fundamental steps to establish the geometrically nonlinear result in \cite{FrieseckeJamesMueller:02}. In fact, as a first approach to the main result it is convenient to replace the nonlinear problem by such a linearized version which is significantly easier since (1) the estimate only involves the function itself and not its derivative and (2) the set of infinitesimal rigid motions is a linear space in contrast to $SO(2)$.  

Afterwards we recall the geometric rigidity result by Friesecke, James, M\"uller \cite{FrieseckeJamesMueller:02} and carry out a careful analysis how the involved constant depends on the shape of the domain.  At this point we notice that easy counterexamples to rigidity estimates in SBD can be constructed if one does not admit a small modification of the deformation (see Section \ref{rig-sec: sub, pre}).

In Section \ref{rig-sec: subsub, rig-prep} we introduce a procedure to modify sets. In this context, we particularly have to assure that we can control the size and the shape of the jump sets.

The rest of the paper contains the main proof of the SBD-rigidity estimate. The main strategy of the proof is to establish local rigidity results on cells of mesoscopic size (Section \ref{rig-sec: sub, first-weak}) which together with the Korn-Poincar\'e inequality allows to replace the deformation by a modification where the least length of the crack components has increased (Section \ref{rig-sec: sub, local}). Repeating the arguments on various mesoscopic scales becoming gradually larger it is possible to show that the modified deformation behaves rigidly on each connected component of the domain (Section \ref{rig-sec: sub, proof}). 

The fact that we analyze the problem on different length scales is indispensable to understand specific size effects correctly such as the accumulation of crack patterns on certain scales. Moreover, we briefly note that similarly as in \cite{FriedrichSchmidt:2014.1}  a mesoscopic localization technique proves to be useful to tackle problems in the framework of brittle materials as hereby effects arising from the bulk and the surface contributions can be separated.

Basically, this is enough the establish the requirements for compactness results in the space of SBD functions (cf. \cite{DalMaso:13}). However, as we are also interested in the derivation of effective linearized models (cf. \cite{Friedrich:15-2}), we have to assure that  we do not change the total energy of the deformation during the modification procedure. In particular, for the surface energy this is a subtle problem and in Section \ref{rig-sec: sub, proof-main} a lot of effort is needed to show that the modified configurations can be constructed in a way such that the crack length does not increase substantially.

\section{The main result and overview of the proof}\label{rig-sec: main}

In this section we present our main rigidity estimates in the framework of brittle materials and give and overview of the proof strategies. 

\subsection{The main setting}
Let $\Omega \subset \R^2$ open,  bounded with Lipschitz boundary and for $M>0$ we define 
\begin{align}\label{rig-eq: SBVfirstdef}
SBV_M(\Omega) = \Big\{ y \in SBV(\Omega,\R^2):  \Vert \nabla y\Vert_{\infty} \le M, \ {\cal H}^1(J_y) < + \infty \Big\}.
\end{align}
For the definition and properties of the space $SBV(\Omega,\R^2)$, frequently abbreviated as $SBV(\Omega)$ hereafter,  we refer to Section \ref{rig-sec: sub, bd}.

 Let $W:\R^{2 \times 2} \to [0,\infty)$ be a frame-indifferent stored energy density with $W(F) = 0$ iff $F \in SO(2)$. Assume that $W$ is continuous, $C^3$ in a neighborhood of $SO(2)$ and scales quadratically at $SO(2)$ in the direction perpendicular to infinitesimal rotations. In other words, we have $W(F) \ge c\dist^2(F,SO(2))$ for all $F \in \R^{2 \times 2}$ and a positive constant $c$.  For $\eps >0$ define the Griffith-energy $E_\eps : SBV_M(\Omega) \to [0,\infty)$ by
\begin{align}\label{rig-eq: Griffith en}
E_\eps(y) =  \frac{1}{\eps}\int_\Omega W(\nabla y(x)) \,dx + {\cal H}^1(J_y).
\end{align}
The main goal of the work at hand is the derivation of uniform rigidity estimates for configurations with $E_\eps(y) \le C$. Performing the passage to the small strain limit $\eps \to 0$ we have to face major challenges including (1) difficulties concerning the coercivity of the functionals due to the frame indifference of the energy density and (2) the possible high irregularity of the jump set rendering the problem subtle from an analytical point of view.

We briefly note that we can also treat inhomogeneous materials where the energy density has the form $W: \Omega \times \R^{2 \times 2} \to [0,\infty)$. Moreover, it suffices to assume $W \in C^{2,\alpha}$, where $C^{2,\alpha}$ is the H\"older space with exponent $\alpha >0$. In the context of discrete systems the small parameter $\eps$, denoting the order of the elastic energy in our model, represents the typical interatomic
distance (compare  \eqref{rig-eq: Griffith en}  with, e.g.,  the Griffith functionals in \cite{FriedrichSchmidt:2014.1, FriedrichSchmidt:2014.2}). Having also applications to discrete systems in mind, we will sometimes refer to $\eps$  as the `atomic length scale'.

Observe that $M$ may be chosen arbitrarily large (but fixed) and therefore the constraint $\Vert \nabla y\Vert_{\infty} \le M$ is not a real restriction as we are interested in the small displacement regime in the regions of the domain where elastic behavior occurs. The uniform bound on the absolute continuous part of the gradient is indeed natural when dealing with discrete energies where the corresponding deformations are piecewise affine on cells of microscopic size (see e.g. \cite{Braides-Gelli:2002-2, FriedrichSchmidt:2014.2}). The condition essentially assures that the elastic energy cannot concentrate on scales being much smaller than $\eps$. This observation already shows that the atomic length scale plays an important role in our analysis since the system shows remarkably different behavior on scales smaller and larger than
the atomistic unit.

For later we also introduce a relaxed energy functional. For $\rho > 0$, $\eps > 0$ and $U \subset \Omega$ define $f_\eps^\rho(x) = \min\lbrace\frac{x}{\sqrt{\eps}\rho} ,1 \rbrace$ and
\begin{align}\label{rig-eq: Griffith en2}
E_\eps^\rho(y,U) =  \frac{1}{\eps}\int_U W(\nabla y(x)) \,dx + \int_{J_y \cap U} f_\eps^\rho(|[y](x)|)\,d{\cal H}^1(x).
\end{align}
Clearly, we have $E_\eps^\rho(y,U) \le E_\eps(y)$ for all $y \in SBV_M(\Omega)$ and $U \subset \Omega$. 

\subsection{Rigidity estimates}

We first observe that for configurations with uniform bounded energy $E_\eps(y_\eps)$ the absolute continuous part of the gradient satisfies $\nabla y_\eps \approx SO(2)$ as the stored energy density is frame-indifferent and minimized on $SO(2)$. Assuming that $y_\eps \to y$ in $L^1$, one can show that $\nabla y \in SO(2)$ a.e. applying lower semicontinuity results for SBV functions (see \cite{Kristensen:1999}) and the fact that the quasiconvex envelope of $W$ is minimized exactly on $SO(2)$ (see \cite{Zhang:2004}). 

A classical result due to Liouville states that a smooth function $y$ satisfying the constraint $\nabla y \in SO(2)$ is a rigid motion. In the theory of fracture mechanics global rigidity can fail if the crack disconnects the body. More precisely, Chambolle, Giacomini and Ponsiglione have proven that for configurations which do not store elastic energy (i.e. $\nabla y \in SO(2)$ a.e.) and have finite Griffith energy (i.e. ${\cal H}^1(J_y) < +\infty$) the only way that rigidity may fail is that the body is divided into at most countably many parts each of which subject to a different rigid motion (see \cite{Chambolle-Giacomini-Ponsiglione:2007}). 

Clearly, it is desirable to establish an appropriate quantitative version of this qualitative statement. In nonlinear elasticity such quantitative estimates are available forming one of the starting points of our analysis. Friesecke, James and M\"uller (see \cite{FrieseckeJamesMueller:02} and Theorem \ref{rig-th: geo rig} below) have  extended the classical Liouville results and showed that, loosely speaking, if the deformation gradient is close to $SO(2)$ (in $L^2$), then it is in fact close to one single rotation $R \in SO(2)$ (in $L^2$).

The overall goal of this work is to `combine' the rigidity results of the pure elastic and pure brittle regime in order to derive a rigidity estimate for general Griffith functionals \eqref{rig-eq: Griffith en} where both energy forms are coexistent. As a preparation recall  the definition of the  \emph{perimeter} $P(E,\Omega)$  of a set $E \subset \R^2$ in $\Omega$ (see \cite[Section 3.3]{Ambrosio-Fusco-Pallara:2000}) and recall that we say that a partition ${\cal P} = (P_j)_j$ of $\Omega$ is called a \textit{Caccioppoli partition} of $\Omega$ if $\sum_j P(P_j,\Omega) < + \infty$. Let $\Omega_\rho = \lbrace x\in\Omega: \dist(x, \partial \Omega) > C\rho \rbrace$ for $\rho>0$ and for some sufficiently large constant $C$.

\begin{theorem}\label{rig-th: rigidity}
Let $\Omega \subset \R^2$  open, bounded with Lipschitz boundary. Let $M>0$ and $0 < \eta, \rho \ll 1$. Then there is a constant $C=C(\Omega,M,\eta)$  and a universal $c>0$ such that the following holds  for $\eps >0$ small enough: \\
For each $y \in SBV_{M}(\Omega) \cap L^2(\Omega)$ with ${\cal H}^{1}(J_y) \le M$ and $\int_\Omega \dist^2(\nabla y,SO(2) )  \le M\eps$, there is an open set $\Omega_y$ with $|\Omega\setminus\Omega_y| \le C\rho$, a modification  $\hat{y} \in SBV_{cM}(\Omega) \cap L^2(\Omega)$  with  $\Vert \hat{y} - y \Vert^2_{L^2(\Omega_y)} +  \Vert \nabla \hat{y} - \nabla y \Vert^2_{L^2(\Omega_y)}\le C\eps\rho$  and 
\begin{align}\label{rig-eq: energy le}
E_\eps^\rho(\hat{y},\Omega_\rho) \le E_\eps(y) + C\rho
\end{align}
with the following properties: We find  a Caccioppoli partition ${\cal P} = (P_j)_j$ of $\Omega_\rho$ with $\sum_j P(P_j,\Omega_{\rho}) \le C$ and for each $P_j$ a corresponding rigid motion $R_j \, x +c_j$, $R_j \in SO(2)$ and $c_j \in \R^2$, such that the function $u: \Omega \to \R^2$ defined by
\begin{align}\label{rig-eq: u def2} 
u(x) := \begin{cases} \hat{y}(x) - (R_j\,\,x +c_j) & \ \ \text{ for } x \in P_j \\
                      0                      & \ \ \text{ for } x \in \Omega \setminus \Omega_\rho \end{cases}
\end{align}
satisfies the estimates
\begin{align}\label{rig-eq: main properties2}
\begin{split}
(i) & \ \, {\cal H}^{1}(J_u) \le C, \ \  \ \ \  \ \ \  \ \ \ \   \ \ \   \ \  \ \ \  \ \ \  \ \  \ \ (ii) \,  \ \Vert u\Vert^2_{L^2(\Omega_\rho)} \le \hat{C}\eps, \\
(iii) & \ \, \sum\nolimits_j \Vert e(R^T_j \nabla u)\Vert^2_{L^2(P_j)} \le \hat{C}\eps,  \  \ \ \  \  \ \ \,  (iv)  \ \, \Vert \nabla u\Vert^2_{L^2(\Omega_\rho)} \le \hat{C}\eps^{1-\eta} 
\end{split}
\end{align}
for some constant $\hat{C}=\hat{C}(\rho)$, where $e(G) = \frac{G + G^T}{2}$ for all $G \in \R^{2 \times 2}$.  
\end{theorem}

Whereas in elasticity theory there is a simple connection between the deformation $y$ and the displacement field $u$, in the present context the description is rather complicated since the deformation is related to a triple  $(P_j)_j$, $(R_j,c_j)_j$ and $u$ consisting of a partition, associated piecewise rigid motion and a suitably rescaled displacement field which is defined separately on each piece of the body. The central estimate \eqref{rig-eq: main properties2} provides the fundamental ingredients to establish a corresponding compactness result (see  \cite{Friedrich:15-2}) by employing a GSBD compactness result proved in \cite{DalMaso:13}. 

We remark that   this estimate might be wrong without allowing for a small modification of the deformation as we show by way of example in Section \ref{rig-sec: sub, pre}. Moreover, we get a sufficiently strong bound only for the symmetric part of the gradient (see (iii)) which is not surprising due to the fact that there is no analogue of Korn's inequality in SBV. However, there is at least a weaker bound on the total absolutely continuous part of the gradient (see (iv)) which will essentially be needed to derive a $\Gamma$-convergence result in the passage from nonlinear to linearized models in \cite{Friedrich:15-2}. We emphasize that also (ii) is highly nontrivial as Poincar\'e's inequality cannot be applied due to the presence of discontinuity sets. 

\begin{rem}
{\normalfont
(i) The proof of Theorem \ref{rig-th: rigidity} shows that the Caccioppoli partition $(P_j)_j$ is in fact a finite partition. In particular, each $P_j$ is the union of squares of sidelength $\sim \rho$ and thus $|P_j| \ge c\rho$ for all $j$.

(ii) In view of \eqref{rig-eq: energy le} and \eqref{rig-eq: main properties2}(i) one also has 
$$E_\eps(\hat{y}) \le C E_\eps(y).$$ 
Moreover, the estimate \eqref{rig-eq: energy le} can even be refined. Indeed, we obtain (see \eqref{rig-eq: part + crack} below)
$$\sum\nolimits_j \tfrac{1}{2} P(P_j,\Omega_\rho) + \int_{J_{\hat{y}} \setminus \partial P} f_\eps^\rho(|[\hat{y}]|) \,d{\cal H}^1 \le {\cal H}^1(J_y) + c\rho,$$
where $\partial P := \bigcup_j \partial P_j$. Whereas on the boundary of the partition $\partial P$ there is a sharp estimate for the surface energy, the passage to to a relaxed functional in the interior of the sets is necessary due to the possible presence of microcracks accumulating on different mesoscopic scales. 

(iii) The assumption $y \in L^2(\Omega)$ may be dropped. In this case we obtain a slightly weaker approximation of the form $\Vert \hat{y} - y\Vert^2_{L^1(\Omega_y)}\le C\eps\rho$ (cf. the approximation schemes in \cite[Theorem 3.1]{Chambolle-Giacomini-Ponsiglione:2007}, \cite[Theorem 2.3]{Friedrich:15-1}). 

(iv) The approximation preserves an $L^\infty$-bound, i.e. $\Vert y \Vert_\infty \le M$ implies $\Vert \hat{y} \Vert_\infty \le cM$.
}
\end{rem}

\subsection{A piecewise Korn-Poincar\'e inequality}

We now discuss  a variant of Theorem \ref{rig-th: rigidity} in the geometrically linear setting which can be interpreted as a `piecewise Korn-Poincar\'e-inequality in SBD'. Let $\R^{2 \times 2}_{\rm skew} = \lbrace A \in \R^{2 \times 2}: A^T = -A \rbrace$ be the set of skew symmetric matrices. Set 
\begin{align}\label{rig-eq: Fenergy}
F_\eps^\rho(y,U) =  \frac{1}{\eps}\int_U V(e(\nabla u)(x)) \,dx + \int_{J_u \cap U} f_\eps^\rho(|[u]|)\,d{\cal H}^1
\end{align}
 for a coercive quadratic form $V$, i.e. $V(G) \ge c|G|^2$ for $c>0$ and $G \in \R^{2 \times 2}_{\rm sym}$. Furthermore, define $F_\eps = F^0_\eps(\cdot,\Omega)$, where $f^0_\eps \equiv 1$.  For the definition of the space SBD we refer to Section \ref{rig-sec: sub, bd}.

\begin{theorem}\label{rig-th: rigidity_lin}
Let $\Omega \subset \R^2$  open, bounded with Lipschitz boundary. Let $M>0$, and $0 < \rho \ll 1$. Then there is a constant $C=C(\Omega,M)$ such that for $\eps>0$ small enough  the following holds:\\
 For each $u \in SBD^2 (\Omega,\R^2) \cap L^2(\Omega,\R^2)$ with ${\cal H}^{1}(J_u) \le M$ and 
 $$\int_\Omega |e(\nabla u)(x)|^2\,dx \le  M\eps,$$
there is an open set $\Omega_u$ with $|\Omega\setminus\Omega_u| \le C\rho$, a modification $\hat{u}: \Omega \to \R^2$ with $\Vert \hat{u} - u \Vert^2_{L^2(\Omega_u)} + \Vert e(\nabla \hat{u}) - e(\nabla u) \Vert^2_{L^2(\Omega_u)} \le C\rho\eps$  and
$$F_\eps^\rho(\hat{u},\Omega_\rho) \le F_\eps(u) + C\rho$$
with the following properties: We find a Caccioppoli partition ${\cal P} = (P_j)_j$ of $\Omega_\rho$ with $\sum_j P(P_j,\Omega_\rho) \le C$ and for each $P_j$ a corresponding infinitesimal rigid motion $A_j \, x +c_j$, $A_j \in \R^{2 \times 2}_{\rm skew}$ and $c_j \in \R^2$, such that ${\cal H}^{1}(J_{\hat{u}}) \le C$ and   
\begin{align}\label{rig-eq: main properties3}
(i) \  \Vert e(\nabla \hat{u})\Vert^2_{L^2(\Omega_\rho)} \le C\eps, \ \  (ii) \  \sum\nolimits_j \Vert \hat{u} - (A_j\, \cdot -c_j)\Vert^2_{L^2(P_j)} \le \hat{C}\eps.
\end{align}
for some constant $\hat{C} = \hat{C}(\rho)$.  
\end{theorem}

To prove Theorem \ref{rig-th: rigidity_lin} one may essentially follow the proof of Theorem \ref{rig-th: rigidity} with some changes, where altogether the proof  is considerably simpler as a lot of estimates and arguments can be skipped. We again observe that estimate \eqref{rig-eq: main properties3} together with the result of \cite{DalMaso:13} is the fundamental ingredient to establish a compactness result.

\subsection{Overview of the proof}\label{rig-sec: sub, proofoverview}

As the proof of Theorem \ref{rig-th: rigidity} is very long and technical, we present here a short overview  for the convenience of the reader and highlight  the principle proof strategies.
   
The main estimates in the rigidity result (see \eqref{rig-eq: main properties2}) provide bounds for both the displacement field $u$ itself and its derivative. The fundamental ingredient to measure the distance of the function from a rigid motion is a (local) Korn-Poincar\'e-type inequality established in \cite{Friedrich:15-1}. The other key point is then the derivation of an estimate for the symmetric part of the gradient. Using the expansion
\begin{align}\label{rig-eq: over:rig3}
|e(R^T (\nabla y - \Id))|^2 = \dist^2(\nabla y,SO(2)) + O(|\nabla y - R|^4)
\end{align}
and recalling that $\Vert \dist(\nabla y,SO(2)) \Vert^2_{L^2(\Omega)}  \sim \eps$ we see that it suffices to establish an estimate of fourth order. Indeed, also in the proof of the geometric rigidity result in nonlinear elasticity (see \cite{FrieseckeJamesMueller:02}) one first derives a bound for $\Vert \nabla y -R\Vert^4_{L^4(\Omega)}$ to control the symmetric part. The control over the full gradient is then obtained by Korn's inequality.

Clearly, in our framework this rigidity result (see  Theorem \ref{rig-th: geo rig} below) cannot  be applied due to the presence of cracks, in particular $\Omega\setminus J_y$ will generically not be a Lipschitz set. Therefore, by a density argument we again first assume that the jump set is contained in a finite number of rectangle boundaries. A careful quantitative analysis shows that the constant  in Theorem \ref{rig-th: geo rig} depends on the quotient of the diameter of the domain, denoted by $k$, and the minimal distance of two cracks, denoted by $s$. In particular, $C=C(k/s) \sim 1$ if $k \sim s$. Provided that $\frac{k}{s}$ is not too large, the principal strategy will be to show that possibly after a modification  we get $\Vert \nabla y - R\Vert^2_{L^\infty(\Omega)} \le (C(k/s))^{-1}$ which then gives 
\begin{align}\label{rig-eq: over:rig1}
\Vert e(R^T (\nabla y - \Id))\Vert^2_{L^2(\Omega)} \le \eps + (C(k/s))^{-1} \Vert \nabla y -R \Vert^2_{L^2(\Omega)} \le C\eps
\end{align}
by  \eqref{rig-eq: over:rig3} and Theorem \ref{rig-th: geo rig}. Of course, in general we cannot suppose that $\frac{k}{s}$ is not large. Moreover, a global rigidity result may fail due to the separation of the domain by the jump set. Consequently, we will apply the presented ideas on a fine partition of the Lipschitz domain $\Omega$ consisting of squares with diameter $k$. This local result will be used to modify the jump set such that the minimal distance of each pair of cracks increases. Then we can repeat the arguments for a larger $k$. The idea is that after an iterative application of the arguments we obtain an estimate for $k \approx \rho$ which then will provide rigid motions on the connected components of the domain (see \eqref{rig-eq: u def2}) with the desired properties.

In Section \ref{rig-sec: subsub, rig-prep}  we introduce a procedure to modify sets and  conduct a thorough analysis on how to control the size and the shape of the jump sets. 

In Section \ref{rig-sec: sub, first-weak} we construct piecewise constant $SO(2)$-valued mappings approximating the deformation gradient. In each square $Q$ of diameter $k$ we may assume that the elastic energy is bounded by $\sim \eps k$ as otherwise it would be energetically favorable to introduce jumps at the boundary of the square and to replace the deformation in the interior by a rigid motion. (The same technique has been used in the proof of the Korn-Poincar\'e inequality.) Similarly as in \cite{FrieseckeJamesMueller:02} we pass to the harmonic part of the deformation  (denoted by $\hat{y}$) and obtain by the mean value property
\begin{align}\label{rig-eq: over:rig2}
\begin{split}
\Vert \nabla  \hat{y} - R_Q \Vert^2_{L^{\infty}( \hat{Q})} &\le Ck^{-2} \Vert \nabla  \hat{y} - R_Q \Vert^2_{L^{2}(Q)} \\ &
\le C(k/s)k^{-2} \Vert \dist(\nabla y,SO(2))\Vert^2_{L^2(Q)}  \le C(k/s) k^{-1}\eps
\end{split}
\end{align} 
for a suitable $R_Q \in SO(2)$,  where $\hat{Q} \subset Q$ is a slightly smaller square.  Consequently, if we can assure that $\frac{\eps}{k} \le (C(k/s))^{-2}$ we obtain the desired $L^\infty$-bound which allows to derive an estimate of the form \eqref{rig-eq: over:rig1}. We note that for this argument we at least have to assume that $k \gg \eps$ which will be denoted as the  `superatomistic regime' (recall the discussion about the signification of $\eps$ after \eqref{rig-eq: Griffith en}).  

In the subsequent Section \ref{rig-sec: subsub,  est-h1} we show that not only the distance of the derivative from a piecewise rigid motion can be controlled  but also the distance of  the function itself. On the one hand this is essential for \eqref{rig-eq: main properties2}, on the other hand such an estimate is crucial for establishing a modification of the deformation and the jump set. The main idea is to apply the Korn-Poincar\'e-type inequality proved in \cite{Friedrich:15-1} on the function $R^T_Q y - \id$. Major difficulties arise from the facts that  the rotation $R_Q$ may vary from one square to another and  that the inequality derived in \cite{Friedrich:15-1} only provides a local estimate (cf. also Corollary  \ref{rig-cor: kornpoin}). Consequently, the arguments have to be repeated for several shifted copies of the fine partition (see Lemma \ref{rig-lemma: weaklocA2}). Moreover, the projections $P_Q$ onto the the space of infinitesimal rigid motions (see Theorem \ref{rig-theo: korn} below) have to be combined with the rotations $R_Q$ in a suitable way to obtain appropriate rigid motions, which do not vary too much on adjacent squares (see Lemma \ref{rig-lemma: weaklocA}).

Having an approximation of the deformation by piecewise rigid motions defined on squares with diameter $k$, we then are able to modify the function such that the minimal distance $\tilde{s}$ of two cracks of the new configuration satisfies $\tilde{s} \sim k$ (see Lemma \ref{rig-th: global estimate}). Now we can repeat the above procedure for some larger $\tilde{k}$ such that $\eps/\tilde{k} \le (C(\tilde{k}/\tilde{s}))^{-2}$ is guaranteed and we can repeat the arguments in \eqref{rig-eq: over:rig2}. 

The strategy is to end up with $k \approx \rho$ after a finite number of iterations. As the number of iteration steps is not bounded but grows logarithmically with $\frac{1}{\eps}$ we have to assure that in each step  the surface and the elastic energy do not increase too much. The crucial point is that during the iteration process the coarseness of the partition $k$ grows much faster than the stored elastic energy $\eps$ such that the argument in \eqref{rig-eq: over:rig2} may be repeated. The details are given in Theorem \ref{rig-thm: V2}.  Having an estimate for $k \approx \rho$ it is then not hard to establish the desired result up to a small exceptional set (see Theorem  \ref{rig-thm: V1}).

Clearly, we cannot assume that initially $s \ge \eps$. In this case the argument in \eqref{rig-eq: over:rig2} can typically not be applied. As a remedy we do not employ the geometric rigidity result directly but first approximate the deformation in each square by an $H^1$-function, where the distance can be  measured by the curl of $\nabla y$. (See Theorem \ref{rig-th: cgp} below which was one of the essential ingredients to prove the qualitative result in \cite{Chambolle-Giacomini-Ponsiglione:2007}.) We address this problem in Lemma \ref{rig-lemma: local estimate2} and subsequently we show that we may modify the configuration such that $\tilde{s} \ge \eps$ (see Theorem \ref{rig-thm: V2.5}).

Finally, by a density argument we can approximate each SBV function by a configuration where the jump set is contained in a finite number of rectangle boundaries (see proof of Theorem \ref{rig-th: rigidity2}). Observe that standard density results as \cite{Cortesani:1997} cannot be applied directly in our framework since in general an  $L^\infty$ bound for the derivative is not preserved. The problem can be circumvented by using a different approximation introduced in \cite{Chambolle:2004} at  the cost of a non exact approximation of the jump set, which suffices for our purposes.

The rigidity result, which we then have established, only holds up to a small exceptional set as  due to the modification of the jump set the deformation might not be defined in the interior of certain rectangles. We emphasize that such an estimate is not enough to obtain good compactness and convergence results, in particular for the convergence of the  surface energy further difficulties arise. Therefore, we eventually have to construct a suitable extension to the whole domain. A major challenge is to determine the surface energy correctly, at least for the relaxed functional \eqref{rig-eq: Griffith en2}. This problem is addressed in Section \ref{rig-sec: sub, proof-main}. 

For small cracks a good extension is already provided by the Korn-Poincar\'e inequality \cite{Friedrich:15-1} which is based on the derivation of  a suitable modification for which jump heights can be controlled. Near large cracks we define the extension as a piecewise constant rigid motion such that the jump heights on the new jump sets are sufficiently small (see the proof of Theorem \ref{rig-th: rigidity}). Consequently, the length of these jumps may possibly be  much larger than ${\cal H}^1(J_y)$, but due to the small jump height their contribution to \eqref{rig-eq: Griffith en2} is considerably small. Finally, for the large cracks in the domain, in particular for the boundary $\bigcup_j \partial P_j$ of the partition $(P_j)_j$, we have to construct an appropriate jump set consisting of Jordan curves which provides the correct crack energy up to a small error (see Lemma \ref{rig-lemma: jordan}).


\section{Preliminaries}\label{rig-sec: pre}

In this preparatory section  we recall first  the definition and basic properties of functions of bounded variation. Then we introduce the notion of \emph{boundary components} and present the Korn-Poincar\'e inequality established in \cite{Friedrich:15-1}. Finally, we recall the geometric rigidity result in nonlinear elasticity and  carefully estimate the involved constant pertaining to its dependence on the shape of the domain.

\subsection{Special functions of bounded variation}\label{rig-sec: sub, bd}

In this section we collect the definitions of SBV and SBD functions. Let $\Omega \subset \R^d$ open, bounded with Lipschitz boundary. Recall that the space $SBV(\Omega, \R^d)$, abbreviated as $SBV(\Omega)$ hereafter,  of \emph{special functions of bounded variation} consists of functions $y \in L^1(\Omega, \R^d)$ whose distributional derivative $Dy$ is a finite Radon measure, which splits into an absolutely continuous part with density $\nabla y$ with respect to Lebesgue measure and a singular part $D^j y$ whose Cantor part vanishes and thus is of the form 
$$ D^j y = [y] \otimes \xi_y {\cal H}^{d-1} \lfloor J_y, $$
where ${\cal H}^{d-1}$ denotes the $(d-1)$-dimensional Hausdorff measure, $J_y$ (the `crack path') is an ${\cal H}^{d-1}$-rectifiable set in $\Omega$, $\xi_y$ is a normal of $J_y$ and $[y] = y^+ - y^-$ (the `crack opening') with $y^{\pm}$ being the one-sided limits of $y$ at $J_y$. If in addition $\nabla y \in L^2(\Omega)$ and ${\cal H}^{d-1}(J_y) < \infty$, we write $y \in SBV^2(\Omega)$. See \cite{Ambrosio-Fusco-Pallara:2000} for the basic properties of this function space. 

Likewise, we say that a function $y \in L^1(\Omega, \R^d)$ is a \emph{special  function of bounded deformation} if the symmetrized distributional derivative $Eu := \frac{(Dy)^T + Dy}{2}$ is a finite $R^{d \times d}_{\rm sym}$-valued Radon measure with vanishing Cantor part. It can be decomposed as 
\begin{align}\label{rig-eq: symmeas}
 Ey = e(\nabla y) {\cal L}^d  + E^j y = e(\nabla y) {\cal L}^d + [y] \odot \xi_y {\cal H}^{d-1}|_{J_y},
 \end{align}
where $e(\nabla y)$ is the absolutely continuous part of $Ey$ with respect to the Lebesgue measure ${\cal L}^d$, $[y]$, $\xi_y$, $J_y$ as before and $a \odot b = \frac{1}{2}(a \otimes b + b \otimes a)$. For basic properties of this function space we refer to \cite{Ambrosio-Coscia-Dal Maso:1997,  Bellettini-Coscia-DalMaso:98}.

The general idea in our analysis will be to establish Theorem \ref{rig-th: rigidity} for a dense subset of SBV for which we can suppose much more regularity of the jump set. For density results in the spaces SBV and SBD we refer to \cite{Cortesani:1997, Cortesani-Toader:1999} and \cite{Chambolle:2004}, respectively. In our framework we cannot use these results directly but have to derive a slightly different variant of \cite{Cortesani:1997} in order to preserve an $L^\infty$-bound for the derivative (see the proof of Theorem \ref{rig-th: rigidity2}). 

Moreover, we recall the property that the distance of an SBV function to Sobolev functions can be measured by the distribution $\text{curl}\,\nabla y$ (see \cite[Proposition 5.1]{Chambolle-Giacomini-Ponsiglione:2007}). 

\begin{theorem}\label{rig-th: cgp}
Let $Q = (0,1)^d$. Let $y \in SBV_\infty(Q) := \lbrace y \in SBV(Q,\R^d): \Vert\nabla y \Vert_\infty < \infty, \ \ {\cal H}^{d-1}(J_y) < \infty\rbrace$. Then $\mu_y:=\text{curl}\, \nabla y$ is a measure concentrated on $J_y$ such that
$$|\mu_y| \le C\Vert\nabla y\Vert_\infty {\cal H}^{d-1}|_{J_y}.$$
Moreover, for $p < \frac{d}{d-1}$ there is a constant $C=C(p)>0$ such that for all $y \in SBV_\infty(Q)$ there is a function $u \in H^1(Q,\R^d)$ such that
$$\Vert \nabla u - \nabla y\Vert_{L^p(Q)} \le C |\mu_y|(Q) \le  C\Vert\nabla y\Vert_\infty {\cal H}^{d-1}(J_y).$$
\end{theorem}

\subsection{Boundary components}\label{rig-sec: modifica}

Using a density result alluded to above it will suffice to prove the main result for configurations where the jump set is contained in the boundary of squares. In this section we recall the necessary notation and definitions for boundary components introduced in \cite{Friedrich:15-1}.

For $s>0$ we partition $\R^2$ up to a set of measure zero into squares $Q^s(p) = p + s(-1,1)^2$ for $p \in I^s := s(1, 1) + 2s\Z^2$.  Let 
\begin{align}\label{rig-eq: calV-s def}
{\cal U}^s := \Big\{ U \subset \R^2: U = \Big(\bigcup\nolimits_{p \in I} \overline{Q^s(p)} \Big)^\circ: \  \ I \subset I^s \Big\}.
\end{align}
Here the superscript $\circ$ denotes the interior of a set. Let $\mu>0$. We will concern ourselves with subsets  $V  \subset Q_\mu:=(-\mu,\mu)^2$ of the form
\begin{align}\label{rig-eq: calV-s def2}
{\cal V}^s := \lbrace V \subset Q_\mu: V = Q_\mu  \setminus \, \bigcup\nolimits^m_{i=1} X_i, \ \ X_i \in {\cal U}^s, \ X_i  \text{ pairwise disjoint} \rbrace
\end{align}
for $s >0$. Note that each set in $V \in {\cal V}^s$ coincides with a set $U \in {\cal U}^s$ up to subtracting a set of zero Lebesgue measure, i.e. $U \subset V$, ${\cal L}^2(V \setminus U) = 0$. The essential difference of $V$ and the corresponding $U$ concerns the connected components of the complements $Q_\mu \setminus V$ and $Q_\mu \setminus U$.  Observe that one may have $Q_\mu  \setminus \, \bigcup\nolimits^{m}_{i=1} X_i = Q_\mu  \setminus \, \bigcup\nolimits^{\hat{m}}_{i=1} \hat{X}_i$ with $(X_1,\ldots,X_{m}) \neq (\hat{X}_1,\ldots,\hat{X}_{\hat{m}})$, e.g. by combination of different sets. In such a case we will regard $V_1 = Q_\mu  \setminus \, \bigcup\nolimits^{m}_{i=1} X_i$ and $V_2 = Q_\mu  \setminus \, \bigcup\nolimits^{\hat{m}}_{i=1} \hat{X}_i$ as different elements of ${\cal V}^s$. For this and the following sections we will always tacitly assume that all considered sets are elements of ${\cal V}^s$ for some small, fixed $s>0$.

Let  $W \in {\cal V}^s$ and arrange the components $X_1, \ldots, X_m$ of the complement such that $\partial X_i \subset Q_\mu$ for $1 \le i \le n$ and $\partial X_i \cap  \partial Q_\mu \neq \emptyset$ otherwise. Define $\Gamma_i(W) = \partial X_i$ for $i=1,\ldots,n$. In the following we will often refer to these sets as \emph{boundary components}. Note that $\bigcup^n_{i=1} \Gamma_i(W)$ might not cover $\partial W \cap Q_\mu$ completely if $n < m$. We frequently drop the subscript and write $\Gamma(W)$ or just $\Gamma$ if no confusion arises. Observe that in the definition we do not require that boundary components are connected. Therefore, we additionally introduce the subset ${\cal V}^s_{\rm con} \subset {\cal V}^s$ consisting of the sets where all  $\overline{X_1},\ldots, \overline{X_n}$ are connected.

Beside the Hausdorff-measure $\vert \Gamma \vert_{\cal H} = {\cal H}^1(\Gamma)$ (we will use both notations) we define the `diameter' of a boundary component by 
$$\vert \Gamma\vert_\infty:= \sqrt{|\pi_1 \Gamma|^2 + |\pi_2 \Gamma|^2},$$
 where $\pi_1$, $\pi_2$ denote the orthogonal projections onto the coordinate axes. We recall that many arguments in the proof of the Korn-Poincar\'e inequality in \cite{Friedrich:15-1} relied on the fact that due to the strict convexity of $|\cdot|_\infty$  it is often energetically favorable if different components are combined to a larger one.

Note that by definition of ${\cal V}^s$ (in contrast to the definition of ${\cal U}^s$) two components in $(\Gamma_i)_i$ might not be disjoint. Therefore, we choose an (arbitrary) order $(\Gamma_i)^n_{i=1} = (\Gamma_i(W))^n_{i=1}$ of the boundary components of $W$, introduce 
\begin{align}\label{rig-eq: Xdef}
\Theta_i = \Theta_i(W) = \Gamma_i \setminus \bigcup\nolimits_{j<i} \Gamma_j
\end{align}
 for $i=1,\ldots,n$ and observe that the boundary components $(\Theta_i)_i$ are pairwise disjoint. With a slight abuse of notation we define
$$\vert \Theta_i \vert_\infty = \vert\Gamma_i \vert_\infty.$$
Again we will often drop the subscript if we consider a fixed boundary component. We now introduce a convex combination of $\vert \cdot \vert_\infty$ and $\vert \cdot \vert_{\cal H}$. For an $h_*>0$ to be specified below we set
\begin{align}\label{rig-eq: h*}
\vert \Theta \vert_* =  h_* \vert \Theta \vert_{\cal H} + (1-h_*) \vert \Theta \vert_\infty.
\end{align}
For sets $W \in {\cal V}^s$ we then define 
\begin{align}\label{rig-eq: h**}
\Vert W \Vert_Z =  \sum\nolimits^n_{j=1} \vert \Theta_j(W) \vert_Z 
\end{align}
for $Z={\cal H}, \infty,*$. Note that $\Vert W\Vert_\infty, \Vert W\Vert_{\cal H}$ and thus also $\Vert W\Vert_*$ are independent of the specific order which we have chosen in \eqref{rig-eq: Xdef}. Indeed, for $\Vert W\Vert_\infty$ this is clear as $\vert \Theta_i \vert_\infty = \vert\Gamma_i \vert_\infty$, for $\Vert W\Vert_{\cal H}$ it follows from the fact that $\Vert W\Vert_{\cal H} = {\cal H}^1(\bigcup^n_{i=1} \Gamma_i)$.

From \cite{Friedrich:15-1} we recall some elementary properties of $\vert \cdot \vert_*$ which will be exploited frequently in the following.  

\begin{lemma}\label{rig-lemma: infty}
Let $W\subset Q_\mu$. Let $\Gamma=\Gamma(W)$ be a  boundary component with $\Gamma = \partial X$ and let $\Theta \subset \Gamma$ be the corresponding set defined in \eqref{rig-eq: Xdef}. Moreover,  let  $V \in {\cal U}^s$ be a rectangle with $\overline{V} \cap \overline{X} \neq \emptyset$. Suppose that $h_*$ is sufficiently small. Then
\begin{itemize} 
\item[(i)] $|\Gamma|_* \ge |\partial R(\Gamma)|_*$  if $\Gamma$ is connected, where $R(\Gamma)$ denotes the smallest (closed) rectangle such that $\Gamma \subset R(\Gamma)$,
\item[(ii)] $\vert \Theta \vert_* = \vert \Gamma \vert_* \Leftrightarrow \vert \Theta \vert_{\cal H} = \vert \Gamma \vert_{\cal H}$,
\item[(iii)] $|\partial (X \setminus \overline{V})|_\infty \le \vert \Theta \vert_\infty$ and $|\Theta \setminus \overline{V}|_{\cal H} \le \vert \Theta \vert_{\cal H}$,
\item[(iv)] $ |\partial (V \cup X)|_* \le |\partial V|_* + |\Gamma|_* $,
\item[(v)] $ \frac{1}{\sqrt{2}}\vert \partial R \vert_{\cal H} \le 2\vert \partial R \vert_\infty \le \vert \partial R \vert_{\cal H}$ if $R \in {\cal U}^s$ is are rectangle.
\end{itemize}
\end{lemma}

As a further preparation, we define $H(W)\supset W \in {\cal V}^s$ as the `variant of $W$ without holes' by
\begin{align}\label{rig-eq: no holes}
H(W) = W \cup \bigcup\nolimits^n_{j=1} X_j.
\end{align} 
Additionally, for $\lambda > 0$ we define $H^{\lambda}(W)\supset W$ as the `variant of $W$ without holes of size smaller than $\lambda$': 
We arrange the sets $(\Gamma_j)_{ j=1,\ldots,n}$ in the way that $\vert \Gamma_j \vert_\infty \le \lambda$ for $j \ge l_\lambda$ and $\vert \Gamma_j \vert_\infty > \lambda$ for $j < l_\lambda$. Define 
\begin{align}\label{rig-eq: no holes2}
H^{\lambda}(W) =  W \cup \bigcup\nolimits^{n}_{  j=l_\lambda} X_j.
\end{align}

\subsection{A Korn-Poincar\'e inequality}\label{rig-sec: subsub,  mainproof}

We start this section with the formulation of the classical Korn-Poincar\'e inequality in BD (see \cite{Kohn:82,Temam:85}). 

\begin{theorem}\label{rig-theo: korn}
Let $\Omega \subset \R^d$ bounded,  connected with Lipschitz boundary and let $P: L^2(\Omega,\R^d) \to  L^2(\Omega,\R^d)$ be a linear projection onto the space of infinitesimal rigid motions. Then there is a constant $C>0$,  which is invariant under rescaling of the domain, such that for all $u \in BD(\Omega,\R^d)$
$$\Vert u - Pu \Vert_{L^{\frac{d}{d-1}}(\Omega)} \le C |Eu|(\Omega),$$
where $Eu = \frac{Du^T + Du}{2}$ is the symmetrized distributional derivative.
\end{theorem}

There is also a corresponding trace estimate.

\begin{theorem}\label{rig-th: tracsbv}
Let $\Omega \subset \R^2$ bounded,  connected with Lipschitz boundary. There exists a constant $C > 0$ such
that the trace mapping $\gamma: BD(\Omega,\R^2) \to  L^1(\partial \Omega, \R^2)$ is well defined and
satisfies the estimate
$$\Vert \gamma u \Vert_{L^1(\partial \Omega)}\le C\big(\Vert u\Vert_{L^1(\Omega)} + |Eu|(\Omega) \big)$$
for each $u \in BD(\Omega,\R^2).$
\end{theorem}

It first appears that this inequality is not adapted for linearized Griffith energies of the form \eqref{rig-eq: Fenergy} (or their nonlinear counterparts \eqref{rig-eq: Griffith en}) as in $|Eu|(\Omega)$ the jump height is involved and in \eqref{rig-eq: Fenergy} we only have control
over the size of the crack. However, in \cite{Friedrich:15-1} we have shown that one can indeed find bounds on the jump heights after a suitable modification of the jump set and the displacement field. Before we can recall the results obtained in \cite{Friedrich:15-1}, we have to introduce a further notation: We fix a sufficiently large universal constant $c$ and let ${\cal W}^s \subset {\cal V}^s$ be the subset consisting of the sets, where  for a specific ordering of the boundary components $(\Gamma_l)^n_{l=1}$  we find for all components $\Gamma_l$ a corresponding rectangle $R_l =  R(\Gamma_l) \in {\cal U}^s$  such that
\begin{align}\label{rig-eq: newnew}
(i) \ \, |\Gamma_l|_\infty \le |\partial R_l|_\infty \le c|\Gamma_l|_\infty, \ \  (ii) \ \,  \vert \Theta_l\vert_{\cal H} \le \vert \partial R_l \vert_{\cal H}, \ \  (iii) \ \, \vert \partial R_l \vert_* \le c \vert \Theta_l\vert_*.
\end{align}
In particular, the diameter of $\Gamma_l$ and the corresponding rectangle $R_l$ are  comparable. (Note that in \cite[Section 5]{Friedrich:15-1} we have defined the set ${\cal W}^s$ in a slightly different way. See also (3.5) and (3.6) in \cite{Friedrich:15-1}.) For given $\bar{\tau}>0$ and a rectangle $R_l \in {\cal U}^s$ we define $\tau_l = \bar{\tau} |\partial R_l|_\infty$ and let $N^{\tau_l}(\partial R_l) \in {\cal U}^s$ be the largest set in ${\cal U}^s$ with $N^{\tau_l}(\partial R_l) \subset \lbrace x \in \R^2 \setminus \overline{R_l}: \dist_\infty(x, \partial R_l) \le \tau_l\rbrace$, where $\dist_\infty(x,A) := \inf_{y \in A} \max_{i=1,2} |(x - y) \cdot \e_i|$ for $A \subset \R^2$, $x \in  \R^2$.  We can now formulate \cite[Theorem 5.2]{Friedrich:15-1} as follows.

\begin{theorem}\label{rig-th: derive prop}
Let $\eps>0$ and $h_* \ge \sigma >0$ sufficiently small. Let $C_1=C_1(\sigma,h_*) \ge 1$ large, $0 < C_2=C_2(\sigma,h_*) <1$ small enough, and $\bar{\tau}>0$ such that $C_2 \ll \bar{\tau} \ll 1$. Moreover, let $c>0$ be a universal constant. Then for all $W \in {\cal V}^{s}_{\rm con}$  and $u \in H^1(W)$ there is a set $U \in {\cal W}^{C_2s}$  with $|U \setminus W| = 0$   and an extension $\bar{u}$ in SBV defined by
\begin{align}\label{rig-eq: extend def_new}
\bar{u}(x) = \begin{cases} A_l\,x  +c_l & x \in X_l \ \ \ \text{for all } \Gamma_l(U) \text{ with } N^{\tau_l}(\partial R_l) \subset H(U), \\
u(x) & \text{else,} \end{cases}
\end{align} 
such that for all $\Gamma_l(U)$  with $N^{\tau_l}(\partial R_l) \subset H(U)$
\begin{align}\label{rig-eq: D1} 
\int\nolimits_{\Theta_l(U)} |[\bar{u}](x))|^2 \,d{\cal H}^1(x) \le C_1  \eps  \vert \Theta_l(U) \vert^2_*.
\end{align}
Moreover, one has   $|W \setminus U| \le c\Vert U \Vert^2_\infty$ and
\begin{align*}
\eps \Vert U \Vert_* + \Vert e(\nabla u)\Vert^2_{L^2(U)} \le (1+ \sigma)\big(\eps \Vert W \Vert_* + \Vert e(\nabla u)\Vert^2_{L^2(W)}\big).
\end{align*}
 \end{theorem}

  \begin{rem}\label{rig_rem: connect}
 {\normalfont
 
(i) During the modification process in Theorem \ref{rig-th: derive prop} the components  $X_{n+1}(W), \ldots, X_{m}(W)$  at the boundary of $Q_\mu$ might be changed and the corresponding components of $U$ are given by $X_j(U) = X_j(W) \setminus \overline{H(U)}$ for $j=n+1,\ldots,m$. In particular, one has $|\partial X_j(U)|_* \le |\partial X_j(W)|_*$   arguing as in Lemma \ref{rig-lemma: infty}.

(ii) Observe that $U \notin {\cal V}^s_{\rm con}$ is possible as components can  be separated by other components in the proof of Theorem \ref{rig-th: derive prop}. However, we can obtain a set $U' \subset U$ with $\Vert U' \Vert_* \le \Vert U \Vert_*$ and $|U\setminus U'| \le C \Vert U' \Vert_\infty^2 \le C\mu \Vert U' \Vert_\infty$ such that all components of $U'$ are pairwise disjoint and rectangular and thus particularly connected.  Moreover, for each $\Gamma(U)$ the corresponding rectangle $R(U)$ given by \eqref{rig-eq: newnew} is contained in a component of $U'$. (Namely in the same component as $\Gamma(U)$.)

 }
 \end{rem}

Recall \eqref{rig-eq: symmeas} and define ${\cal E}(V) = \int_V |e(u)| + |D^j u|(V)$. Observe that ${\cal E}(V)$ differs from $|Eu|(V)$ as we consider the measure $D^j u$ instead of $E^j u$. We then obtain the following corollary  (cf.\ \cite[Corollary 5.7]{Friedrich:15-1}).

\begin{corollary}\label{rig-cor: kornpoin}
Let $\eps, \mu, h_* >0$. Let $U\subset Q_\mu=(-\mu,\mu)^2$, $U \in {\cal W}^{C_2 s}$ and $u \in H^1(U)$. Assume there is a square $\tilde{Q} = (-\tilde{\mu},\tilde{\mu})^2 \subset Q_\mu$ such that \eqref{rig-eq: D1} is satisfied for all components $\Theta_l(U)$ having nonempty intersection with $\tilde{Q}$, where $\bar{u}$ is the extension of $u$ defined in \eqref{rig-eq: extend def_new}. Then there is a universal constant $C$ such that

\begin{align*} 
|E\bar{u}|(\tilde{Q})^2  \le ({\cal E}(\tilde{Q}))^2\le C \tilde{\mu}^2 \Vert e(\nabla u) \Vert^2_{L^2(U\cap \tilde{Q})} + CC_1\mu\eps|\partial U \cap \tilde{Q}|_{\cal H} |\partial U \cap Q_\mu|_{\cal H},
\end{align*}
where $C_1$ is the constant in Theorem \ref{rig-th: derive prop}.
\end{corollary}

Now observe that by combination of Theorem  \ref{rig-th: derive prop}, Corollary \ref{rig-cor: kornpoin} and Theorem \ref{rig-theo: korn} one may estimate the distance of $u$ from an infinitesimal rigid motion. We will exploit this property in Section \ref{rig-sec: subsub,  est-h1}.  \cite[Lemma 6.7]{Friedrich:15-1} provides the following estimate for the skew symmetric matrices  involved in \eqref{rig-eq: extend def_new}. 

\begin{lemma} \label{rig-lemma: A neigh}
Let be given the situation of Theorem \ref{rig-th: derive prop}  for a function $u \in H^1(W)$ and define $y = \bar{R}\, (\id + u)$, where $\id$ denotes the identity function and $\bar{R} \in SO(2)$. Let  $V \subset Q_\mu$ be a rectangle and  let ${\cal F}(V)$  be the boundary components $(\Gamma_l)_l =  (\Gamma_l(U))_l$ satisfying $N^{\tau_l}(\partial R_l) \subset V$ and \eqref{rig-eq: D1}. Then there is a $C_3 = C_3(\sigma,h_*)$ such that   
$$\sum\nolimits_{\Gamma_l \in {\cal F}(V)}  |X_l|_\infty^2 |A_l|^p \le C_3 \big( \Vert \nabla y - \bar{R}\Vert^p_{L^p(V  \cap W)} +  (\eps s^{-1})^{\frac{p}{2}-1} \eps|\partial U \cap V|_{\cal H} \big)$$
for $p=2,4$, where  $X_l \subset Q_\mu$, $A_l \in \R^{2 \times 2}_{\rm skew}$ is given in  \eqref{rig-eq: extend def_new}.
\end{lemma}
 
 We close this section with a short remark about the constants involved in the above estimates.

\begin{rem}\label{rig-rem: z}
 {\normalfont
(i) The constants $C_i=C_i(\sigma,h_*)$, $i=1,2,3$, have polynomial growth in $\sigma$: We find $z \in \N$ large enough such that $C_1(\sigma,h_*), C_3(\sigma,h_*) \le C(h_*) \sigma^{-z}$  and $C_2(\sigma,h_*) \ge C(h_*) \sigma^{z}$.

(ii) The constant $C_2(\sigma,h_*)$ can be chosen small with respect to $\sigma$ (see (5.12) in \cite{Friedrich:15-1}). In particular, we can assume $C_2(\sigma,h_*) \ll \sigma$ as well as $\bar{C} C_2(\sigma,h_*)  \le \sigma$ for constants $\bar{C} = \bar{C}(h_*)$.

(iii) We find a constant $\bar{C} = \bar{C}(h_*)$ such that $\bar{\tau} \le \bar{C}C_2$ (cf. (5.2) in \cite{Friedrich:15-1}).

(iv) If we apply Theorem \ref{rig-th: derive prop} on sets $W \in {\cal V}^{\bar{s}}_{\rm con}$ for some $\bar{s}\ll s$, where  the length of all boundary components of $W$ is bounded from below by $s$, we still obtain $U \in {\cal V}^{C_2s}$.
 }
 \end{rem}

\subsection{Geometric rigidity in nonlinear elasticity}\label{rig-sec: sub, georid}

The following geometric rigidity result in nonlinear elasticity proved by Friesecke, James and M\"uller (see \cite{FrieseckeJamesMueller:02}) is one of the starting points for our analysis.

\begin{theorem}\label{rig-th: geo rig}
Let $\Omega \subset \R^d$ a (connected) Lipschitz domain and $1 < p < \infty$. Then there exists a constant $C = C(\Omega,p)$ such that for any $y \in W^{1,p}(\Omega,\R^d)$ there is a rotation $R \in SO(d)$ such that
\begin{align*}
\left\|\nabla y - R\right\|_{L^p(\Omega)} \leq C \left\|\dist(\nabla y, SO(d))\right\|_{L^p(\Omega)}. 
\end{align*}
\end{theorem}

One ingredient in the proof  is the following decomposition into a harmonic and a rest part.

\begin{theorem}\label{rig-thm: harmonic}
Let $\Omega\subset \R^2$ open and $1 < p <\infty$. There is a constant $C=C(p)$ such that all $y \in W^{1,p}(\Omega,\R^2)$  can be split into $y = w + z$, where $w$ is a harmonic function and $z$ satisfies
$$\Vert \nabla y - \nabla w \Vert_{L^p(\Omega)} = \Vert \nabla z \Vert_{L^p(\Omega)} \le C \Vert \dist(\nabla y,SO(2))\Vert_{L^p(\Omega)}.$$
\end{theorem}
Note that the constant $C$ is independent of the domain $\Omega$. In higher dimensions one additional needs $\Vert \nabla y\Vert_\infty \le M$ for $M >0$.

\Proof Following the singular-integral estimates in \cite[Section 2.4]{ContiSchweizer:06} we find $\Vert \nabla z\Vert_{L^p(\Omega)} \le c \Vert \text{cof} \nabla y - \nabla y \Vert_{L^p(\Omega)}.$ The assertion follows from the fact that $|\text{cof}A - A|^p \le C_p \dist^p(A,SO(2))$ for all $A \in \R^{2 \times 2}$ (see also (3.11) in \cite{FrieseckeJamesMueller:02}). \eop

For sets which are related through bi-Lipschitzian homeomorphisms with Lipschitz constants of both
the homeomorphism itself and its inverse uniformly bounded the constant in Theorem \ref{rig-th: geo rig}  can be chosen independently of these sets, see e.g. \cite{FrieseckeJamesMueller:02}.

\subsection{Geometric rigidity: Dependence on the set shape}\label{rig-sec: sub, pre} 

In general, the constant of the inequality stated in Section \ref{rig-sec: sub, georid} depends crucially on the set shape. This will be discussed in detail in this section. As an introductory example we consider the deflection of a thin elastic beam.

\begin{example}\label{rig-example: FMJ rectangle} \upshape Let $U = (0,1) \times (0,\delta)$ and let $y: U \to \R^2$ be given by $y(x_1,x_2) = (x_2+1) (\sin(x_1), \cos(x_1))$. Then 
$$\nabla y (x_1,x_2) = \begin{pmatrix}  (x_2+1) \cos(x_1) & \sin(x_1) \\ -(x_2+1) \sin(x_1) & \cos(x_1) \end{pmatrix}$$ 
and therefore $\dist^2(\nabla y, SO(2)) = |\sqrt{\nabla y^T \nabla y} - \Id|^2 = x_2^2$, i.e. 
$$\Vert\dist(\nabla y, SO(2))\Vert^2_{L^2(U)} = \tfrac{1}{3}\delta^3.$$
Let $R_{\phi} \in SO(2)$, $R_{\phi} = \begin{pmatrix}  \cos\phi & \sin\phi \\ -\sin\phi & \cos\phi \end{pmatrix}$ for $\phi \in [0,2\pi]$. Then $|\nabla y(x) - R|^2 \geq |\sin(x_1) - \sin\phi|^2 + |\cos(x_1) - \cos\phi|^2$. It is not hard to see that it exists a $C>0$ such that $\int^{1}_{0} |\nabla y(x) - R|^2 \, dx_1 \ge C$ for all $\phi \in [0,2\pi]$ and $x_2 \in (0,\delta)$. We conclude that 
$$\Vert \nabla y - R\Vert^2_{L^2(U)} \ge C\delta \ge \frac{C}{\delta^2}\Vert\dist(\nabla y, SO(2))\Vert^2_{L^2(U)}$$
for all $R \in SO(2)$.
A similar argument shows 
$$\Vert y - (R\,\cdot + c) \Vert^2_{L^2(U)} \ge C\delta \ge \frac{C}{\delta^2}\Vert\dist(\nabla y, SO(2))\Vert^2_{L^2(U)}$$
for all $R \in SO(2)$ and  $c \in \R^2$.
\end{example}
Similar examples can be constructed in the linearized framework for the Korn-Poincar\'e inequality given in Theorem \ref{rig-theo: korn}. As a direct consequence we get that the estimate \eqref{rig-eq: main properties2} might be wrong without allowing for a small modification of the deformation.

\begin{example}\label{rig-example: FMJ rectangle2} \upshape
Let $\eps> 0$. Assume without restriction that the set $U = (0,1) \times (0,\eps^{\frac{1}{3}})$ considered above satisfies $ \overline{U} \subset \Omega$. Define $y : \Omega \to \R^2$ by $y(x) = \id + \e_2$ for $x \in \Omega \setminus U$ and $y(x)  = (x_2+1) (\sin(x_1), \cos(x_1))$ for $x \in U$. Then $y \in SBV^2(\Omega)$ with $J_y = (0,1) \times \lbrace 0, \eps^{\frac{1}{3}}\rbrace \cup \lbrace 1\rbrace \times (0, \eps^{\frac{1}{3}})$ and $\Vert \dist(\nabla y, SO(2))\Vert^2_{L^2(\Omega)} = \frac{\eps}{3}$. However, for all $R \in SO(2)$ and $c \in \R^2$ we have
$$\Vert \nabla y - R\Vert^2_{L^2(\Omega)} \ge C\eps^{\frac{1}{3}}, \ \ \ \ \Vert y - (R\,\cdot + c)\Vert^2_{L^2(\Omega)} \ge C\eps^{\frac{1}{3}}.$$
Although omitted here, a similar estimate can be derived for the symmetric part of the gradient. 
\end{example}

Recall the definition of ${\cal U}^s$ in \eqref{rig-eq: calV-s def}. In order to quantify how the constant in Theorem \ref{rig-th: geo rig}  depends on the set shape we will estimate the variation from a square $Q^s(a)$ to a neighboring square $Q^s(b)$, $b= a + 2s\nu$ for $\nu = \pm \e_i$, $i=1,2$ proceeding similarly as in \cite{FrieseckeJamesMueller:02}. Consider $y \in H^1(U)$ with $U \in {\cal U}^s$. On a square $Q^s(p) \subset U$ and for subsets $V \subset U$, $V \in {\cal U}^{s}$ we define for shorthand (we drop the integration variable if no confusion arises)
\begin{align*}
\gamma(p) = \int_{Q^s(p)} \dist^2(\nabla y,SO(2)), \ \ \ \gamma(V) = \sum\nolimits_{p \in I^s(V)} \gamma(p), 
\end{align*}
where $I^s (V) := \lbrace p \in I^s: Q^s(p) \subset V\rbrace.$ Applying Theorem \ref{rig-th: geo rig} we obtain $R(a), R(b) \in SO(2)$ such that
\begin{align}\label{rig-eq: R rig}
\int_{Q^s(p)} |\nabla y - R(p)|^2 \le C \gamma(p) \ \ \ \text{ for } p = a,b.
\end{align}
Likewise on the rectangle $Q^s(a,b) := (\overline{Q^s(a)} \cup \overline{Q^s(b)})^\circ$ we obtain $R(a,b) \in SO(2)$ such that
$$\int_{Q^s(a,b)} |\nabla y - R(a,b)|^2 \, dx \le C \int_{Q^s(a,b)} \dist^2(\nabla y,SO(2)) \le C(\gamma(a) + \gamma(b)).$$
Combining these estimates we see $|Q^s(p)| |R(p) - R(a,b)|^2 \le C(\gamma(a) + \gamma(b))$ for $p=a,b$ and therefore
\begin{align}\label{rig-eq: R difference}
s^2 |R(a) - R(b)|^2 \le  C(\gamma(a) + \gamma(b)).
\end{align}
More general, we consider a difference quotient with two arbitrary points $a,b \in I^s(U)$. We assume that there is a path $\xi=(\xi_0,\ldots,\xi_m)$ such that
\begin{align}\label{rig-eq: path}
\begin{split}
&\xi_1 = a, \ \  \xi_{m} = b,  \\ 
&\xi_j - \xi_{j-1} = \pm 2s \e_i \text{ for some } i=1,2, \ \  \forall j=2, \ldots,m.
\end{split}
\end{align}
Then iteratively applying the above estimate \eqref{rig-eq: R difference} we obtain
\begin{align}\label{rig-eq: R diff}
s^2|R(a) - R(b)|^2 \, dx \le Cm \sum\nolimits^m_{j=1} \gamma(\xi_j).
\end{align}
We now state a first weak rigidity result.

\begin{lemma}\label{rig-lemma: weak rig}
Let $\mu,s > 0$ such that $l:=\mu s^{-1} \in \N$. Then there is a constant $C > 0$ independent of $\mu$, $s$ such that for all connected sets $U \in {\cal U}^s$,  $U \subset (-\mu,\mu)^2$, the following holds:  For all $y \in H^{1}(U)$ there is a rotation $R \in SO(2)$ such that
$$\int_U |\nabla y- R|^2 \le C(s^{-2}|U|)^2 \int_U \dist^2(\nabla y, SO(2)) \le C l^4 \int_U \dist^2(\nabla y, SO(2)).$$
\end{lemma}

\Proof The second inequality is obvious as $|U| \le 4\mu^2$. To see the first inequality we fix $p_0 \in I^s(U)$ and consider an arbitrary $p \in I^s(U)$. As $U$ is connected there is a path $\xi=(\xi_1=p_0,\ldots,\xi_m=p)$ with $m \le |U|(2s)^{-2} $. We first apply \eqref{rig-eq: R rig} on each square and then by \eqref{rig-eq: R diff} we obtain
$$\int_{Q^s(p)} |R(p) - R(p_0)|^2 \le C |U|s^{-2} \sum\nolimits^m_{j=1} \gamma(\xi_j) \le C |U|s^{-2} \gamma(U).$$
Then setting $R = R(p_0)$ and summing over all $p \in I^s(U)$ we derive
\begin{align*}
\int_U |\nabla y - R|^2  & \le C\sum\nolimits_{p \in I^s(U)} \int_{Q^s(p)}\Big(|\nabla y - R(p)|^2 + |R(p) - R(p_0)|^2 \Big)  \\
& \le C \sum\nolimits_{p \in I^s(U)} (\gamma(p) + |U|s^{-2} \gamma(U)) \le C \# I^s(U) \, |U|s^{-2} \gamma(U) \\
& \le C (|U|s^{-2})^2 \gamma(U).
\end{align*}
 \eop

\begin{rem}\label{rig-rem: 1} 
{\normalfont
\begin{itemize}
\item[(i)] Let $U=(0,1) \times (0,\delta)$. If we choose $s=\frac{\delta}{2}$, Lemma \ref{rig-lemma: weak rig} provides a constant $\sim \delta^{-2}.$ Example \ref{rig-example: FMJ rectangle} shows that this estimate is sharp in the sense that the exponent of $\delta$ cannot be improved. 
\item[(ii)] Following the above arguments we find  that in Lemma \ref{rig-lemma: weak rig} one can replace $p=2$ by any $1 < p < \infty$ replacing $l^4$ suitably by $l^{2p}$. 
\item[(iii)] In view of the proof  in the choice of $R$ we have the freedom to select any of the rotations which are given on each square $Q^s(p)\subset U$ by application of \eqref{rig-eq: R rig}.
\end{itemize}
}
\end{rem}

We briefly note that similar calculations may be provided to estimate the difference of rigid motions. Consider $b_1,b_2 \in \R^2$,  and the rectangles $B_i = b_i + (-l_i, l_i) \times (-m_i,  m_i)  \in {\cal U}^s$ for $i=1,2$, where we assume without restriction that $l_1\ge m_1>0$, $l_2 \ge m_2>0$. Suppose that there is  a point $b_{12} \in \overline{B_1} \cap \overline{B_2}$. For given $R_1, R_2, R_{12} \in SO(2)$ and $c_1,c_2, c_{12} \in \R^2$ we set $E_i := \Vert y - (R_i\, \cdot + c_i)\Vert^2_{L^2(B_i)}$ for $i=1,2$ and assume that 
\begin{align*}
\Vert y - (R_{12} \, \cdot + c_{12})\Vert^2_{L^2(B_1 \cup B_2)}  \le C(E_1+ E_2).
\end{align*}
Then we find 
\begin{align}\label{rig-eq: A difference}
|B_1 \cup B_2| (l_1+ l_2)^2\ |R_1 - R_2|^2 \le C\kappa(E_1 + E_2),
\end{align}
as well as 
\begin{align}\label{rig-eq: A,c difference2}
\Vert y - (R_{1} \, \cdot + c_{1})\Vert^2_{L^2(B_1 \cup B_2)} + \Vert y - (R_{2} \, \cdot + c_{2})\Vert^2_{L^2(B_1 \cup B_2)}  \le C\kappa(E_1+ E_2),
\end{align}
where $\kappa = \frac{|B_1 \cup B_2|}{\min_j |B_j|}\big(\frac{ l_1 + l_2}{ \min_j l_j}\big)^2$. This estimate follows similarly as in the geometrically linear setting treated in \cite[Section 2.2]{Friedrich:15-1} and we therefore omit the details. Indeed, in all the calculations, in particular in (2.10) of \cite{Friedrich:15-1}, one may replace $\R^{2 \times 2}_{\rm skew}$ by $SO(2)$ since the estimates essentially rely on the fact that  $|R \e_1| = |R \e_2|$ which is satisfied for both $\R^{2 \times 2}_{\rm skew}$ and $SO(2)$. Moreover, although we stated this property only  for two rectangles  for the sake of simplicity, we remark that an estimate of the above form also holds for sets with more general geometries.  

Similarly as in \eqref{rig-eq: path}, considering  two arbitrary points $a,b \in I^s(U)$ connected by a path $\xi = (\xi_1,\ldots,\xi_m)$ with corresponding estimates
\begin{align*}
\Vert  (R(\xi_j) -  R(\xi_{j-1})) \, \cdot + c(\xi_j) - c(\xi_{j-1}) \Vert_{L^2(Q^s_{j-1,j})}  \le C E_{j-1,j}, 
\end{align*}
(here we defined $Q^s_{j-1,j} = (\overline{Q^s(\xi_{j-1}) } \cup \overline{Q^s(\xi_{j})})^\circ$) we obtain (cf. (2.20) in \cite{Friedrich:15-1})
\begin{align}\label{rig-eq: A diff}
\begin{split}
 \Vert y - (R(a)\, \cdot + c(a))\Vert^2_{L^2(Q^s(b))}  &\le Cm^2\Big(\sum\nolimits^m_{j=2} E_{j-1,j}\Big)^2 \\ &\le Cm^3\sum\nolimits^m_{j=2} \big(E_{j-1,j}\big)^2.
\end{split}
\end{align}
In the last step we used H\"older's inequality. Similarly as before, \eqref{rig-eq: A diff} also holds for any of the other rigid motions $R(\xi_j)\, x + c(\xi_j)$ (cf. Remark \ref{rig-rem: 1}(iii)).

\section{Modification of sets}\label{rig-sec: subsub, rig-prep}

Before we start with the proof of Theorem \ref{rig-th: rigidity}, we first introduce a procedure to modify sets. In particular, it will be fundamental to assure that during the modification process boundary components do not become too large or are separated by other components.

Recall the definition of the sets ${\cal U}^s$, ${\cal V}^s$ in Section \ref{rig-sec: modifica}. We consider a Lipschitz domain $\Omega \subset \R^2$  and choose $\mu_0$ so  large that $ \overline{\Omega} \subset Q_{\mu_0} = (-\mu_0,\mu_0)^2$.
 We let $\Omega^k$ be the largest set in ${\cal V}^{\bar{c}k}$ satisfying $\Omega^k \subset \lbrace x \in \Omega: \dist(x, \partial \Omega) \ge \bar{c}k\rbrace$ for $k \ge 0$  for some $\bar{c} \ge \sqrt{2}$ large enough.   
 
For sets $W \subset \Omega^k$, $W \in {\cal V}^s$,  we assume that one component in definition \eqref{rig-eq: calV-s def2} is given by $X = Q_{\mu_0} \setminus \Omega^k$. In particular, all other components $X_1,\ldots,X_n$ satisfy $\partial X_i \subset Q_{\mu_0}$ as $\overline{\Omega} \subset Q_{\mu_0}$.  We again choose an (arbitrary) order of $(\Gamma_j)_{ j=1,\ldots,n}$ and define $(\Theta_j)_j$ as in \eqref{rig-eq: Xdef}. Recall the definition of $\Vert \cdot \Vert_X$, $X=*,\infty,{\cal H},$ in \eqref{rig-eq: h*} and \eqref{rig-eq: h**}.  Moreover, we recall that  ${\cal V}^s_{\rm con} \subset {\cal V}^s$ was defined as the subset consisting of the sets where all  $\overline{X_1},\ldots, \overline{X_n}$ are connected.

We now introduce a modification procedure for sets. Given a set $W = Q_\mu \setminus \bigcup^{m}_{i=1} X_i \in {\cal V}^s$ and  some $V \in {\cal U}^s$ we consider the modification   
 \begin{align}\label{rig-eq: new1}
W' = Q_\mu \setminus \bigcup\nolimits^{m}_{i=0} X'_i,
\end{align}
where $X'_i = X_i \setminus \overline{V}$ for $i=1, \ldots, m$ and $X'_{0} = V$. (It is convenient to start with index $0$.) We observe that $W'  = (W \setminus V) \cup \partial V$ (as a subset of $\R^2$). Therefore, for shorthand we will write $W'  = (W \setminus V) \cup \partial V$ to indicate the element of ${\cal V}^s$ which is given by \eqref{rig-eq: new1}. We briefly note that then the boundary components of $W'$ are given by $\Gamma_0(W') = \Theta_0(W') = \partial V$  as well as by  $\Gamma_j(W') = \partial (X_j \setminus \overline{V})$ and  $\Theta_j(W') = \Theta_{j}(W) \setminus \overline{V}$ for $j\ge 1$ (cf. also Lemma \ref{rig-lemma: infty}(iii)).

Having several pairwise disjoint sets $(V_j)_j \subset {\cal U}^s$ the modification is  defined analogously by $W'' = (W \setminus \bigcup_j V_j) \cup \bigcup_j \partial V_j$.

As large surfaces of general shape may not be measured adequately in terms of $|\cdot|_\infty$, in what follows we have to assure that boundary components do not become too large. For $0 < s \le \lambda\le k$ we introduce 
$${\cal V}^s_{(\lambda,k)} := \lbrace W \in {\cal V}^s_{ \rm con}:  2\lambda \le \max\lbrace |\pi_1 \Gamma_j(W)|, |\pi_2 \Gamma_j(W)| \rbrace\le 2k  \  \text{ for all } \Gamma_j(W)\rbrace.$$   
By definition we have $\max\lbrace |\pi_1 \Gamma_j(W)|, |\pi_2 \Gamma_j(W)| \rbrace \ge 2s$ for all $\Gamma_j(W)$ and therefore we write for shorthand ${\cal V}^s_k = {\cal V}^s_{(s,k)}$.

Although we have to avoid that boundary components become to large, it is essential to combine small components. To this end, it is convenient to alter configurations on sets of negligible measure.

\begin{lemma}\label{rig-lemma: small bc}
Let $t \ge 2k$, $t'>0$ and  $W \in {\cal V}^s_t$. 

(i) Then there is a set $\tilde{W} \in {\cal V}^s_{t}$ with $\tilde{W} \subset W$, $|W \setminus \tilde{W}| = 0$ and $\Vert \tilde{W} \Vert_* \le \Vert W \Vert_*$ such that
\begin{align}\label{rig-eq: small bc2}
\Gamma_{j_1}(\tilde{W}) \cap \Gamma_{j_2}(\tilde{W}) = \emptyset \ \ \  \text{ if } |\Gamma_{j_i}(\tilde{W})|_\infty \le k \ \text{ for } i=1,2.
\end{align}

(ii) Then there is a set $U \in {\cal V}^s_{t + k}$ with  $U \subset W$, $|W \setminus U| = 0$ and $\Vert U \Vert_* \le \Vert W \Vert_*$ such that
\begin{align}\label{rig-eq: small bc}
\Gamma(U) \cap \Gamma_j(U) = \emptyset \ \ \ \text{ for all } \Gamma_j(U) \neq \Gamma(U)
\end{align}
for all $\Gamma(U)$ with $|\Gamma(U)|_\infty \le k$.
\end{lemma}

\Proof (i) The strategy is to combine iteratively different boundary components. Clearly, if $|\Gamma_{j_i}(W)|_\infty \le k$ for $i=1,2$ with $\Gamma_{j_1}(W) \cap \Gamma_{j_2}(W) \neq \emptyset$ we may replace $W$ by $W' = W \setminus (\overline{X_{j_1} \cup X_{j_2}})^\circ$  and note that $W' \in {\cal V}^s_{t}$ as well as $|W \setminus W'| = 0$ and $\Vert W'\Vert_* \le \Vert W \Vert_*$ similarly as in Lemma \ref{rig-lemma: infty}. (Recall that $\partial X_{j_i} = \Gamma_{j_1}(W)$ for $i=1,2$.) We proceed in this way until we obtain a set $\tilde{W} \in {\cal V}^s_{t}$ with $|W \setminus \tilde{W}| = 0$ and $\Vert \tilde{W} \Vert_* \le \Vert W \Vert_*$ such that \eqref{rig-eq: small bc2} holds.

(ii)    We apply (i) and then proceed to combine two components $\Gamma_{j_1}(\tilde{W}),  \Gamma_{j_2}(\tilde{W})$  if $\Gamma_{j_1}(\tilde{W}) \cap \Gamma_{j_2}(\tilde{W}) \neq  \emptyset$ and  $\min \lbrace |\Gamma_{j_1}(\tilde{W})|_\infty, |\Gamma_{j_2}(\tilde{W})|_\infty \rbrace \le k$. Arguing as before we end up with a set $U$ satisfying $|W \setminus U| = 0$,   $\Vert U \Vert_* \le \Vert W \Vert_*$ and \eqref{rig-eq: small bc}. It remains to show that $U \in {\cal V}^s_{t + k}$.  Consider some $\Gamma(U) = \partial X$ with $|\Gamma(U)|_\infty > k$ and observe that there are $\Gamma(\tilde{W}) = \partial X'$ with $|\Gamma(\tilde{W})|_\infty > k$ and $\Gamma_{j_i}(\tilde{W}) = \partial X_{j_i}$, $i=1,\ldots,m$, with $|\Gamma_{j_i}(\tilde{W})|_\infty \le k$, $\Gamma_{j_{i_1}}(\tilde{W}) \cap \Gamma_{j_{i_2}}(\tilde{W}) = \emptyset$ for $i_1 \neq i_2$ and $\Gamma_{j_i}(\tilde{W}) \cap \Gamma(\tilde{W}) \neq \emptyset$ such that $\overline{X} = \overline{X' \cup \bigcup^m_{i=1} X_{j_i}}$. But this implies $|\pi_i \Gamma(U)| \le 2k + |\pi_i \Gamma(\tilde{W})| \le  2k + 2t$ for $i=1,2$, as desired.\eop

In what follows we often modify sets by subtracting rectangular neighborhoods of boundary components. In this context it is particularly important that the components remain connected and do not become too large. By  $\triangle$ we denote the symmetric difference of two sets.

\begin{lemma}\label{rig-lemma: modifica}
Let $k,t,t'>0$  with $t, t' \le Ck$ and  $\nu \ge 0$. Let $V\subset \Omega^k$ with  $V \in {\cal V}^s_{\rm con}$.

(i) Assume that for each component $X_j  = X_j(V)$, $j=1,\ldots,n$, there is a rectangle $Z_j \in {\cal U}^s$ with $X_j \subset Z_j$, $|\pi_i \partial Z_j| \le |\pi_i\partial X_j| +  \nu|\partial X_j|_\infty$ for $i=1,2$ and $\max_{i=1,2}|\pi_i\partial Z_j| \le 2t'$ for all $j=1,\ldots,n$. Moreover, assume that $Z_{j_1} \setminus Z_{j_2}$ or $Z_{j_2} \setminus Z_{j_1}$ is connected for all $1 \le j_1 < j_2 \le n$. Then there is a set $U \in {\cal V}^s_{t'}$,  $U \subset \Omega^k$, with $\bigcup^n_{j=1} \overline{X_j(U)} = \bigcup^n_{j=1} \overline{Z_j}  \cap \Omega^k$ and $\Vert U \Vert_* \le (1+c\nu)\Vert V\Vert_*$ for a universal constant $c>0$.

(ii) In addition let $V' \in {\cal V}^s_t$  be given and define  $\hat{W} = V' \setminus \bigcup^{n}_{j=1} Z_j$.  Then there is a set $W \in {\cal V}^{s/2}_{ t + 2t'}$ with  $|W \setminus \hat{W}| = 0$, $|\hat{W} \setminus W| \le ct \Vert V' \Vert_*$  and $\Vert W \Vert_* \le (1+c\nu)\Vert V \Vert_* + \Vert V' \Vert_*$  . 
\end{lemma}

As the proof of this result is very technical and in principle not relevant to
understand the proof of the main result in the subsequent sections, it may be omitted on first
reading.

\Proof (i) Let $V \subset \Omega^k$ with components $(X_j)^n_{j=1}$ and rectangles $(Z_j)^{n}_{j=1}$ be given.  It suffices to show the following: There are  connected, pairwise disjoint $(X'_j)^{n}_{j=1}$ with $X'_j \subset Z_j$, $\bigcup_{j=1}^{n} \overline{X_j'} = \bigcup_{j=1}^{n} \overline{Z_j}$ and  
\begin{align}\label{rig-eq: mayblast}
\big|\bigcup\nolimits^{n}_{j=1} \partial X_j'\big|_{\cal H}& \le  \sum\nolimits^{n}_{j=1} |\Theta_j(V)|_{\cal H} + c\nu\sum\nolimits^n_{j=1} |\Gamma_j|_{\cal H}.
\end{align}
 Moreover, we have $X_j' = R_j \setminus \overline{(A^1_j \cup A^2_j)}$. Here $R_j \in {\cal U}^s$ is a rectangle and $A^i_j \in {\cal U}^s$, $i=1,2$, are (if nonempty) unions of rectangles whose closure intersect the corner $c^i_j \in \partial R_j$, where $c^1_j$, $c^2_j$ are adjacent corners of $R_j$.

Then the claim of the lemma follows for $U=\Omega^k \setminus \bigcup^n_{j=1} X_j'$. Indeed, to see $\Vert U\Vert_* \le (1+c\nu)\Vert V \Vert_*$ we first observe $\sum_j |\partial X'_j|_\infty \le \sum_j |\partial Z_j|_\infty \le (1+c\nu) \sum_j |\partial X_j|_\infty$. Moreover, by \eqref{rig-eq: mayblast} we get 
\begin{align}\label{rig-eq: mayblast3}
\Vert U \Vert_{\cal H} \le |\bigcup\nolimits^n_{j=1} \partial X'_j|_{\cal H}  \le (1+c\nu)\Vert V \Vert_{\cal H} = (1+c\nu)|\bigcup\nolimits^n_{j=1} \partial X_j|_{\cal H}.
\end{align}
In the first inequality we also used $|\partial X_j'|_\infty \le |\partial Z_j|_\infty \le Ck$ and $\Omega^{k} \in {\cal V}^{\bar{c}k}$ for  $\bar{c} \gg 1$. (Arguments of this form will be used frequently in the following and from now on we will omit the details.) Finally, we conclude $U \in {\cal V}^s_{t'}$ as $\max_{i=1,2}|\pi_i \partial Z_j|\le 2t'$ for $j=1,\ldots,n$.

We prove the above assertion by induction. Clearly, the claim holds for $n=1$ for $X_1' = Z_1$. Now assume the assertion holds for sets with at most $n-1$ components and consider $V \subset \Omega^k$ with components $(X_j)^{n}_{j=1}$ and corresponding $(Z_j)^n_{j=1}$. Without restriction we assume that $\max_{x \in \overline{Z_n}} x_2 = \max_{x \in \bigcup_{j=1}^n\overline{Z_j}} x_2$.  By hypothesis we obtain pairwise disjoint, connected sets $X_j''$, $j=1,\ldots,n-1$, fulfilling the above properties, in particular  $\bigcup^{n-1}_{j=1} \overline{X''_j}  = \bigcup^{n-1}_{j=1} \overline{Z_j}$. 

 Given $Z_n = (z^1_1,z^2_1) \times (z^1_2,z^2_2)$ we set $\tilde{Z}_n = (z^1_1,z^2_1) \times (z^1_2,z^2_2]$. For $j=1,\ldots,n-1$ let $Z'_{j,i} \in {\cal U}^s$ be the largest rectangle in $Z_n$ satisfying $Z_j \cap Z_n \subset Z'_{j,i} \subset \bigcup_{j=1}^{n-1} \overline{Z_j}$ with $z^i_1 \in \overline{Z'_{j,i}}$ for $i=1,2$. If $Z'_{j,i} \neq \emptyset$ for some $i$, we let $Z'_j = Z'_{j,i}$, otherwise we set $Z'_j = Z_j \cap Z_n$. (Note that $Z'_{j,1} = Z'_{j,2}$ if $Z'_{j,1}, Z'_{j,2} \neq \emptyset.$ )  

Let $J_0 \subset \lbrace 1, \ldots,n-1 \rbrace$ such that $Z_j \cap Z_n = \emptyset$ for $j \in J_0$. Let $J_1 \subset \lbrace 1, \ldots,n-1 \rbrace \setminus J_0$ such that  $(\overline{Z'_j} \setminus Z_n) \cap \lbrace z_1^1, z_1^2 \rbrace = \emptyset$ for $j \in J_1$ and $J_2 \subset \lbrace 1, \ldots,n-1 \rbrace \setminus J_0$ such that $ \tilde{Z}_n \setminus Z'_j$ is a rectangle for $j \in J_2$. (Observe that $J_1 \cap J_2 = \emptyset$.)  Let $J_3 = \lbrace 1, \ldots,n-1 \rbrace \setminus (J_0 \cup J_1 \cup J_2)$. Define $X_n' = Z_n \setminus \bigcup_{j \in J_2 \cup J_3} \overline{Z'_j}$. Moreover, we let $X_j' = X_j''$ for $j \in J_0 \cup J_2 \cup J_3$ and $X_j' = X_j'' \setminus \overline{X_n'}$ for $j \in J_1$. Clearly, by construction the sets are pairwise disjoint and fulfill $\bigcup_{j=1}^n \overline{X_j'} = \bigcup_{j=1}^n \overline{Z_j}$. 

Moreover, we observe that the sets are connected and have the special shape given above. In fact, for $j \in J_0 \cup J_2 \cup J_3$ this is clear. For $X_n'$ we first note that $J_3 = J^1_3 \dot{\cup} J_3^2$, where $\overline{Z'_j}$ intersects the lower right and the lower left corner of $Z_n$ for $j \in J_3^1$ and $j \in J_3^2$, respectively. (It cannot happen that $\overline{Z'_j}$ intersects only the other corners due to the choice of $Z_n$.) We observe $ X_n' = R_n \setminus \overline{(A^1_n \cup A^2_n)}$ is connected, where $R_n = Z_n \setminus \bigcup_{j \in J_2} \overline{Z'_j}$ and $A^i_n = \bigcup_{j \in J_3^i} Z'_j$ for $i=1,2$.

Finally, to see the properties for  $j \in J_1$ we first observe that $S_j := Z_j \setminus \overline{X'_{n}}$ is a rectangle. In fact, otherwise due to the special shape of $X_n'$ it is elementary to see that  $(\overline{Z'_j} \setminus Z_n) \cap \lbrace z_1^1, z_1^2 \rbrace \neq \emptyset$ and thus $j \notin J_1$.  We get  $X_j' = X_j'' \cap S_j = \big( R_j \cap S_j \big) \setminus \overline{(A_j^1 \cup A_j^2)}$ is connected and $X_j' = \hat{R}_j \setminus \overline{(\hat{A}^1_j \cup \hat{A}^2_j)}$, where $\hat{R}_j = S_j$ and $\hat{A}^i_j = A_j^i  \cap S_j$ for $i=1,2$. 

It remains to confirm \eqref{rig-eq: mayblast}. We first observe that
\begin{align}\label{rig-eq: mayblast2}
\sum\nolimits^{n}_{j=1} |\Theta_j(V)|_{\cal H} = \tfrac{1}{2} \sum\nolimits^n_{j=1} |\Gamma_j|_{\cal H} + \tfrac{1}{2}\big|\partial\big( \bigcup\nolimits^n_{j=1}  \overline{X_j}\big)\big|_{\cal H}.
\end{align}
(Recall that different boundary components may have nonempty intersection.)
Similarly, for the components $(X_j')_j$ we find
$$\big|\bigcup\nolimits^{n}_{j=1} \partial X_j'\big|_{\cal H} = \tfrac{1}{2} \sum\nolimits^n_{j=1} |\partial X_j'|_{\cal H} + \tfrac{1}{2} \big|\partial \big(\bigcup\nolimits^n_{j=1} \overline{X_j'}\big)\big|_{\cal H}.$$
We now treat the two terms on the right separately. By $T_j \in {\cal U}^s$ we denote the smallest rectangle containing $X_j$ and observe that $|\partial T_j|_\infty = |\Gamma_j|_\infty$, $|\partial T_j|_{\cal H}\le |\Gamma_j|_{\cal H}$. Recall $|\partial Z_j|_{\cal H} \le  |\partial T_j|_{\cal H} + c\nu |\partial T_j|_{\infty} \le (1 + c\nu) |\Gamma_j|_{\cal H}$ for $j=1,\ldots,n$. Due to the special shape of the components $X_j'$ we find $|\partial X_j'|_{\cal H} \le |\partial Z_j|_{\cal H}$ and thus 
\begin{align}\label{rig-eq: mayblast2****}
\sum\nolimits^n_{j=1} |\partial X_j'|_{\cal H} \le (1 + c\nu)\sum\nolimits^n_{j=1} |\Gamma_j|_{\cal H}.
\end{align} 
Moreover, it is elementary to see that we can find a connected set $\tilde{X}_j \supset X_j$ such that $\tilde{\Gamma}_j := \partial \tilde{X}_j$ satisfies $|\tilde{\Gamma}_j|_{\cal H} \le (1+c\nu)|\Gamma_j|_{\cal H}$ and $Z_j \in {\cal U}^s$ is the smallest rectangle containing $\tilde{X}_j$. By a projection argument it is then not hard to see that 
\begin{align*}
\big|\partial \big(\bigcup\nolimits^n_{j=1} \overline{X_j'}\big)\big|_{\cal H} &= \big|\partial \big(\bigcup\nolimits^n_{j=1} \overline{Z_j}\big)\big|_{\cal H} \le \big| \partial \big( \bigcup\nolimits^n_{j=1}  \overline{\tilde{X}_j} \big) \big|_{\cal H}\\
&\le \big| \partial \big( \bigcup\nolimits^n_{j=1}  \overline{X_j} \big) \big|_{\cal H}+ c\nu \sum\nolimits |\Gamma_j|_{\cal H}.
\end{align*}
Consequently, we derive by \eqref{rig-eq: mayblast2} and \eqref{rig-eq: mayblast2****}
\begin{align*}
\big|\bigcup\nolimits^{n}_{j=1} \partial X_j'\big|_{\cal H} &\le \tfrac{1}{2} \sum\nolimits^n_{j=1} |\Gamma_j|_{\cal H} + \tfrac{1}{2}\big|\partial \big( \bigcup\nolimits^n_{j=1} \overline{X_j} \big) \big|_{\cal H}+  c\nu\sum\nolimits^n_{j=1} |\Gamma_j|_{\cal H} \\
&  = \sum\nolimits^{n}_{j=1} |\Theta_j(V)|_{\cal H} +  c\nu\sum\nolimits^n_{j=1} |\Gamma_j|_{\cal H},
\end{align*}
as desired.

\vspace{0.2cm}
(ii) Let  $(Y_j)^{n'}_{j=1}$ be the components of $V'$ and let $T_j \in {\cal U}^s$ be the smallest rectangle containing $Y_j$. It is elementary to see that $T_{j_1} \setminus T_{j_2}$ is connected for $1 \le j_1, j_2 \le n'$. Thus, by (i) we obtain pairwise disjoint, connected sets $(Y_j')_j$ with $\bigcup_j \overline{Y_j'} = \bigcup_j \overline{T_j}$ and define $V'' =  \Omega^k \setminus \bigcup^{n'}_{j=1} Y'_j$. By (i) for $\nu = 0$ we then also obtain $\Vert V'' \Vert_* \le \Vert V' \Vert_*$. Moreover, the isoperimetric inequality yields $|V'\setminus V''| \le c t \Vert V' \Vert_*$ since $|\partial T_j|_\infty \le 2\sqrt{2}t$ for all $j=1,\ldots,n'$.

 Let $(X'_j)^n_{j=1}$ and $U \in {\cal V}^s_{t'}$ as given in (i). We define $W' = (U \setminus \bigcup_{j=1}^{n'} Y'_j) \cup \bigcup_{j=1}^{n'} \partial Y'_j$. Clearly, we have $|W' \setminus \hat{W}|=0$,  $|\hat{W} \setminus W'| \le ct \Vert V' \Vert_*$ and $\Vert W' \Vert_* \le (1+c\nu) \Vert V \Vert_* + \Vert V' \Vert_*$ arguing similarly as in  Lemma \ref{rig-lemma: infty}. Observe that possibly $W' \notin {\cal V}^s_{\rm con}$ as the components $(X'_j)^{n}_{j=1}$ of $U$ may have become disconnected. Thus, we now construct a set $W \in {\cal V}^{s/2}_{\rm con}$ with $|W' \triangle W|=0$.
 
By $R_j \in {\cal U}^s$ we denote  the smallest rectangle such that $X'_j \subset R_j$ for $j=1,\ldots,n$ and observe $\bigcup_j \overline{R_j} = \bigcup_j \overline{X'_j}$. To simplify the exposition we assume that each of the components $(X'_j)_j$ has become disconnected as otherwise we do not have to alter the boundary component in the  modification procedure described below. Moreover, we can suppose that for each pair $Y_{j_1}'$, $X_{j_2}'$, $1 \le j_1 \le n'$, $1 \le j_2 \le n$, with $R_{j_2} \setminus Y'_{j_1}$ is not disconnected we have $X'_{j_2} \setminus Y'_{j_1}$ is not disconnected. In fact, otherwise we can pass to some $Y^*_{j_1} \subset Y'_{j_1}$ with $|\partial Y^*_{j_1}|_* \le |\partial Y'_{j_1}|_*$ such that $X'_{j_2} \setminus Y^*_{j_1}$ is not disconnected and $\bigcup_j \overline{Y'_j} \cup \bigcup_j \overline{X'_j} = \bigcup_j \overline{Y^*_j} \cup \bigcup_j \overline{X'_j}$.

We now proceed by induction. Let $W_0 = V''$ and $T_j^0 = Y_j'$ for $j=1,\ldots,n'$. Assume there are pairwise disjoint, connected sets $T^{l-1}_j \in {\cal U}^{\frac{s}{2}}$, $j=1,\ldots,n'$ such that 
\begin{align}\label{rig-eq: modifica1}
(i) \ \ \bigcup\nolimits^{n'}_{j=1} \overline{T^{l-1}_j} = \bigcup\nolimits^{n'}_{j=1} \overline{Y'_j} \cup \bigcup\nolimits^{l-1}_{j=1} \overline{X'_j}, \ \ \ (ii) \ \  T^{l-1}_{j_1} \cap  \overline{X'_{j_2}} = T^{0}_{j_1} \cap  \overline{X'_{j_2}}
\end{align}
for all $1 \le j_1 \le n'$, $l \le j_2 \le n$. Moreover, assume that the set $W_{l-1} := \Omega^k \setminus \bigcup_j T^{l-1}_j$ satisfies $\Vert W_{l-1} \Vert_\infty \le \sum_j |\partial T^0_{j}|_\infty + \sum^{l-1}_{i=1} |\partial X_i'|_\infty$ and 
\begin{align}\label{rig-eq: modifica2}
\Vert W_{l-1}\Vert_{\cal H} \le |\bigcup\nolimits_{j} \partial T^0_{j} |_{\cal H} + |\bigcup\nolimits^{l-1}_{i=1} \partial X_i' \setminus \bigcup\nolimits_j T^{0}_j|_{\cal H}  +  \frac{1}{2}\sum^{l-1}_{i=1}|\partial X_i' \cap \bigcup\nolimits_{j} T^{0}_j|_{\cal H}.
\end{align}
We now construct $W_l$.  Let $J^l \subset \lbrace 1,\ldots,n' \rbrace$ such that $T_j^{l-1} \cap X_l' \neq \emptyset$ with $J^l = J^l_1 \dot{\cup} J^l_2$, where $j \in J^l_2$ if and only if $R_l \setminus T_j^{l-1}$ is disconnected.

If $j \in J^l_1$, we define $T_j^l = T_j^{l-1} \setminus \hat{X}_l'$, where $\hat{X}_l'\in  {\cal U}^{\frac{s}{2}}$ is the largest set with $\overline{\hat{X}_l'} \subset X_l'$. It is not hard to see that $|\partial T_j^l|_{\infty} \le |\partial T_j^{l-1}|_{\infty}$ for all $j \in J^l_1$ and $|\partial T_j^l|_{\cal H} \le |\partial T_j^{l-1} \setminus X_l'|_{\cal H} + \frac{1}{2} |\partial T_j^{l-1} \cap X'_l|_{\cal H} + \frac{1}{2} |\partial X_l' \cap T_j^{l-1}|_{\cal H}$. As each $x \in \R^2$ is contained in at most two different $\partial T^{l-1}_j$, we find $\sum_{j \in J_1^l} \frac{1}{2} |\partial T^{l-1}_j \cap X_l'|_{\cal H} \le |\bigcup_{j \in J_1^l} \partial T_j^{l-1} \cap X_l'|_{\cal H}$. Therefore, taking the union over all components    we derive
\begin{align}\label{rig-eq: modifica3}
\big|\bigcup\nolimits_{j \in J^l_1} \partial T_j^l \cup \bigcup\nolimits_{j \notin J^l_1} \partial T_j^{l-1}\big|_{\cal H}  \le \big|\bigcup\nolimits_{j} \partial T_j^{l-1}\big|_{\cal H}  + \tfrac{1}{2} \big|\partial X_l' \cap \bigcup\nolimits_{j \in J^l_1}T_j^{0}\big|_{\cal H}.
\end{align} 
Here we used \eqref{rig-eq: modifica1}(ii)  and the fact that the sets $(T^{l-1}_j)_j$ are pairwise disjoint. Observe that the above construction together with \eqref{rig-eq: modifica1}(ii) and the special shape of $T_j^0$ (see proof of (i)) implies that the sets $T^l_j$, $j \in J_1^l$, are connected. Moreover, \eqref{rig-eq: modifica1}(ii) holds for $j_1 \in J_1^l$. 

 We define $\tilde{X}_l' = X_l' \setminus \bigcup_{j \in J^l_1} \overline{T_j^l} \in {\cal U}^{\frac{s}{2}}$. Due to the fact that $\hat{X}_l' \neq \emptyset$ we observe that the number of connected components  of the sets $X_l' \setminus \bigcup_{j \in J^l_2} T^{l-1}_j$ and $\tilde{X}_l' \setminus \bigcup_{j \in J^l_2} T^{l-1}_j$ coincide. Therefore, letting  $A_1,\ldots,A_{m}$ be the connected components of $\tilde{X}_l' \setminus \bigcup_{j \in J^l_2} \overline{T^{l-1}_{j}}$ it is elementary to see that $m = \# J^l_2 +1$. 
 
 Up to a rotation by $\frac{\pi}{2}$ we can assume that each $\overline{A_i}$ intersects the upper and lower boundary of $R_l$ and that $\overline{A_1}$ intersects the left boundary. For convenience we denote the components $(T^{l-1}_j)_{j \in J^l_2}$ by $(T^{l-1}_{j_i})^{m-1}_{i=1}$. Suppose $R_l= (0,l_1) \times (0,l_2)$. Let $a_i = \inf_{x \in A_i} x_1$ and $d_i = a_{i+1} - a_{i}$, where $a_{m+1} = l_1$. Define $T^l_{j_1} = (\overline{T^{l-1}_{j_1}} \cup \overline{(A_1\cup A_2) })^\circ$ and $T^l_{j_i} = (\overline{T^{l-1}_{j_i}} \cup \overline{A_{i+1}})^\circ$ for $i=2,\ldots,m-1$. Observe that the sets are pairwise disjoint, connected  and that \eqref{rig-eq: modifica1}(ii) holds for $j_i \in J_2^l$.  It is elementary to see that $|T^l_{j_1}|_\infty \le |T^{l-1}_{j_1}|_\infty + d_1 + d_2$ and $|T^l_{j_i}|_\infty \le |T^{l-1}_{j_i}|_\infty + d_{i+1}$ for $i=2,\ldots,m-1$. Thus, we have 
 \begin{align}\label{rig-eq: mayblast9}
 \sum\nolimits^{m-1}_{i=1} |T^l_{j_i}|_\infty \le \sum\nolimits^{m-1}_{i=1} |T^{l-1}_{j_i}|_\infty + |X_l'|_\infty.
 \end{align}
For $j \notin J^l$ we define $T^l_j = T_j^{l-1}$ and observe that  \eqref{rig-eq: modifica1}(i) holds by construction and the assumptions before \eqref{rig-eq: modifica1}. Together with \eqref{rig-eq: modifica3} and \eqref{rig-eq: modifica1}(ii) we then also get 
\begin{align*}
\big|\bigcup\nolimits_{j} \partial T^l_{j} \big|_{\cal H} \le \big|\bigcup\nolimits_{j} \partial T^{l-1}_{j} \big|_{\cal H} + \big|\partial X_l' \setminus \big(\bigcup\nolimits^{l-1}_{i=1} \partial X_i' \cup \bigcup\nolimits_{j} T_j^{0}\big)\big|_{\cal H} + \tfrac{1}{2}\big|\partial X_l' \cap \bigcup\nolimits_{j} T^{0}_j\big|_{\cal H}.
\end{align*}
This in conjunction with \eqref{rig-eq: modifica2} for $W_{l-1}$ implies that \eqref{rig-eq: modifica2} holds for $W_l$. Moreover, by \eqref{rig-eq: mayblast9} it is elementary to see that $\Vert W_{l} \Vert_\infty \le \sum_j |\partial T^0_{j}|_\infty + \sum^{l}_{i=1} |\partial X_i'|_\infty$.

Finally, we define $W = W_n$ and observe that $W$ has the desired properties. In fact, by \eqref{rig-eq: modifica1}(i) we have $|W \triangle W'|= 0$ and thus $|\hat{W} \setminus W| \le ct\Vert V' \Vert_*$. Moreover, we clearly get $\Vert W \Vert_\infty \le \Vert U \Vert_\infty + \Vert V'' \Vert_\infty \le (1+c\nu) \Vert V \Vert_\infty+ \Vert V'\Vert_\infty$.  As each $x \in \R^2$ is contained in at most two different $\partial X'_l$, we find by \eqref{rig-eq: modifica2}
\begin{align*}
\Vert W \Vert_{\cal  H} &\le \Vert V'' \Vert_{\cal H} + \big|\bigcup\nolimits^n_{i=1} \partial X'_i \setminus \bigcup\nolimits_j T_j^0\big|_{\cal H} + \big|\bigcup\nolimits^n_{i=1} \partial X'_i \cap\bigcup\nolimits_j T_j^0\big|_{\cal H} \\ & = \Vert V'' \Vert_{\cal H} + \Vert U \Vert_{\cal H} \le  \Vert V' \Vert_{\cal H} + (1+c\nu)\Vert V \Vert_{\cal H}, 
\end{align*}
as desired. Finally, similarly as in Lemma \ref{rig-lemma: small bc}(ii) we obtain $|\pi_i X_j(W)|\le 2t + 4t'$ for $i=1,2$ for all $j$ and thus $W \in {\cal V}^{s/2}_{t+2t'}$.  \eop

\section{A local rigidity estimate}\label{rig-sec: sub, first-weak}

We now establish a local rigidity estimate on a fine partition of the Lipschitz domain $\Omega$. As a preparation we  introduce some further notions.  Recall the point set $I^s = s(1,1) +  2 s\Z^2$, $s>0$, introduced in Section \ref{rig-sec: modifica} and the definitions of ${\cal U}^s, {\cal V}^s$ in \eqref{rig-eq: calV-s def}, \eqref{rig-eq: calV-s def2}  with respect to the square $Q_{\mu_0}$.  We define additional partitions. Set $z_1 = (0,0)$, $z_2 = (1,0)$, $z_3=(0,1)$, $z_4=(1,1)$ and let $I^s_i = s z_i + 2s\Z^2$ as well as $Q^s_i(p) = p + s(-1,1)^2$ for $p \in I^s_i$, $i=1,\ldots,4$.  Moreover, for $U\subset \Omega$ let 
\begin{align*}
I^s_i(U) = \lbrace p \in I^s_i: Q^s_i(p) \subset U \rbrace
\end{align*}
for $i=1,\ldots, 4$. For shorthand we also write $I^s = I^s_4$ and $Q^s = Q^s_4$.

In the following, constants which are much smaller than $1$ will frequently appear. For the sake of convenience we introduce one universal parameter. For  given $l \ge 1$ and $0 <s,\epsilon, m \le 1$ we let
\begin{align}\label{rig-eq: vartheta def}
\vartheta = l^{ 9} C^{2}_m s^{-1} \epsilon,   
\end{align}
where $C_m  = C_1(m,{h_*})  + C_3(m,{h_*}) +  m^{-4} C^{-2}_2(m,{h_*})$ with the constants of Theorem \ref{rig-th: derive prop} and Lemma \ref{rig-lemma: A neigh} (for fixed $h_*$). By Remark \ref{rig-rem: z}(i)  we find  some $z \in \N$ such that $C_m \le C(h_*)m^{-z}$. Moreover, for later let $\hat{m} =  C_2(m,{h_*})$ and recall that by Remark \ref{rig-rem: z}(ii) we  can assume $\hat{m} \ll m$  as well as $ \bar{C}\hat{m} \le m$ for constants $\bar{C} = \bar{C}(h_*)$.  Using only one universal parameter the estimates we establish are often not sharp. However, this will not affect our analysis. 

\begin{rem}\label{rig-rem: h_*}
{\normalfont
All the constants $C$ in the following may depend on $h_*$ unless they are universal constants indicated as $C_{\rm u}$. However, to avoid further notation we drop the dependence here. Only at the end of the proof in Section \ref{rig-sec: sub, proof-main}, when we pass to the limit $h_*\to 0$, we will take the $h_*$ dependence of the constants into account. 
}
\end{rem}

 In the following, $\epsilon$ will represent the stored elastic energy. We first construct piecewise constant $SO(2)$-valued mappings approximating the deformation gradient. Afterwards, we employ Theorem \ref{rig-th: derive prop} and Corollary \ref{rig-cor: kornpoin} to find a piecewise rigid motion being a good approximation of the deformation.

\subsection{Estimates for the derivatives}\label{rig-sec: subsub,  est-deriv}

We divide our investigation into two regimes, the `superatomistic' $k \ge \epsilon$ and the `subatomisic' $k \le \epsilon$. Here, recall that we call the $\epsilon$-regime the `atomistic regime' as in discrete fracture models $\epsilon$ is of the same order as the typical interatomic distance  (c.f. \cite{FriedrichSchmidt:2014.1, FriedrichSchmidt:2014.2}). We begin with the superatomistic regime.

\begin{lemma}\label{rig-lemma: weaklocR}
Let $k > s\ge \epsilon > 0$ with $1 \ll l:= \frac{k}{s}$. Let $m^{-1} \in \N$ and assume that $\frac{km}{s} \in \N$.  Then for a constant $C>0$ we have the following: \\
For all $U \in {\cal V}^s_k$ with $U \subset \Omega^k$ and for all $y \in H^1(U)$ with $\Delta y = 0$ in  $U^\circ$ and
\begin{align}\label{rig-eq: weakloc gamma}
\gamma := \Vert\dist(\nabla y, SO(2))\Vert^2_{L^2(U)},
\end{align}
there is a set  $W \in {\cal V}^{sm}_{(s,3k)}$ with $W \subset \Omega^{3k}$, $|W \setminus U|=0$,  $|(U\setminus W) \cap \Omega^{3k} | \le C_{u}k\Vert  W \Vert_* $ and
\begin{align}\label{rig-eq: weaklocR2}
\Vert W \Vert_* \le (1+C_{u} m)\Vert U \Vert_* + C\epsilon^{-1} \gamma.
\end{align}
Moreover, there are mappings $\hat{R}_i:  W^\circ \to SO(2)$, $i=1,\ldots,4$, which are constant on the connected components of $Q^{k}_i(p) \cap  W^\circ$, $p \in I^{k}_i(\Omega^k)$, such that
\begin{align}\label{rig-eq: weaklocR1}
\begin{split}
(i) & \ \ \Vert \nabla y  - \hat{R}_i\Vert^2_{L^2(W)} \le C  l^4 \gamma,\\
(ii) & \ \ \Vert \nabla y  - \hat{R}_i\Vert^4_{L^4(W)} \le C  \vartheta \gamma.
\end{split}
\end{align}
\end{lemma}

\Proof We first construct the set $W$. Let $J \subset I^k(\Omega^k)$ such that
\begin{align}\label{rig-eq: weaklocR3}
\Vert \dist(\nabla y,SO(2))\Vert^2_{L^2(Q^k(p) \cap U)} >  \epsilon k 
\end{align}
for all $p \in J$. Define 
$$\hat{W} = \Big(U \setminus \bigcup\nolimits_{p \in J} Q^k(p)\Big)  \cup \bigcup\nolimits_{p \in J} \partial Q^k(p)$$
and note that $\hat{W} \in {\cal V}^s_k$.  In particular, the property $\hat{W} \in {\cal V}^s_{\rm con}$ holds since $\max\lbrace |\pi_1\Gamma_t(U)|,|\pi_2\Gamma_t(U)| \rbrace \le 2k$. The fact that we add the union of the boundary on the right hand side assures that we do not `combine' boundary components.  Moreover, we derive $\Vert \hat{W} \Vert_* \le \Vert U \Vert_* + C\epsilon^{-1} \gamma$. Indeed, $\sum_{p\in J} |\partial Q^k_p|_* \le 8k \cdot \# J \le 8k \frac{\gamma}{\epsilon k}$  by \eqref{rig-eq: weakloc gamma}. For all other $\Gamma_t(\hat{W})$ we find a corresponding $\Gamma_t(U)$ (without restriction we use the same index) such that $\Theta_t(\hat{W}) = \Theta_t(U)\setminus \bigcup_{p\in J} \overline{Q^k(p)}$  and thus $|\Theta_t(\hat{W})|_* \le |\Theta_t(U)|_*$. (Arguments of this form will be used frequently in the following and from now on we will omit the details.) Furthermore, we easily deduce $|U \setminus \hat{W}| \le C_{u} k\Vert \hat{W} \Vert_*$. 

Then we can find a set $W \in {\cal V}^{sm}_{2k}$ with $\Vert W \Vert_* \le (1+C_{u}m)\Vert\hat{W} \Vert_*$, $|U \setminus  W| \le C_{u}k\Vert  W \Vert_*$ and  $W^\circ \subset \lbrace x \in \Omega^{3k} \cap \hat{W}: \dist_\infty(x, \partial \hat{W})    \le 2 sm \rbrace$, where $\dist_\infty(x,A) := \inf_{y \in A} \max_{i=1,2} |(x - y) \cdot \e_i|$ for $A \subset \R^2$, $x \in  \R^2$. 

Indeed, let $M(\Gamma_j) \in {\cal U}^{sm}$ be the smallest rectangle satisfying $M(\Gamma_j) \supset \lbrace x \in \R^2: \dist_\infty(x, X_j) \le 2sm \rbrace$, where $X_j$ denotes the component corresponding to $\Gamma_j(\hat{W})$. Clearly, we obtain $|\pi_i \partial M(\Gamma_j)| \le |\pi_i \Gamma_j(\hat{W})|+C_u m |\Gamma_j(\hat{W})|_\infty$ for $i=1,2$, $j=1,\ldots,n$ as $\hat{W} \in {\cal V}^s$. Moreover, it is elementary to see that $M(\Gamma_{j_1}) \setminus M(\Gamma_{j_2})$ is connected for $1 \le j_1 , j_2 \le n$ since $sm \ll s$. Then by Lemma \ref{rig-lemma: modifica}(i) with $Z_j = M(\Gamma_j)$ we obtain a  set $W \in {\cal V}^{sm}_{2k}$ which coincides with  
\begin{align}\label{rig-eq: new111}
\Omega^{3k} \cap  \big(\hat{W} \setminus \bigcup\nolimits^n_{j=1} M(\Gamma_j)\big) = \Omega^{3k} \setminus \bigcup\nolimits^n_{j=1} M(\Gamma_j)
\end{align}
up to  a set of negligible measure. Here we used $sm \ll k$. Moreover, we have $|(U \setminus W) \cap \Omega^{3k}| \le C_{u}k\Vert W\Vert_*$ and $\Vert W \Vert_* \le (1+C_{u}m)\Vert\hat{W} \Vert_*$. 

Boundary components of $W$ are possibly smaller than $2s$ due to the modification in \eqref{rig-eq: new111}. Therefore, we apply Lemma \ref{rig-lemma: small bc}(ii) to get a (not relabeled) set $W \in {\cal V}^{sm}_{3k}$ such that \eqref{rig-eq: weaklocR2} still holds and \eqref{rig-eq: small bc} is satisfied. Now the fact that $U \in {\cal V}^s_{(s,k)}$ and \eqref{rig-eq: small bc} imply $W \in {\cal V}^{sm}_{(s,3k)}$.

Fix $i=1,\ldots,4$ and let $F \subset Q^{k}_i(p) \cap  W^\circ$ be a connected component of $Q^{k}_i(p) \cap  W^\circ$. Define $\hat{F} \in {\cal U}^{s}$ as the smallest  (connected) set satisfying 
$$\hat{F} \supset \lbrace x: \dist_\infty(x,F) <  2sm \rbrace.$$
Due to the construction of $W$ we get $\hat{F} \subset  \hat{W}^\circ \subset U$. As $|\hat{F}|\le C_{u}k^2$, Lemma \ref{rig-lemma: weak rig}  for $\mu = 2k$  implies that there is a rotation $R \in SO(2)$ such that
$$\Vert \nabla y  - R\Vert^2_{L^2(\hat{F})}  \le C k^4  s^{-4} \Vert \dist(\nabla y,SO(2))\Vert^2_{L^2(\hat{F})} = C l^4 \gamma(\hat{F}),$$
where for shorthand we write $\gamma(\hat{F}) = \Vert \dist(\nabla y,SO(2))\Vert^2_{L^2(\hat{F})}$. As $\nabla y - R$ is harmonic in $\hat{F}$, the mean value property of harmonic functions for $r = sm$ and Jensen's inequality yield
\begin{align}\label{rig-eq: uniform estimate}
\begin{split}
|\nabla y(x) - R|^4 & \le \Big|\frac{1}{|B_r(x)|} \int_{B_r(x)} (\nabla y(t) - R) \, dt \Big|^4  \\
&  \le C\Big((sm)^{-2} \int_{\hat{F}} |\nabla y  -R|^2 \Big)^2   \le C l^{ 8} m^{-4} s^{-4} \gamma(\hat{F})^2
\end{split}
\end{align}
for all $x \in F$. Consequently, as $\hat{F}$ intersects at most nine squares $Q^{k}(p)$, $p \in I^{k}(\Omega^k) \setminus J$, by \eqref{rig-eq: weaklocR3} and $l = \frac{k}{s}$ we get $\Vert \nabla y  - R\Vert^2_{L^\infty(F)} \le C l^{ 4} m^{-2} s^{-2} \cdot  k \epsilon \le C l^{-4}\vartheta $ as well as
\begin{align*}
\begin{split}
\Vert \nabla y  - R\Vert^4_{L^4(F)}  \le C \vartheta l^{-4} \Vert \nabla y  - R\Vert^2_{L^2(\hat{F})} \le C  \vartheta \gamma(\hat{F}). 
\end{split}
\end{align*}
Proceeding in this way for every connected component $F$ of all $Q^{k}_i(p)$, $p \in I^{k}_i(\Omega^k)$, and noting that every $Q^s(q)$, $q \in I^s(\Omega^k)$, intersects at most four different associated enlarged sets $\hat{F}$  ($Q^s(q)$ can intersect more than one set if it lies at the boundary of some $Q^{k}_i(p)$) we obtain a function $\hat{R}_i:  W^\circ \to SO(2)$ with the desired properties \eqref{rig-eq: weaklocR1}. \eop

We now concern ourselves with the subatomistic regime. 

\begin{lemma}\label{rig-lemma: local estimate2}
Let $M \ge 0$, $\epsilon>0$ and $s \le k \le  \epsilon$.  Then for a fixed constant $C = C(M)>0$ we have the following: \\
For all $U \in {\cal V}^s_k$ with $U \subset \Omega^k$ and for all $y \in H^1(U)$ with $\gamma$ as defined in \eqref{rig-eq: weakloc gamma} and  $\Vert \nabla y\Vert_\infty \le M$ there is a set  $W \in {\cal V}^s_k$ with  $W\subset \Omega^{3k}$, $|W \setminus U|=0$, $|(U \setminus W) \cap \Omega^{3k}|\le C_{u}k\Vert W \Vert_* $ and
\begin{align}\label{rig-eq: locest2}
\Vert W \Vert_* \le \Vert U \Vert_* + C\epsilon^{-1} \gamma
\end{align}
as well as mappings $\hat{R}_i: \Omega^{3k} \to SO(2)$, $i=1,\ldots,4$, which are constant on $Q^{k}_i(p) \cap W$, $p \in I^{k}_i(\Omega^k)$, such that
\begin{align}\label{rig-eq: locest1}
\begin{split}
\Vert \nabla y  - \hat{R}_i\Vert^2_{L^2(W)}  \le C \gamma + C\epsilon \Vert U \Vert_*.
\end{split}
\end{align}
\end{lemma}

\Proof  Similarly as in \eqref{rig-eq: weaklocR3} we let $J \subset I^{k}(\Omega^k)$ such that
\begin{align}\label{rig-eq: locest3}
\epsilon{\cal H}^1(\partial U \cap Q^k(q)) + \Vert \dist(\nabla y,SO(2))\Vert^2_{L^2(Q^{k}(q) \cap U)} > c_*\epsilon k
\end{align}
for all $q \in J$. Define $W = \Omega^{3k} \cap \big(\big(U \setminus \bigcup_{p \in J} Q^{k}(q)\big) \cup \bigcup_{p \in J} \partial Q^{k}(q)\big)$ and note that the $\Vert W \Vert_* \le \Vert U \Vert_* + C\epsilon^{-1} \gamma$ for $c_* = c_*(h_*)>0$ sufficiently large.  Indeed, for the subset $J_1 \subset J$, for which \eqref{rig-eq: weaklocR3} holds, we argue as in the previous proof. Then with $J_2 = J \setminus J_1$ we note $\Vert W \Vert_\infty \le \Vert U \Vert_\infty + C\epsilon^{-1} \gamma + 2\sqrt{2} k \cdot \# J_2$ and $\Vert W \Vert_{\cal H} \le \Vert U \Vert_{\cal H} + C\epsilon^{-1} \gamma + 8 k \cdot \# J_2 - c_* k \cdot \# J_2$. This gives the desired result for $c_*$ large. Moreover, we get $W \in {\cal V}^s_k$ and $|(U \setminus W) \cap \Omega^{3k}|\le C_{u}k\Vert W \Vert_*$.

Consider some $\tilde{Q} := Q^k_i(q)$, $q \in I^k_i(\Omega^k)$. We extend $y$ from $\tilde{Q}\cap W$ to $\tilde{Q}$ by setting $\bar{v} = y$ on $W\cap \tilde{Q}$ and $\bar{v}(x) = x$ for all $x \in \tilde{Q} \setminus W$. Note that $\bar{v} \in SBV(\tilde{Q})$ with $J_{\bar{v}} = \partial W \cap \tilde{Q}$. By Theorem \ref{rig-th: cgp} we obtain a function $v \in H^1(\tilde{Q})$ such that by a rescaling argument 
$$\Vert \nabla \bar{v} - \nabla v\Vert_{L^p(\tilde{Q})} \le C k^{\frac{2}{p}-1} \Vert\nabla \bar{v}\Vert_\infty {\cal H}^{1}(J_{\bar{v}}\cap \tilde{Q}) \le CM k^{\frac{2}{p}-1} k^{1 - \frac{1}{p}}\beta^{\frac{1}{p}} \le CM \epsilon^{\frac{1}{p}} \beta^{\frac{1}{p}}$$
for $p < 2$, where $\beta = {\cal H}^1(\partial W \cap \tilde{Q})$. In the second step we used $\beta \le Ck$ by  \eqref{rig-eq: locest3} and applied $k \le \epsilon$ in the last step. Consequently, we obtain 
$$\Vert\dist(\nabla v,SO(2))\Vert^p_{L^p(\tilde{Q})} \le C\Vert\dist(\nabla \bar{v},SO(2))\Vert^p_{L^p(\tilde{Q})} + C\epsilon \beta.$$
Thus, since $\gamma(\tilde{Q}) := \Vert \dist(\nabla  \bar{v},SO(2)) \Vert^2_{L^2(\tilde{Q})} = \Vert \dist(\nabla y,SO(2)) \Vert^2_{L^2(\tilde{Q} \cap W)}$, the  rigidity estimate in Theorem \ref{rig-th: geo rig} yields a rotation $R \in SO(2)$ such that
\begin{align*}
\Vert \nabla v - R\Vert^p_{L^p(\tilde{Q})} & \le C\Vert\dist(\nabla v,SO(2))\Vert^p_{L^p(\tilde{Q})} \le C |\tilde{Q}|^{1-\frac{p}{2}}\gamma(\tilde{Q})^{\frac{p}{2}} + C\epsilon\beta\\
& \le C \epsilon^{2-p}\gamma(\tilde{Q})^{\frac{p}{2}-1} \gamma(\tilde{Q}) +C\epsilon\beta\le C \epsilon^{2 - p} \epsilon^{p-2}\gamma(\tilde{Q}) + C\epsilon\beta \\
& \le  C\gamma(\tilde{Q}) + C\epsilon\beta.
\end{align*}
In the second step we used H\"older's inequality and we applied \eqref{rig-eq: locest3} in the  fourth step. This implies $\Vert \nabla y - R\Vert^p_{L^p(W \cap \tilde{Q})} \le \Vert \nabla \bar{v} - R\Vert^p_{L^p(\tilde{Q})} \le C\gamma(\tilde{Q}) + C\epsilon\beta$ and proceeding in this way for all $Q^{k}_i(q)$, $q \in I^{k}_i(\Omega^k)$, we obtain a function $\hat{R}_i: \Omega^{3k} \to SO(2)$ such that  for a constant $C=C(h_*)$
$$\Vert \nabla y - \hat{R}_i\Vert^p_{L^p(W)} \le C \gamma + C\epsilon \Vert U \Vert_*,$$
where $\hat{R}_i$ is constant on $Q^{k}_i(p) \cap W$, $p \in I^{k}_i(\Omega^k)$. Finally, by $\Vert \nabla y \Vert_\infty \le M$ we derive 
\begin{align*}
\Vert \nabla y - \hat{R}_i\Vert^2_{L^2(W)}  & \le (M + \sqrt{2})^{2-p}\Vert \nabla y - \hat{R}_i\Vert^p_{L^p(W)}  \le C \gamma + C\epsilon \Vert U \Vert_*,
\end{align*}
as desired. \eop

Given a deformation $y: F \to \R^2$ for $F \subset \R^2$ and a rotation $R \in SO(2)$ we define the displacement field $u_R := R^T y - \id$, where $\id$ denotes the identity function. We introduce the linear elastic strain by 
$$
\bar{e}_R(\nabla y) :=    e(\nabla u_R) =  \frac{R^T \nabla y + (\nabla y)^T R}{2} -\Id,
$$
where $\Id$ denotes the identity matrix. For a general function $\hat{R}:F \to SO(2)$ we then define for shorthand $\alpha_{\hat{R}}(F) = \Vert \bar{e}_{\hat{R}}(\nabla y)\Vert^2_{L^2(F)}$.  Applying the linearization formula
\begin{align}\label{rig-eq: linearization}
\dist(G,SO(2)) = |\bar{e}_R(G)| + O(|G - R|^2)
\end{align}
for $R \in SO(2)$ and $G \in \R^{2 \times 2}$ we get 
\begin{align}\label{rig-eq: linearization2}
\alpha_{\hat{R}}(F) = \int_{F} |\bar{e}_{\hat{R}}(\nabla y)|^2 \le C_u \int_{F} \dist^2(\nabla y,SO(2))
 + C_u \int_{F} |\nabla y - \hat{R}|^4 .
\end{align}
Here we already see that it suffices to establish a rigidity estimate of fourth order as in Lemma \ref{rig-lemma: weaklocR} in order to bound the symmetric part of the gradient. One of the main ideas in the following will be to choose $l=l(s,\epsilon,m)$ in \eqref{rig-eq: weaklocR1} such that $\vartheta \le 1$ which will imply $\alpha_{\hat{R}}(W)  \le C_u \gamma$.

\subsection{Estimates in terms of the H$^1$-norm}\label{rig-sec: subsub,  est-h1}

We now show that not only the distance of the derivative from a rigid motion can be controlled as derived  in \eqref{rig-eq: weaklocR1} and \eqref{rig-eq: locest1}, respectively,  but also the distance of  the function itself. Once we have such estimates we will be in a position to  `heal' cracks (see Section \ref{rig-sec: sub, local} below). After the modification of the deformation $\nu = sd$ will stand for the minimal distance of two different cracks, where $d$ represents the corresponding increase factor. It will turn out that the least crack length  will be given by $\lambda = \nu m^{-1}$. Moreover, $k = \lambda m^{-1}$ will denote the size of the cell on which we apply Theorem \ref{rig-th: derive prop}. Define 
$$S_i := \bigcup\nolimits_{p \in I^k_i(\Omega^{3k})} Q_i^{\frac{5}{8}k}(p)$$
and note that $\Omega^{5k}  \subset \bigcup^4_{i=1} S_i$. Recall   \eqref{rig-eq: symmeas},  \eqref{rig-eq: no holes}, \eqref{rig-eq: no holes2} and the definition  $\hat{m} = C_2(m, h_*)$  (see below \eqref{rig-eq: vartheta def}).  We will proceed in two steps first deriving an estimate for the total variation of the distributional derivative (cf. Corollary \ref{rig-cor: kornpoin}) and then employing Theorem \ref{rig-theo: korn}. For shorthand we will write  $\gamma(F) = \Vert\dist(\nabla y, SO(2))\Vert^2_{L^2(F)}$, $\delta_p(F) = \sum\nolimits^4_{i=1}\Vert \nabla y  - \hat{R}_i\Vert^p_{L^p(F)}$ for $p=2,4$ and subsets $F \subset W$.

\begin{lemma}\label{rig-lemma: weaklocA2}
Let $k > s  > 0$,  $\epsilon >0$ such that $l:=\frac{k}{s} = dm^{-2}$ for $m^{-1},d \in \N$ with $m^{-1}, d \gg 1$. Let $\lambda = sdm^{-1} = km$. Then for constants $C,c>0$ we have the following: \\
For all $W \in  {\cal V}^{sm}_{(s,3k)}$ with $W \subset \Omega^{3k}$ and for all $y \in H^1(W)$ with 
\begin{align*}
\gamma:= \Vert \dist(\nabla y,SO(2)\Vert^2_{L^2(W)}, \ \ \ \delta_4 := \sum\nolimits^4_{i=1}\Vert \nabla y  - \hat{R}_i\Vert^4_{L^4(W)}  
\end{align*}
for mappings $\hat{R}_i:  W^\circ\to SO(2)$, $i=1,\ldots,4$, which are constant on the connected components of $Q^{k}_i(p) \cap  W^\circ$, $p \in I^{k}_i(\Omega^{3k})$, we obtain:\\
We find sets $U \in {\cal V}^{s\hat{m}}_{70k}$,  $U_Q \in  {\cal V}^{s\hat{m}}$ with $U \subset  U_Q \subset \Omega^{5k}$, $| U_Q \setminus W|=0$ and $|(W \setminus U) \cap \Omega^{5k}| \le C_{u}k\Vert U \Vert_*$ such that
\begin{align}\label{rig-eq: weaklocA1}
\begin{split}
\Vert U \Vert_* \le (1+C_{u}m)\Vert W \Vert_* +  C\epsilon^{-1}(\gamma+ \delta_4)
\end{split}
\end{align}
as well as 
\begin{align}\label{rig-eq: weaklocA2.13}
\begin{split}
 |Q^{\lambda}(p) \cap  U_Q| \ge  c m\lambda^2 \ \ \  \text{ for all } p \in J( U_Q),
\end{split}
\end{align}
where $J( U_Q) := \lbrace p \in I^{\lambda}(\Omega^{3k}): Q^{\lambda}(p) \cap  U_Q \neq \emptyset \rbrace$.  \\ Moreover, letting $U_J = \bigcup_{p \in J( U_Q)} \overline{Q^{\lambda}(p)} $, for $i=1,\ldots,4$ we find extensions  $\bar{y}_i \in SBV^2(U_J\cap S_i, \R^2)$ with $\bar{y}_i = y$ on $ U_Q \cap S_i$ such that for all $\tilde{Q} := Q_j^{ 3\lambda}(p) \cap U_J$, $p \in I^{\lambda}_j(\Omega^{3k})$, $j=1,\ldots,4$, with $\tilde{Q} \subset S_i$ we have   that  $R_i = \hat{R}_i|_{W^\circ \cap \tilde{Q}}$ is constant on $W^\circ \cap \tilde{Q}$ and
\begin{align}\label{rig-eq: weaklocA17}
\begin{split}
(|E ( R^T_i\bar{y}_i - \id)|(\tilde{Q}))^2  &\le Ck^2 C_m\min\Big\{ \epsilon k, \gamma(W \cap Q_i^{2k}(q)) \\ & \quad \quad \quad \quad + \delta_4(W \cap Q_i^{2k}(q)) + \epsilon|\partial W \cap Q_i^{2k}(q)|_{\cal H}\Big\},
\end{split}
\end{align}
where  $q \in I^k_i(\Omega^{3k})$ such that $\tilde{Q} \subset Q_i^{k}(q)$. 
\end{lemma}

\Proof  Similarly as in the previous proof we let $J \subset I^{3k}(\Omega^{3k})$ such that
\begin{align}\label{rig-eq: weaklocA6}
\begin{split}
\epsilon{\cal H}^1(Q^{3k}(p) \cap  \partial W)&+\Vert \dist(\nabla y,SO(2))\Vert^2_{L^2(Q^{3k}(p) \cap W)}\\
& + \sum\nolimits^4_{i=1}\Vert \nabla y  - \hat{R}_i\Vert^4_{L^4(Q^{3k}(p) \cap W)} > c_*\epsilon  k
\end{split}
\end{align}
for all $p \in J$. Define $\hat{W} = \big(W \setminus \bigcup_{p \in J} Q^{3k}(p)\big) \cup \bigcup_{p \in J} \partial Q^{3k}(p)$ and note that choosing  $c_*$ sufficiently large and arguing as in the previous proof
\begin{align}\label{rig-eq: weaklocA7}
\Vert \hat{W} \Vert_* \le  \Vert W \Vert_* + C\epsilon^{-1} ( \gamma + \delta_4), 
\end{align}
$\hat{W} \in {\cal V}^{sm}_{(s,3k)}$ as well as $|(W \setminus \hat{W})  \cap \Omega^{5k}| \le C_{u}k\Vert \hat{W} \Vert_*$. We now subsequently construct sets $\hat{U}_1 \supset \ldots \supset \hat{U}_4$  (the inclusions hold up to sets of negligible measure) by application of Theorem \ref{rig-th: derive prop} on connected components of $\hat{W}$ (Step (I)). Afterwards, since in Theorem \ref{rig-th: derive prop} the trace estimate cannot be derived for components near the boundary, we will further modify the sets in a neighborhood of large boundary components (Step (II)). A final modification procedure will then assure property \eqref{rig-eq: weaklocA2.13} (Step (III)).

\vspace{.5cm}
(I) Begin with $i=1$ and fix $q \in I^{k}_1(\Omega^{3k})$. Consider a connected component $F$ of $Q^{k}_1(q) \cap  \hat{W}^\circ$. As $\hat{R}_1 = R$ is constant on $F$ we obtain $\alpha_R(F) \le C(\gamma(F) + \delta_4(F))$ by \eqref{rig-eq: linearization2}.  Define $Q_\mu := Q^{k}_1(q)$  and recall \eqref{rig-eq: no holes}. Passing to the closure of $F$ (not relabeled)  we can regard $F$ as an element of ${\cal V}^{sm}$ with respect to $Q_\mu$ (recall \eqref{rig-eq: calV-s def2}), where one component is given by $X = Q_\mu \setminus H(F) \in {\cal U}^{sm}$.  (Observe, however, that $Q_\mu \setminus H(F)$ may intersect several components of $\hat{W}$.) We apply Theorem \ref{rig-th: derive prop} on $F \subset Q_\mu$ for $\eps=\epsilon$, $\sigma=m$ to obtain a set $G \in  {\cal W}^{s\hat{m}}$  with $|G \setminus F|=0$ and  
\begin{align}\label{rig-eq: new101}
\epsilon \Vert G \Vert_* + \alpha_R(G) \le (1 + C_{u}m)(\epsilon \Vert F \Vert_* + \alpha_R(F)).
\end{align}
(Recall that the sum in $\Vert F \Vert_*$ runs only over the boundary components having  empty intersection with $\partial Q_\mu$.)  Moreover, similarly as before we have
\begin{align}\label{rig-eq: F hat}
|F \setminus G| \le C_{u} k \Vert G \Vert_*
\end{align} 
and  using \eqref{rig-eq: D1}, \eqref{rig-eq: newnew}(i),(ii) for all  $\Gamma_t(G) \in {\cal T}(G):=\lbrace \Gamma_t: N^{\tau_t}(\partial R_t) \subset H(G) \rbrace$
\begin{align}\label{rig-eq: weaklocA3} 
\int\nolimits_{\Theta_t(G)} |[\bar{y}_1](x))|^2 \,d{\cal H}^1(x) \le CC_m  \epsilon \vert \Gamma_t(G) \vert^2_\infty,
\end{align}
where $\bar{y}_1$ is the extension (cf. \eqref{rig-eq: extend def_new})
\begin{align}\label{rig-eq: construction}
\bar{y}_1(x) = \begin{cases}  y & x \in \hat{W},\\ R\,(\Id + A_t)\, x + R\,c_t & x \in X_t \ \text{ for } \Gamma_t(G) \in {\cal T}(G).  \end{cases}
\end{align}
Here recall  that $\partial R_t$ are the rectangles  given by \eqref{rig-eq: newnew}  as well as $\tau_t = \bar{\tau}|\partial R_t|_\infty \ll |\partial R_t|_\infty$.

 We proceed in this way for every connected component $(F_j)_j$ of all $Q^{k}_1(q)$, $q \in I^{k}_1(\Omega^{3k})$ and define $\hat{U}_1 = (\hat{W} \setminus \bigcup_j F_j) \cup \bigcup_j  G_j \in {\cal V}^{s\hat{m}} $. (Observe that one may have $H(F_{j_1}) \subset H(F_{j_2})$. In this case the above arguments can be omitted for $F_{j_1}$.)  By ${\cal G}$ we denote the set of boundary components $\Gamma(\hat{U}_1)$ which do not coincide with some  $\Gamma_t(G_j)$. Note that  by \eqref{rig-eq: linearization2} and \eqref{rig-eq: weaklocA7} 
\begin{align}\label{rig-eq: weaklocA4}
\begin{split}
\epsilon\Vert \hat{U}_1 \Vert_* \le \epsilon\Vert \hat{U}_1 \Vert_* + \alpha_{\hat{R}_1}(\hat{U}_1) &\le (1 + C_{u}m)(\epsilon\Vert \hat{W} \Vert_* + \alpha_{\hat{R}_1}(\hat{W})) \\
& \le (1 + C_{u}m)\epsilon\Vert W \Vert_* + C(\gamma+ \delta_4).
\end{split}
\end{align}
The second step follows as by construction for each  $\Gamma(\hat{U}_1) \in {\cal G}$ there is a $\Gamma(\hat{W}) = \partial X$ such that $\Gamma(\hat{U}_1) \subset \overline{X}$ (recall Remark \ref{rig_rem: connect}(i)). By \eqref{rig-eq: F hat}  we also get $|\hat{W} \setminus \hat{U}_1| \le C_{u}k\Vert \hat{U}_1 \Vert_*$.  Moreover, by Remark \ref{rig_rem: connect}(ii) we can replace the components of $G_j \in {\cal V}^{s\hat{m}}$ by rectangles such that the resulting set $G_j'$ lies in ${\cal V}^{s\hat{m}}_{\rm con}$. Recall that the (rectangular) components of $G_j'$ satisfy $\max_{i=1,2} |\pi_i \Gamma(G_j')| \le 2k$. 

Then we define $\hat{U}''_1 := (\hat{W} \setminus \bigcup_j F_j) \cup \bigcup_j G_j' \in {\cal V}^{s\hat{m}}$. We now apply Lemma \ref{rig-lemma: modifica}(ii) for $\nu=0$, $(Z_j)_j$ the rectangular components of $(G_j')_j$ and $V'$ the set whose boundary components are given by the elements of ${\cal G}$. We obtain a  set $\hat{U}'_1 \in {\cal V}^{s\hat{m}}_{5k}$ with $\Vert \hat{U}'_1 \Vert_* \le \Vert \hat{U}_1 \Vert_*$ and $|\hat{U}''_1 \setminus \hat{U}'_1| \le C_u k \Vert \hat{U}'_1 \Vert_*$. (Strictly speaking, we need to pass from ${\cal V}^{s\hat{m}}$ to ${\cal V}^{s\hat{m}/2}$, but do not include it in the notation for convenience.) Likewise we observe $|\hat{U}'_1 \setminus \hat{W}| = 0$ and $|\hat{W} \setminus \hat{U}'_1| \le C_{u}k\Vert \hat{U}'_1 \Vert_*$.   Additionally, we apply Lemma \ref{rig-lemma: small bc}(ii) and get a (not relabeled) set $\hat{U}'_1 \in {\cal V}^{s\hat{m}}_{6k}$ such that \eqref{rig-eq: small bc}  and \eqref{rig-eq: weaklocA4} hold.  As in the proof of Lemma \ref{rig-lemma: weaklocR} this implies $\hat{U}'_1 \in {\cal V}^{s\hat{m}}_{ (s,6k)}$ since $\hat{W} \in {\cal V}^{sm}_{(s,3k)}$, i.e. the least length of components is bounded from below by $s$.

 In the following, by a slight abuse of notation, we say that a component $\Gamma_t(\hat{U}'_1)$, which coincides with some $\partial X_t = \Gamma_t(G')$ for some component $G'$, satisfies \eqref{rig-eq: weaklocA3} if all corresponding $(\Gamma_{t_s}(G))_s$ with $\Gamma_{t_s}(G) \subset \overline{X_t}$ satisfy \eqref{rig-eq: weaklocA3}. 
It is not hard to see that \eqref{rig-eq: weaklocA3} is satisfied for all boundary components with  (recall \eqref{rig-eq: no holes2})
\begin{align*}
\Gamma_t(\hat{U}'_1) \cap S_1 \neq \emptyset, \ \ \ \vert \Gamma_t(\hat{U}'_1) \vert_\infty \le \tfrac{k}{8}, \ \ \  N_*(\Gamma_t(\hat{U}_1)) \subset H^{\frac{k}{8}}(\hat{U}'_1),
\end{align*}
 where $N_*(\Gamma_t(\hat{U}_1)) = \lbrace x: \dist(x,\Gamma_t(\hat{U}'_1)) \le \bar{C}\hat{m} |\Gamma_t(\hat{U}'_1)|_\infty \rbrace$ for some large constant $\bar{C} = \bar{C}(h_*) >0$. Indeed, assume that  there is some  $ \Gamma_s = \Gamma_{t_s}(G) \subset Q^{k}_1(q)$ such that for the corresponding rectangle $R_s$ one has $N^{\tau_s}(\partial R_s) \not\subset H(G)$ although the corresponding $\Gamma_t(G') = \partial X_t$ fulfills the above three properties.  First, we observe $R_s \subset X_t$ by Remark \ref{rig_rem: connect}(ii) and thus $R_s \subset Q^{\frac{3}{4}k}_1(q)$. By \eqref{rig-eq: newnew}(i) we get $|\partial R_s|_\infty \le C|\Gamma_s|_\infty$.  Consequently, since $\tau_s \ll \frac{1}{C}|\partial R_s|_\infty$ (recall the assumption in Theorem \eqref{rig-th: derive prop})  we have $N^{\tau_s}(\partial R_s) \subset Q^{\frac{7}{8}k}_1(q)$. Since by assumption $N^{\tau_s}(\partial R_s) \not\subset H( G)$, this would imply $\vert \partial H( G) \cap Q^k_1(q)\vert_\infty > \frac{k}{8}$.  
 
Consequently, there is a chain of components $(\Gamma_{t_i}(\hat{U}'_1))^n_{i=1} = (\partial X_{t_i}(\hat{U}'_1))^n_{i=1}$ such that $\Gamma_{t_1}(\hat{U}'_1) \cap \partial Q_\mu \neq \emptyset$, $ X_{t_n}(\hat{U}'_1) \cap N^{\tau_s}(\partial R_s) \neq \emptyset$ and $\Gamma_{t_{i-1}}(\hat{U}'_1) \cap \Gamma_{t_i}(\hat{U}'_1) \neq \emptyset$. Thus,  by \eqref{rig-eq: small bc} there is  one $\Gamma_*(\hat{U}'_1)$ with $|\Gamma_*(\hat{U}'_1)|_\infty > \frac{k}{8}$ such that $N^{\tau_s}(\partial R_s) \cap X_*(\hat{U}'_1) \neq \emptyset$.  Recalling that $R_s \subset X_t$ and $\tau_s < \bar{C} \hat{m}|\Gamma_t(\hat{U}'_1)|_\infty$ for $\bar{C}$ sufficiently large by Remark \ref{rig-rem: z}(iii) we find $N_*(\Gamma_t(\hat{U}_1)) \cap  X_*(U_1') \neq \emptyset$. This, however, is a contradiction to $N_*(\Gamma_t(\hat{U}_1)) \subset H^{\frac{k}{8}}(\hat{U}'_1)$.

We now iteratively repeat the above construction for $i=2,3,4$  for $\hat{U}'_{i-1}$ instead of $\hat{W}$ and obtain extensions $\bar{y}_2,\bar{y}_3, \bar{y}_4$  as well as  $(\hat{U}_i)^4_{i=1}$ and sets $\hat{U}'_4 \subset \ldots \subset \hat{U}'_1 \subset \hat{W}$ (the inclusions hold up to a set of negligible measure) with $\hat{U}'_4 \in {\cal V}^{s\hat{m}}_{(s,15k)}$ such that \eqref{rig-eq: weaklocA4} holds for a possibly larger constant replacing $\hat{U}_1$ by $\hat{U}_4$.  We briefly note that the sets are elements of ${\cal V}^{s\hat{m}}$ due to Remark \ref{rig-rem: z}(iv). Moreover, for $i=1,\ldots,4$, \eqref{rig-eq: weaklocA3} is satisfied for $\bar{y}_i$ and all boundary components $\Gamma_t(\hat{U}'_i)$ with
$\Gamma_t(\hat{U}'_i) \cap S_i \neq \emptyset$, $\vert \Gamma_t(\hat{U}'_i) \vert_\infty \le \frac{k}{8}$ and $N_*(\Gamma_t(\hat{U}'_i)) \subset H^{\frac{k}{8}}(\hat{U}'_{i})$.

 For later we also observe that due to the local nature of the  modification process and \eqref{rig-eq: new101} we get 
\begin{align}\label{rig-eq: new104}
\begin{split}
|\partial \hat{U}_i \cap Q^k_i(q)|_{\cal H} &\le C|\partial \hat{W} \cap Q^{2k}_i(q)|_{\cal H} \\
& \ \ \  +C\epsilon^{-1}\big(\gamma(\hat{W} \cap Q^{2k}_i(q)) + \delta_4(\hat{W} \cap Q^{2k}_i(q))\big).
\end{split}
\end{align}
 Although the inclusions for $(\hat{U}'_i)^4_{i=1}$ only hold up to segments, we observe that the sets are `nested' concerning small boundary components in the following sense: Letting $\hat{U}^*_i = \hat{U}'_i \cap \overline{(H^{\frac{k}{8}}(\hat{U}'_i))^\circ}$ we obtain 
\begin{align}\label{rig-eq: new 513}
\hat{U}^*_4 \subset \ldots \subset \hat{U}^*_1.
\end{align}
Indeed, assume e.g. there was a component $X(\hat{U}^*_1)$ and components $X_1, \ldots, X_n$ of $\hat{U}^*_2$ with $X(\hat{U}^*_1) \subset  \bigcup^n_{j=1} \overline{X_j}$ and $\bigcup^n_{j=1} \partial X_j \cap X(\hat{U}^*_1) \neq \emptyset$. Then by construction of the sets we clearly find some $X_i$ with $\partial X_i \cap X(\hat{U}^*_1) \neq \emptyset$, $|X(\hat{U}^*_1) \setminus X_i| >0$ and $|\partial X_i|_\infty \le \frac{k}{8}$. This, however, together with \eqref{rig-eq: small bc} gives a contradiction to $X(\hat{U}^*_1) \subset  \bigcup^n_{j=1} \overline{X_j}$. In particular, \eqref{rig-eq: new 513} implies $H^{\frac{k}{8}}(\hat{U}'_4) \subset \ldots \subset H^{\frac{k}{8}}(\hat{U}'_1)$ up to sets of negligible measure and thus  for $i=1,\ldots,4$, \eqref{rig-eq: weaklocA3} is satisfied for $\bar{y}_i$ and all boundary components $\Gamma_t(\hat{U}'_i)$ with
$\Gamma_t(\hat{U}'_i) \cap S_i \neq \emptyset$, $\vert \Gamma_t(\hat{U}'_i) \vert_\infty \le \frac{k}{8}$, $N_*(\Gamma_t(\hat{U}'_i)) \subset H^{\frac{k}{8}}(\hat{U}'_{ 4})$.  We want to remove the third condition. For that reason, we subtract neighborhoods of large boundary components as follows.

\vspace{.5cm}

(II) Let $U^* = H^{\frac{k}{8}}(\hat{U}'_4)$ and let $\Gamma_1 (U^*), \ldots, \Gamma_n (U^*)$ be the boundary components. For $\Gamma_j (U^*)$ let $M(\Gamma_j)$ be the smallest rectangle in ${\cal U}^{s\hat{m}}$ satisfying $M(\Gamma_j) \supset \lbrace x\in\R^2: \dist_\infty(x, X_j) \le \bar{C}k\hat{m} \rbrace$ for the constant $\bar{C}>0$ introduced above, where $X_j$ denotes the component corresponding component to $\Gamma_j(U^*)$. Clearly,  using the fact that $\bar{C}\hat{m} \le m$ (see \eqref{rig-eq: vartheta def}) one has    $|\pi_i \partial M(\Gamma_j)| \le |\pi_i \Gamma_j(U^*)| + C_u m |\Gamma_j(U^*)|_\infty \le 31k$ for $i=1,2$.  As the components $(X_j)_j$ are pairwise disjoint and connected, we obtain $Z(\Gamma_{j_1}) \setminus Z(\Gamma_{j_2})$ is connected for all $1 \le j_1, j_2 \le n$, where $Z(\Gamma_j)$ denotes the smallest rectangle containing $X_j$. Consequently, since the neighborhoods $M(\Gamma_j) \setminus Z(\Gamma_j)$ all have the same thickness $\sim \bar{C}k\hat{m}$, we get that $M(\Gamma_{j_1}) \setminus M(\Gamma_{j_2})$ is connected for all $1 \le j_1, j_2 \le n$.

Then by Lemma \ref{rig-lemma: modifica}(ii) applied on $V=U^*$, $V' = \Omega^{5k} \setminus \bigcup_{|\Gamma_j(\hat{U}'_i)|\le \frac{k}{8}} X_j(\hat{U}'_i)$ we obtain sets $\tilde{U}_i$ with $|(\hat{U}'_i \setminus \bigcup^n_{j=1} M(\Gamma_j)) \setminus \tilde{U}_i| \le C_u k \Vert V' \Vert_*$. In particular, we set $\tilde{U} = \tilde{U}_4$ and observe that $\tilde{U} \in {\cal V}^{s\hat{m}}_{32k}$. Moreover, we obtain $\Vert \tilde{U} \Vert_* \le (1+C_u m) \Vert V\Vert_* + \Vert V'\Vert_*$. As $\hat{U}'_4$ satisfies \eqref{rig-eq: small bc}, we derive $(\partial V \cap \partial V') \cap (\Omega^{5k})^\circ =\emptyset$ and therefore $\Vert \tilde{U} \Vert_* \le (1+C_u m) \Vert \hat{U}'_4\Vert_*$, i.e.  \eqref{rig-eq: weaklocA4} holds replacing $\hat{U}_1$ by $\tilde{U}$ (possibly for a larger constant). Applying Lemma \ref{rig-lemma: small bc}(ii) we get (not relabeled) sets $\tilde{U}_i \in {\cal V}^{s\hat{m}}_{33k}$ satisfying \eqref{rig-eq: small bc}. For later we note that $\tilde{U}_4 \subset \ldots \subset \tilde{U}_1$  up to sets of negligible measure. This follows from \eqref{rig-eq: new 513} and the fact that in Lemma \ref{rig-lemma: modifica}(ii) the components of $V'$ are replaced by corresponding rectangles.  Arguing as in \eqref{rig-eq: new 513} we also find 
\begin{align}\label{rig-eq: new 513**}
\text{$\tilde{U}^*_4 \subset \ldots \subset \hat{U}^*_1$, \ \ where \ \  $\tilde{U}^*_i = \tilde{U}_i \cap \overline{(H^{\frac{k}{8}}(\tilde{U}_i))^\circ}$} 
\end{align}
In particular, this also implies $H^{\frac{k}{8}}(\tilde{U}_4) \subset \ldots \subset H^{\frac{k}{8}}(\tilde{U}_1)$ up to sets of negligible measure.

We now see that for $i=1,\ldots,4$, \eqref{rig-eq: weaklocA3} holds  for $\bar{y}_i$ for all components satisfying 
\begin{align}\label{rig-eq: weaklocA11}
\Gamma_t(\tilde{U}_i) \cap S_i \neq \emptyset, \ \ \  \vert \Gamma_t(\tilde{U}_i ) \vert_\infty \le \tfrac{1}{8}k.
\end{align}
(Strictly speaking \eqref{rig-eq: weaklocA3} holds for the corresponding components of $\hat{U}_i$.)  In fact, since $\bar{C}\hat{m}\vert \Gamma_t(\hat{U}'_i) \vert_\infty \le \frac{k}{8}\bar{C}\hat{m}$ for $\vert \Gamma_t(\hat{U}'_i) \vert_\infty \le \frac{k}{8}$, due to the construction of $\tilde{U}_i$ components with $\vert \Gamma_t(\hat{U}'_i) \vert_\infty \le \frac{k}{8}$ and $N_*(\Gamma_t(\hat{U}'_i)) \not\subset H^{\frac{k}{8}}(\hat{U}'_4)$ are `combined' with a boundary component of $ \hat{U}'_4$ which is larger than $\frac{k}{8}$. 

 We apply Lemma \ref{rig-lemma: small bc}(i) to obtain a (not relabeled) set $\tilde{U} \in {\cal V}^{s\hat{m}}_{33k}$ satisfying \eqref{rig-eq: small bc2}. For each $\Gamma_t(\tilde{U})$, $t=1,\ldots,n$,  let $N_1(\Gamma_t)$, $N_2(\Gamma_t)$ be the smallest rectangles in ${\cal U}^{s\hat{m}}$ satisfying  
\begin{align*}
& N_1(\Gamma_t) \supset \lbrace x \in \R^2: \dist_\infty(x , X_t) \le  \min \lbrace Bm |\Gamma_t(\tilde{U})|_\infty, 2\lambda \rbrace \rbrace ,\\
& N_2(\Gamma_t) \supset \lbrace x \in \R^2: \dist_\infty(x , X_t) \le Bm  \, \min\lbrace|\Gamma_t(\tilde{U})|_\infty,\lambda\rbrace \rbrace
\end{align*}
for some $B>0$  (independent of $h_*$) and $\lambda = km$, where $X_t$ is the component corresponding to $\Gamma_t(\tilde{U})$.  It is not restrictive to assume that 
\begin{align}\label{rig-eq: isoper}
{\cal H}^1 \big(N_2(\Gamma_t) \cap (\partial \tilde{U}  \setminus (\Gamma_t(\tilde{U}) \cup \partial H^{\frac{k}{8}}(\tilde{U}))\big) \le CBm \,  \min\lbrace|\Gamma_t(\tilde{U})|_\infty,\lambda\rbrace 
\end{align}
 for all  $\Gamma_t(\tilde{U})$ with $|\Gamma_t(\tilde{U})|_\infty \le \frac{k}{8}$. Indeed, otherwise we replace $\tilde{U}$ by  $\tilde{U}' := \big(\tilde{U} \setminus N^*_2(\Gamma_t)\big)\cup \partial N^*_2(\Gamma_t)$, where $N^*_2(\Gamma_t) = (N_2(\Gamma_t) \cap H^{\frac{k}{8}}(\tilde{U}))^\circ$, and arguing similarly as in \eqref{rig-eq: weaklocA6}   and Lemma \ref{rig-lemma: infty} we  get $\Vert \tilde{U}' \Vert_* \le \Vert \tilde{U} \Vert_*$.   Let $(X_{t'})_{t'}$, $X_{t'} \neq X_t$, be the components of $\tilde{U}$ having nonempty intersection with $N^*_2(\Gamma_t)$. Clearly, we have $|\partial X_{t'}|_\infty \le \frac{k}{8}$. We define $T = \overline{N^*_2(\Gamma_t)} \cup \bigcup_{t'} \overline{X_{t'}}$ and modify $\tilde{U}'$ on a set of measure zero by letting $\tilde{U}'' = (\tilde{U}' \setminus T) \cup \partial T$.  Arguing similarly as in the proof of Lemma \ref{rig-lemma: small bc} we find $\tilde{U}'' \in {\cal V}^{s\hat{m}}_{33k}$ and $\Vert \tilde{U}'' \Vert_* \le \Vert \tilde{U} \Vert_*$. Then by Lemma \ref{rig-lemma: small bc}(i) we find a (not relabeled) set $\tilde{U}''$ which additionally satisfies \eqref{rig-eq: small bc2}.  We continue with this iterative modification process until \eqref{rig-eq: isoper} is satisfied for all components smaller than $\frac{k}{8}$.   
 Finally, by  Lemma \ref{rig-lemma: small bc}(ii) we obtain a (not relabeled) set $\tilde{U}'' \in {\cal V}^{s\hat{m}}_{34k}$ satisfying \eqref{rig-eq: small bc}. Noting that during the modification procedure  components larger than $\frac{k}{8}$ do not become smaller than $\frac{k}{8}$  we also find $H^{\frac{k}{8}}(\tilde{U}'') \subset H^{\frac{k}{8}}(\tilde{U})$.  For convenience the set will still be denoted by $\tilde{U}$ in the following.

\vspace{.5cm}
 
(III) We now finally construct  the sets $U_Q$ and $U$. For each $t=1,\ldots,n$  define the rectangle 
\begin{align}\label{rig-eq: new802}
 Z_t = \bigcup\nolimits_{p \in I^\lambda(N_{1}(\Gamma_t))} \overline{Q^{\lambda}(p)}.
\end{align}
We find  $Z_t \subset N_{1}(\Gamma_t)$ and for sufficiently small components one has $Z_t =  \emptyset$. Choosing $B$ sufficiently large we get $X_t \subset Z_t$ if $|\partial X_t|_\infty > \tfrac{k}{8}$.  Rearrange the components in a way that $Z_t = \emptyset$ for $t>n'$.  This implies
 \begin{align}\label{rig-eq: new801}
\Omega^{5k} \setminus H^{\frac{k}{8}}(\tilde{U}) \subset \bigcup\nolimits^{n'}_{t=1}  Z_t.
\end{align}
Let $Y_t \in {\cal U}^{s\hat{m}}$ be the smallest rectangle containing $Z_t \cup X_t$. By the definition of $N_1(\Gamma_t)$ and $Z_t$ we obtain  
\begin{align}\label{rig-eq: Ycrack}
 |\pi_i \partial Y_t| = |\pi_i \partial (\overline{Z_t \cup X_t})| \le  |\pi_i \Gamma_t(\tilde{U})| + C_{u}m| \Gamma_t(\tilde{U})|_\infty , \ \ i=1,2
\end{align}
 for some  $C_u=C_u(B)$ large enough. As $(X_t)_t$ are pairwise disjoint and connected, it is elementary to see that $Z_{t_1} \setminus Z_{t_2}$ or $Z_{t_2} \setminus Z_{t_1}$ is connected for all $1 \le t_1, t_2 \le n'$. In fact, assume there were $t_1 \neq t_2$ such that $\overline{\pi_1 Z_{t_2}} \subset \pi_1 Z_{t_1}$ and $\overline{\pi_2 Z_{t_1}} \subset \pi_2 Z_{t_2}$. Then due to the definition of the neighborhoods we find $\overline{\pi_1  X_{t_2}} \subset \pi_1 X_{t_1}$ and $\overline{\pi_2  X_{t_1}} \subset \pi_2 X_{t_2}$. This, however, implies $X_{t_1} \cap X_{t_2} \neq \emptyset$ and yields a contradiction.  A similar argument yields that $Y_{t_1} \setminus Y_{t_2}$ or $Y_{t_2} \setminus Y_{t_1}$ is connected for all $1 \le t_1, t_2 \le n'$.

Define $U_Q' = \tilde{U} \setminus \bigcup^{n'}_{j=1} Z_j$ and  let $\hat{J} \subset I^\lambda(\Omega^{3k})$ such that (cf. also \eqref{rig-eq: weaklocA6}) 
\begin{align}\label{rig-eq: weaklocA6*}
{\cal H}^1(Q^\lambda(p) \cap \partial  U_Q') > c_*\lambda
\end{align}
for all $p \in \hat{J}$. Then let $U_Q = (\Omega^{5k} \cap U_Q') \setminus \bigcup\nolimits_{p \in \hat{J}} Q^\lambda(p)$. Observe that possibly $U_Q \notin {\cal V}^{s\hat{m}}_{\rm con}$. Therefore, we now define a set $U \subset U_Q$ with connected boundary components.

By Lemma \ref{rig-lemma: modifica}(ii) for $V= \Omega^{5k} \setminus \bigcup^{n'}_{t=1} X_t$, $V' = \Omega^{5k} \setminus \bigcup^{n}_{t=n' +1} X_t$  we obtain a set $U'$ with $|(\tilde{U} \setminus \bigcup^{n'}_{t=1} Y_t) \setminus U'| \le C_u k \Vert V' \Vert_*$ such that $U' \in {\cal V}^{s\hat{m}}_{\rm con}$. Moreover, recalling \eqref{rig-eq: Ycrack} as well as $|\partial X_t|_\infty \le \frac{k}{8}$ for $t > n'$, we get $U' \in {\cal V}^{s\hat{m}}_{69k}$ for $m$ sufficiently small. Using \eqref{rig-eq: Ycrack} and the fact that $\tilde{U}$ satisfies \eqref{rig-eq: small bc} we have $\Vert U'\Vert_* \le (1+C_u m)\Vert V \Vert_* + \Vert V' \Vert_* \le  (1 +C_{u}m)\Vert \tilde{U}\Vert_*$. Finally, again using Lemma \ref{rig-lemma: modifica}(ii) we find a set $U \in {\cal V}^{s\hat{m}}_{70k}$ with 
\begin{align} \label{rig-eq: U const}
\Big|\Big(\Omega^{5k} \setminus \big( \bigcup\nolimits^{n'}_{t=1} Y_t \cup \bigcup\nolimits^{n}_{t=n'+1} X_t \cup \bigcup\nolimits_{p \in \hat{J}} Q^\lambda(p) \big) \Big) \setminus U\Big| \le C_u k \Vert U\Vert_*
\end{align}
Arguing similarly as in \eqref{rig-eq: locest3}, \eqref{rig-eq: weaklocA6} we find $\Vert U \Vert_* \le \Vert U' \Vert_* \le (1 +C_{u}m)\Vert \tilde{U}\Vert_*$. This implies \eqref{rig-eq: weaklocA1} since $\tilde{U}$ satisfies \eqref{rig-eq: weaklocA4}. Moreover,  we derive $|(W \setminus U)\cap \Omega^{5k}| \le C_{u}k\Vert U \Vert_*$.

 Define $U_J$ as in the assertion of Lemma \ref{rig-lemma: weaklocA2}. We see that all $\Gamma_t(\tilde{U}_i)=\partial X_t$ with $\Gamma_t(\tilde{U}_i)\cap U^{\circ}_J  \neq \emptyset$ satisfy $\vert \Gamma_t(\tilde{U}_i)\vert_\infty \le \frac{k}{8}$.  In fact, if $\vert \Gamma_t(\tilde{U}_i)\vert_\infty > \frac{k}{8}$,  we would have $X_t \subset \Omega^{5k} \setminus (H^{\frac{k}{8}}(\tilde{U}_i))^\circ$ and thus $X_t \subset \Omega^{5k} \setminus (H^{\frac{k}{8}}(\tilde{U}))^\circ$, where we used  $H^{\frac{k}{8}}(\tilde{U}'') \subset H^{\frac{k}{8}}(\tilde{U}) \subset H^{\frac{k}{8}}(\tilde{U}_4) \subset H^{\frac{k}{8}}(\tilde{U}_i)$ up to a set of negligible measure (see \eqref{rig-eq: new 513**}). (Recall that the set $\tilde{U}''$ given by the modification described below \eqref{rig-eq: isoper} is also denoted by $\tilde{U}$ for convenience.) Therefore, by \eqref{rig-eq: new801} we get $\Gamma_t(\tilde{U}_i) \subset \overline{X_t} \subset \bigcup^{n'}_{j=1}\overline{Z_j}$ and thus $\Gamma_t(\tilde{U}_i) \cap U_J^\circ = \emptyset$ giving a contradiction. Consequently, by \eqref{rig-eq: weaklocA11} 
 \begin{align}\label{rig-eq: new800}
\text{\eqref{rig-eq: weaklocA3} holds for $\bar{y}_i$ for all $\Gamma_t(\tilde{U}_i)$ with $\Gamma_t(\tilde{U}_i)\cap U^{\circ}_J \cap S_i \neq \emptyset$.}
\end{align}
For later we recall that the corresponding components $(\Gamma_{t_s}(G))_s$ with $\Gamma_{t_s}(G)\subset \overline{X_t(\tilde{U}_i)}$ (which satisfy \eqref{rig-eq: weaklocA3}) also satisfy \eqref{rig-eq: newnew} since $G \in {\cal W}^{s\hat{m}}$. Consider $\tilde{Q} := Q_j^{ 3\lambda}(p) \cap  U_J \subset S_i$.  We  observe that $\tilde{Q}$ consists of a  bounded number of squares and that $\tilde{Q} \cap  U_Q$ is contained in a connected component $F$ of $Q^k_i(q) \cap  \hat{W}^\circ$. Indeed, this follows from the fact that due to the construction of $U_Q$, in particular \eqref{rig-eq: new802}, two connected components $F_1 \neq F_2$, $F_t \cap S_i \neq \emptyset$ for $t=1,2$, for which $H(F_t)$ is not completely contained in another component $H(F_{t'})$, fulfill  $\dist(F_1 \cap U_J,F_2\cap U_J) \ge 2\lambda$. This observation also implies that $\tilde{Q}^\circ$ is connected, i.e. each $Q \subset \tilde{Q}$ shares at least one face with the rest of $\tilde{Q}$.  Consequently, Corollary \ref{rig-cor: kornpoin}  together with \eqref{rig-eq: weaklocA3} yield
$$(|E ( R^T_i\bar{y}_i - \id)|(\tilde{Q}))^2 \le C\lambda^2\alpha_{ R_i}(U_Q \cap \tilde{Q}) + Ck C_m \epsilon|\partial \hat{U}_i \cap Q^k_i(q)|^2_{\cal H},$$
where $R_i$ is the value of the constant function $\hat{R}_i|_F$.  Then \eqref{rig-eq: new104} and  \eqref{rig-eq: weaklocA6} imply $|\partial \hat{U}_i \cap Q^k_i(q)|_{\cal H} \le Ck $ which together with \eqref{rig-eq: linearization2} yields \eqref{rig-eq: weaklocA17}. For later we note that Corollary \ref{rig-cor: kornpoin} also yields 
\begin{align}\label{rig-eq: mayblast7}
(|D^j ( \bar{y}_i - R_i\, \id)|(\tilde{Q}))^2 \le C\lambda^2\alpha_{ R_i}(U_Q \cap \tilde{Q}) + Ck C_m \epsilon|\partial \hat{U}_i \cap Q^k_i(q)|^2_{\cal H}.
\end{align}
It remains to show  \eqref{rig-eq: weaklocA2.13}. Consider $\hat{Q} = Q^\lambda(p)$ with $\hat{Q} \cap  U_Q \neq \emptyset$ and show that $|\hat{Q} \cap  U_Q| \ge c m\lambda^2$. First note that $\hat{Q} \cap  U_Q = \hat{Q} \cap \tilde{U}$.  Let $\Gamma=\Gamma(\tilde{U})  = \partial X$ be the boundary component maximizing $\vert X \cap \hat{Q}\vert_\infty$. If $\vert \Gamma \vert_\infty \ge  \frac{k}{8}$ we get a contradiction for  $B$ large enough as then $\hat{Q} \cap U_J = \emptyset$. Assume $\vert X \cap \hat{Q}\vert_\infty \ll \lambda$. Then \eqref{rig-eq: weaklocA6*} and the isoperimetric inequality imply $|\hat{Q} \setminus  U_Q| \le C_{u}\sum_t \vert X_t(\tilde{U}) \cap \hat{Q}\vert^2_\infty \ll C_{u}\lambda\sum_t \vert X_t(\tilde{U}) \cap \hat{Q}\vert_\infty\le C_{u}\lambda^2$ and thus $|\hat{Q} \cap  U_Q| \ge  c m\lambda^2$ for $m$ small enough. Therefore, we may assume that 
\begin{align}\label{rig-eq: infty sandwich}
\tfrac{1}{8}k =\tfrac{1}{8} m^{-1} \lambda \ge \vert \Gamma \vert_\infty \ge \vert X \cap \hat{Q} \vert_\infty \ge \bar{c}\lambda
\end{align}
for $\bar{c}>0$ small enough. It is not hard to see that $| (N_2(\Gamma) \setminus X)  \cap \hat{Q}| \ge CBm  \bar{c}^2\lambda^2$.  Indeed,  an elementary argument yields  $ |N_2(\Gamma)  \cap \hat{Q}| \ge CBm  \bar{c}^2\lambda^2$. Moreover, if we  had  $|\hat{Q} \setminus X| \ll Bm  \bar{c}^2\lambda^2$, we would get $\hat{Q} \subset N_1(\Gamma)$ and thus $\hat{Q} \cap U_Q = \emptyset$ by  the construction of $U_Q$.  We can assume that $ N_2(\Gamma) \cap \partial H^{\frac{k}{8}}(\tilde{U}) = \emptyset$ since otherwise a component larger than $\frac{k}{8}$ intersects $\hat{Q}$ and we derive $\hat{Q} \cap U_J = \emptyset$ as before. By \eqref{rig-eq: small bc2}  this also implies that all components $X_j(\tilde{U})$ with $X_j(\tilde{U}) \cap N_2(\Gamma) \neq \emptyset$ satisfy $\overline{X_j(\tilde{U})} \cap \Gamma = \emptyset$. Thus by the isoperimetric inequality and by \eqref{rig-eq: isoper}  we  get  $|N_2(\Gamma)\cap \hat{Q}\cap \tilde{U}|  \ge | (N_2(\Gamma) \setminus X) \cap \hat{Q}| - C(B\lambda m)^2$. This implies
\begin{align*}
|\hat{Q} \cap  U_Q| & = |\hat{Q} \cap \tilde{U}| \ge |\hat{Q} \cap \tilde{U} \cap N_2(\Gamma)| \ge | (N_2(\Gamma) \setminus X)  \cap \hat{Q}| - C(B\lambda m)^2 \\
& \ge -CB^2\lambda^2 m^2 +  C\bar{c}^2Bm \lambda^2 \ge cm \lambda^2 
\end{align*}
for $m$ sufficiently small.  \eop

\begin{rem}\label{rig-rem: crac bd}
{\normalfont
(i) For later we observe that there is  a set $U^H \in {\cal V}^{\lambda}_{35k}$ with 
\begin{align}\label{rig-eq: abs crac}
(i) \   \Vert U^H\Vert_* \le  (1+C_{u}m)\Vert W \Vert_* + C\epsilon^{-1}(\gamma + \delta_4), \  \  (ii) \   \Vert U^H\Vert_{\cal H} \le  C_{u}\Vert U^H\Vert_*
\end{align}
which coincides with the set $U_J$ considered in the previous lemma up to a set of negligible measure.  In fact, we apply Lemma \ref{rig-lemma: modifica}(i) on the rectangles $(Z_t)^{n'}_{t=1}$ considered in \eqref{rig-eq: new802} and find pairwise disjoint  $(Z'_t)^{n'}_{t=1}$ with  $\bigcup^{n'}_{j=1} \overline{Z_j} = \bigcup^{n'}_{j=1} \overline{Z'_j}$. We define 
$$U^H := \Omega^{5k} \setminus \Big( \bigcup\nolimits^{n'}_{j=1} Z'_j \cup  \bigcup\nolimits_{p \in \hat{J}} Q^\lambda(p)  \Big) ,$$
where $\hat{J}$ as in \eqref{rig-eq: U const}. By  Lemma \ref{rig-lemma: modifica}(i) we get (i) and  $U^H \in {\cal V}^{\lambda}_{35k}$ since $|\pi_i \partial Z_t| \le 2\cdot 34k + C_u m k \le 70k$ for $i=1,2$. Moreover, \eqref{rig-eq: abs crac}(ii) is a consequence of Lemma \ref{rig-lemma: infty} and the fact that $(Z_t)_t$, $(Q^\lambda(p))_{p \in J}$ are rectangles.

Clearly $U_J \subset U^H$. Moreover, we see that $|U^H \setminus U_J|>0$ can only happen if there is a square $Q^\lambda(p) \subset U^H$ and components $(X_{t}(\tilde{U}))_t$ of $\tilde{U}$ such that $Q^\lambda(p) \subset \bigcup_t \overline{X_{t}(\tilde{U})}$. Since we can suppose $|\partial X_{t}(\tilde{U})|_\infty \le \frac{k}{8}$ (otherwise the components are contained in some rectangle $Z_t$), this yields a contradiction to \eqref{rig-eq: small bc2}. 

(ii) For $i=1,\ldots,4$ we have 
$$|\partial \hat{U}_i \cap U^\circ_J|_{\cal H} \le C_{u} (\Vert W \Vert_* + C\epsilon^{-1}(\gamma + \delta_4)).$$
In fact, recalling \eqref{rig-eq: new800} we get that all $\Gamma_t(\hat{U}_i)$ with $\Gamma_t(\hat{U}_i) \cap U_J^\circ \neq \emptyset$ fulfill \eqref{rig-eq: weaklocA3}  and \eqref{rig-eq: newnew}. Thus, we obtain $| \Theta_t(\hat{U}_i)|_{\cal H} \le C_u| \Theta_t(\hat{U}_i)|_*$ and the claim follows from \eqref{rig-eq: weaklocA4} replacing $\hat{U}_1$ by $\hat{U}_i$.

}

\end{rem}

We are now in a position to prove the main result of this section. Recall the definition $\lambda = sdm^{-1} = km$ and \eqref{rig-eq: vartheta def}.

\begin{lemma}\label{rig-lemma: weaklocA}
Let $ k > s,  \epsilon   > 0$ such that $l:=\frac{k}{s} = dm^{-2}$ for $m^{-1},d \in \N$ with $m^{-1}, d \gg 1$.  Then for a fixed constant $C>0$ we have the following: \\
For all $W \in  {\cal V}^{sm}_{(s,3k)}$ with $W \subset \Omega^{3k}$ and for all $y \in H^1(W)$ with $\Vert \nabla y \Vert_\infty\le C$, $\gamma$ as defined in \eqref{rig-eq: weakloc gamma} and 
\begin{align}\label{rig-eq: weaklocA5}
\delta_4 := \sum\nolimits^4_{i=1}\Vert \nabla y  - \hat{R}_i\Vert^4_{L^4(W)},  \ \ \  \delta_2 := \sum\nolimits^4_{i=1}\Vert \nabla y  - \hat{R}_i\Vert^2_{L^2(W)}
\end{align}
for mappings $\hat{R}_i:  W^\circ \to SO(2)$, $i=1,\ldots,4$, which are constant on the connected components of $Q^{k}_i(p) \cap  W^\circ$, $p \in I^{k}_i(\Omega^{3k})$, we obtain:\\
We find sets $V \in {\cal V}^{s\hat{m}^2}_{71k}$, $U_J \in  {\cal V}^\lambda$ with $V \subset U_J$ and  $V \subset \Omega^{6k}$, $|V \setminus W| = 0$, $|(W \setminus V) \cap \Omega^{6k}| \le C_{u}k\Vert V \Vert_*$ such that 
\begin{align}\label{rig-eq: weaklocA1.2}
\begin{split}
\Vert V \Vert_* \le (1 + C_{u}m)\Vert W \Vert_* + C\epsilon^{-1}(\gamma+ \delta_4)
\end{split}
\end{align}
as well as  mappings $\bar{R}_j: U_J\to SO(2)$ and $\bar{c}_j: U_J\to  \R^2$,  which are constant on $Q^{\lambda}_j(p)$, $p \in I^{\lambda}_j(\Omega^{3k})$, such that
\begin{align}\label{rig-eq: weaklocA2}
(i) & \ \ \Vert y  - (\bar{R}_j \, \cdot + \bar{c}_j)\Vert^2_{L^2(V)} \le C C_m^2\lambda^2 \min\nolimits_{p=2,4}(1+ \vartheta_p)(\gamma + \delta_p + \epsilon \Vert W \Vert_*),\notag \\
(ii) & \ \ \Vert \nabla y - \bar{R}_j\Vert^p_{L^p(V)} \le  C C^2_m\big(\delta_p + \vartheta_p(\gamma+ \delta_4 +\epsilon \Vert W \Vert_*)\big) , \ p=2,4, \\
(iii) & \ \  \Vert \bar{R}_{j_1} - \bar{R}_{j_2}\Vert^p_{L^p(U_J)} \le  C C^2_m\big(\delta_p + { \vartheta_p}(\gamma+ \delta_4 +\epsilon \Vert W \Vert_*)\big) , \ p=2,4, \notag \\
(iv) & \ \ \Vert (\bar{R}_{j_1} \, \cdot + \bar{c}_{j_1})  - (\bar{R}_{j_2} \, \cdot + \bar{c}_{j_2})\Vert^2_{L^2(U_J )}\le CC_m^2 \lambda^2  \min_{p=2,4}(1+ \vartheta_p) (\gamma + \delta_p + \epsilon \Vert W \Vert_*)\notag 
\end{align}
for $j_1,j_2 = 1,\ldots,4$, $j=1,\ldots,4$, where $\vartheta_4 = \vartheta$ and $\vartheta_2 = 1$. Moreover, we have
\begin{align}\label{rig-eq: weaklocA9}
\lambda^{-2}\Vert (\bar{R}_{j_1} \, \cdot + \bar{c}_{j_1})  - (\bar{R}_{j_2} \, \cdot + \bar{c}_{j_2})\Vert^2_{L^\infty(U_J )}+ \Vert  \bar{R}_{j_1} - \bar{R}_{j_2}\Vert^4_{L^\infty(U_J)}\le C \bar{\vartheta}
\end{align}
for $\bar{\vartheta} = \min \lbrace \vartheta(1+\vartheta), C_m^3\rbrace$ and under the additional assumption that $\Delta y = 0$ in  $W^\circ$ we obtain
\begin{align}\label{rig-eq: weaklocA9.2}
\lambda^{-2}\Vert y  - (\bar{R}_j \, \cdot + \bar{c}_j)\Vert^2_{L^\infty(V)} \le C   \vartheta (1+\vartheta).
\end{align}

\end{lemma}

\Proof Apply Lemma \ref{rig-lemma: weaklocA2} to obtain $U \in {\cal V}^{s\hat{m}}_{70k}$, $U_Q\in {\cal V}^{s\hat{m}}$ with $| U_Q \setminus W|=0$, $U_J$ and extensions $\bar{y}_i: S_i \cap U_J \to \R^2$ such that \eqref{rig-eq: weaklocA1}, \eqref{rig-eq: weaklocA2.13} and \eqref{rig-eq: weaklocA17} hold. Consider $Q = Q^{\lambda}_j(p)$, $p \in I^{\lambda}_j(\Omega^{3k})$, $j=1,\ldots,4$, with $Q \cap U_J \neq \emptyset$. Moreover, let $\tilde{Q} = Q^{{3}\lambda}_j(p) \cap U_J$.  As $6\lambda < \frac{k}{4}$ by $m \ll 1$, we find some $Q^{k}_i(q)$  for some $i=1,\ldots,4$ with $\tilde{Q} \subset Q^{\frac{5}{8}k}_i(q) \subset S_i$ and therefore we can apply \eqref{rig-eq: weaklocA17}. Recall that $\hat{R}: =  \hat{R}_i|_{W^\circ \cap \tilde{Q}}$ is constant due to the construction in Lemma \ref{rig-lemma: weaklocA2} (see below \eqref{rig-eq: new800}). By Theorem \ref{rig-theo: korn} we find $A \in \R^{2 \times 2}_{\rm skew}$ and $c \in \R^2$ such that 
\begin{align}\label{rig-eq: CG}
\begin{split}
\Vert \bar{y}_i - \hat{R}\,(\Id + A)\,\cdot - \hat{R}\,c\Vert^2_{L^2( \tilde{Q})} & =  \Vert \hat{R}^T \bar{y}_i - \cdot - (A\, \cdot + c)\Vert^2_{L^2(\tilde{Q})} \\ & \le C(|E(\hat{R}^T\bar{y}_i - \id)|(\tilde{Q}))^2 \le Ck^2 G,
\end{split}
\end{align}
where
$$G:= C_m\min\Big\{ \epsilon k, \gamma(W \cap Q_i^{2k}(q))  + \delta_4(W \cap Q_i^{2k}(q)) + \epsilon|\partial W \cap Q_i^{2k}(q)|_{\cal H}\Big\}.$$
The constant $C$ is independent of $\tilde{Q}$ as there are (up to rescaling) only a finite number of different shapes of $\tilde{Q}$. (Also recall that  each $Q \subset \tilde{Q}$ shares at least one face with the rest of $\tilde{Q}$.)

In the proof of Lemma \ref{rig-lemma: weaklocA2} we have seen that  all $\Gamma_t = \Gamma_t(\tilde{U}_i)$ with $\tilde{Q} \cap \Gamma_t \neq \emptyset$ satisfy \eqref{rig-eq: weaklocA3}  for $\bar{y}_i$ and $|\Gamma_t|_\infty \le \frac{k}{8}$ as well as $ N^{\tau_l}(\partial R_t) \subset Q^k_i(q)$  (cf. \eqref{rig-eq: new800}). Thus, by Lemma \ref{rig-lemma: A neigh}  for $V = Q_i^k(q)$ we get
\begin{align}\label{rig-eq: A2--3***}
\begin{split}
\Vert \nabla \bar{y}_i- \hat{R}  \Vert^p_{L^p(\tilde{Q})} & \le  C\Vert \nabla y - \hat{R} \Vert^p_{L^p(\tilde{Q} \cap   \hat{W})} + C\sum\nolimits_{\Gamma_t \in {\cal F}(Q^k_i(q))} |X_t|_\infty^2 |A_t|^p \\
&  \le CC_m\delta_p(Q^k_i(q)\cap  \hat{W}) + CC_m (\epsilon s^{-1})^{\frac{p}{2} -1}  \epsilon |\partial \hat{U}_i \cap  Q^k_i(q)|_{\cal H}
\end{split}
\end{align}
for $p=2,4$,  where $\hat{W},\hat{U}_i$ as defined in the previous proof and $X_t, A_t$ as in \eqref{rig-eq: construction}. Recall that the factor $s^{-1}$ appearing in the estimate is related to the fact that the least length of boundary components of $\hat{U}_i$ is $s$.  Thus, recalling that $\hat{U}_i$ fulfills \eqref{rig-eq: new104} we obtain by the definition of $G$ 
\begin{align}\label{rig-eq: A2--3}
\Vert \nabla \bar{y}_i- \hat{R}  \Vert^p_{L^p(\tilde{Q})} \le CC_m\delta_p(Q^k_i(q)\cap  \hat{W}) + C (\epsilon s^{-1})^{\frac{p}{2} -1}  G  =: H_p.
\end{align}
We repeat the estimate \eqref{rig-eq: CG} with  the Poincar\'e inequality in SBV (see \cite[Remark 3.50]{Ambrosio-Fusco-Pallara:2000}) instead of Theorem \ref{rig-theo: korn} and obtain by \eqref{rig-eq: mayblast7} and H\"older's inequality
\begin{align*}
\Vert \bar{y}_i - \hat{R}\,\,\cdot - \tilde{c}\Vert^2_{L^2( \tilde{Q})} &\le C\Vert \nabla \bar{y}_i - \hat{R}\Vert^2_{L^1(\tilde{Q})} +  C(|D^j(\bar{y}_i - \hat{R}\, \id)|(\tilde{Q}))^2 \\
& \le C\lambda^{4(1-\frac{1}{p})}H_p^{\frac{2}{p}} + Ck^2 G,
\end{align*}
for $\tilde{c} \in \R^2$ for $p=2,4$. This together with \eqref{rig-eq: CG} and an argumentation similar to \eqref{rig-eq: A difference} (see also (2.11) in \cite{Friedrich:15-1}, where such an estimate is derived in the geometrically linear setting) yields $\lambda^4|A|^2 \le C\lambda^{4-4/p}H_p^{2/p} + Ck^2G $ and therefore by \eqref{rig-eq: A2--3}
\begin{align}\label{rig-eq: A2--3*****}
\begin{split}
\lambda^2 |A|^2 &\le CH_2 + Cm^{-2}G \le CC_m\delta_2(Q^k_i(q)\cap  \hat{W}) + Cm^{-2}G =: \hat{H}_2, \\
\lambda^2 |A|^4 &\le CH_4 + C\lambda^{-2} m^{-4}G^2 \le CH_4 +C\lambda^{-1} m^{-5}C_m \epsilon G \\
& \le C C_m\delta_4(Q^k_i(q)\cap  \hat{W}) +C\vartheta G =: \hat{H}_4.
\end{split}
\end{align}
Observe that $\hat{H}_4 \le C(1+ \vartheta) G$. By \eqref{rig-eq: linearization} there is a rotation $\bar{R} \in SO(2)$ such that
\begin{align}\label{rig-eq: R,A diff}
\begin{split}
|\bar{R} - \hat{R}(\Id + A)|^2 & = \dist^2(\hat{R}(\Id + A),SO(2)) \\
& \le  0 + C|\hat{R}(\Id + A) - \hat{R}|^4   = C|A|^4  \le  C\lambda^{-2} \hat{H}_4,
\end{split}
\end{align}
as $\bar{e}_{\hat{R}} (\hat{R}(\Id + A)) = 0$. Likewise,  as $|A| \le C$ by  $\Vert \nabla y\Vert_\infty \le C$ we get $|\bar{R} - \hat{R}(\Id + A)|^2 \le  C|A|^2  \le  C\lambda^{-2}\hat{H}_2$. Consequently, the Poincar\'e inequality, \eqref{rig-eq: CG}  and \eqref{rig-eq: A2--3*****}  yield
\begin{align}\label{rig-eq: weaklocA14}
\Vert \bar{y}_i - (\bar{R}\,\cdot +\bar{c})\Vert^2_{L^2(\tilde{Q})} \le Ck^2 G + C\lambda^4 |A|^4\le Ck^2 G +  Ck^2 \min\nolimits_{p=2,4} \hat{H}_p 
\end{align}
for some possibly different $\bar{c} \in \R^2$. Moreover, we get
\begin{align}\label{rig-eq: R diff2}
\begin{split}
\lambda^2|\hat{R}- \bar{R}|^4 & \le C\lambda^2|\bar{R} - \hat{R}(\Id + A)|^4 + C\lambda^2|A|^4  \\ &\le C\lambda^2|\bar{R} - \hat{R}(\Id + A)|^2 + C\lambda^2|A|^4 \le   C  \hat{H}_4.
\end{split}
\end{align}
and likewise 
\begin{align}\label{rig-eq: weaklocA15}
\lambda^2|\hat{R}- \bar{R}|^2 \le     C\hat{H}_2.
\end{align}
For fixed $j=1,\ldots,4$ we proceed in this way on each $Q_t= Q^{\lambda}_j(p)$, $p \in I^{\lambda}_j(\Omega^{3k})$, with $Q_t \cap U_J \neq \emptyset$ and for the corresponding $\tilde{Q}_t = Q^{ 3\lambda}_j(p) \cap U_J$ we obtain constants $\hat{R}_t, \bar{R}_t \in SO(2)$ and $\bar{c}_t \in  \R^2$  as given in \eqref{rig-eq: weaklocA14}-\eqref{rig-eq: weaklocA15}. Consequently, we find mappings $\bar{R}_j:  U_J \to SO(2)$ and $\bar{c}_j: U_J \to  \R^2$ being constant on each $Q_t$, where on each $Q_t \subset \tilde{Q}_t$ we choose  $\bar{R}_j = \bar{R}_t$ and $\bar{c}_j = \bar{c}_t$.  By \eqref{rig-eq: weaklocA14} and the observation that every $Q^{2k}_i(q)$ is intersected only by $\sim m^{-2}$ squares $\tilde{Q}_t$ we obtain 
\begin{align}\label{rig-eq: weaklocA21}
\begin{split}
\Vert y  - (\bar{R}_j \, \cdot + \bar{c}_j)\Vert^2_{L^2(U)} & \le Ck^2  \min_{p=2,4}(1+\vartheta_p) m^{-2} C_m m^{-2} ( \gamma + \delta_{p} + \epsilon\Vert  W\Vert_*)  \\
& \le C \lambda^2 \min\nolimits_{p=2,4}(1+\vartheta_p) m^{2} C^2_m (\gamma + \delta_p + \epsilon \Vert W \Vert_*)
\end{split}
\end{align}
where $\vartheta_2 = 1$ and $\vartheta_4 = \vartheta$. Here we used that $\delta_4 \le C\delta_2$. Likewise, applying \eqref{rig-eq: weaklocA5}, \eqref{rig-eq: A2--3*****}, \eqref{rig-eq: R diff2}, \eqref{rig-eq: weaklocA15} as well as the triangle inequality we  get 
\begin{align}\label{rig-eq: weaklocA22}
\begin{split}
\Vert \nabla y - \bar{R}_j\Vert^p_{L^p(U)}  &\le   Cm^{-2}  C_m \big( \delta_p +  m^{-2}\vartheta_p (\gamma+ \delta_4 + \epsilon \Vert W \Vert_*)\big) \\ &\le   Cm C^2_m\big( \delta_p +  \vartheta_p (\gamma+ \delta_4 + \epsilon \Vert W \Vert_*)\big)
\end{split}
\end{align}
for $p=2,4$. We now consider  $Q_1 := Q^{\lambda}_{j_1}(p_1)$, $Q_2 :=Q^{\lambda}_{j_2}(p_2)$ with $Q_1 \cap  Q_2 \neq \emptyset$ and $Q_1, Q_2 \cap U_J \neq \emptyset$. Moreover, let $\tilde{Q}_i = Q^{3\lambda}_{j_i}(p_i) \cap U_J$ be the corresponding enlarged sets. It is not hard to see that there is some $Q^{\lambda}(p)$, $p \in  J(U_Q)$, with $Q^{\lambda}(p) \subset \tilde{Q}_1,\tilde{Q}_2$ and therefore by the definition of $U_J$, in particular \eqref{rig-eq: weaklocA2.13}, we derive $|\tilde{Q}_1 \cap \tilde{Q}_2 \cap  U_Q| \ge  cm\lambda^2$. Let $\bar{R}_{j_i} \in SO(2)$, $\bar{c}_{j_i} \in \R^2$, $i=1,2$, be the constants constructed above. We compute
\begin{align}\label{rig-eq: weaklocA22.2.2}
\begin{split}
\lambda^2\Vert \bar{R}_{j_1} - \bar{R}_{j_2}\Vert^p_{L^\infty(Q_1 \cap Q_2)} &\le Cm^{-1}\Vert \bar{R}_{j_1} - \bar{R}_{j_2}\Vert^p_{L^p(\tilde{Q}_1 \cap \tilde{Q}_2 \cap  U_Q)} \\
& \le Cm^{-1} \sum\nolimits^4_{j=1}\Vert \nabla y - \bar{R}_j\Vert^p_{L^p(\tilde{Q}_1 \cap \tilde{Q}_2 \cap  U_Q)}
\end{split}
\end{align}
and summing over all squares we get by \eqref{rig-eq: weaklocA22} 
\begin{align}\label{rig-eq: weaklocA22.2}
 \Vert \bar{R}_{j_1} - \bar{R}_{j_2}\Vert^p_{L^p(U_J)} \le  C C^2_m\big( \delta_p +  \vartheta_p (\gamma+ \delta_4 + \epsilon \Vert W \Vert_*)\big)
\end{align}
for  $1 \le j_1, j_2 \le 4$ and $p=2,4$. Here we used that each $Q^{ 3\lambda}_j(p) \cap U_J$ only appears in a finite number of addends. Note that $\frac{|\pi_1(Q_1 \cap Q_2)| + |\pi_2(Q_1 \cap Q_2)| }{ \max_{i=1,2} |\pi_i( \tilde{Q}_1 \cap \tilde{Q}_2 \cap  U_Q )|} \le Cm^{-1/2}$ and $\frac{|Q_1 \cap Q_2|  }{ | \tilde{Q}_1 \cap \tilde{Q}_2 \cap U_Q |} \le Cm^{-1}$.  Consequently, arguing similarly as in \eqref{rig-eq: A,c difference2} we find 
 \begin{align}\label{rig-eq: weaklocA22.3.2}
 \begin{split}
\lambda^2\Vert (\bar{R}_{j_1} \, & \cdot + \bar{c}_{j_1})  - (\bar{R}_{j_2} \, \cdot + \bar{c}_{j_2})\Vert^2_{L^\infty(Q_1 \cap Q_2)} \\& \le C(m^{-\frac{1}{2}})^2 m^{-1} \Vert (\bar{R}_{j_1} \, \cdot + \bar{c}_{j_1})  - (\bar{R}_{j_2} \, \cdot + \bar{c}_{j_2})\Vert^2_{L^2(\tilde{Q}_1 \cap \tilde{Q}_2 \cap  U_Q)}.
 \end{split}
\end{align}
Replacing \eqref{rig-eq: weaklocA22} by \eqref{rig-eq: weaklocA21} in the above argument we then get
 \begin{align}\label{rig-eq: weaklocA22.3}
 \begin{split}
\hspace{-0.05cm}\Vert (\bar{R}_{j_1} \, \cdot + \bar{c}_{j_1})  - (\bar{R}_{j_2} \, \cdot + \bar{c}_{j_2})\Vert^2_{L^2(U_J)} &\le Cm^{-2}\sum\nolimits^4_{j=1} \Vert y  - (\bar{R}_j \, \cdot + \bar{c}_j)\Vert^2_{L^2( U_Q)} \\
& \le C C_m^2 \lambda^2    \min\nolimits_{p=2,4}(1+\vartheta_p) (\gamma + \delta_{p} + \epsilon \Vert W \Vert_*).
\end{split}
\end{align}
Similarly as in the proof of Lemma \ref{rig-lemma: weaklocR} (see the construction in \eqref{rig-eq: new111}) we can define $V \in {\cal V}^{s\hat{m}^2}_{71k}$   with $|V\setminus U|=0$, $ V^\circ \subset \lbrace x \in U  \cap \Omega^{6k}: \dist_\infty(x, \partial U)  \ge 2s \hat{m}m \rbrace$, $\Vert V\Vert_* \le (1 + C_{u}m) \Vert U \Vert_*$  and  $|( W \setminus V) \cap \Omega^{6k}| \le C_{u}{ k}\Vert  V \Vert_*$. By \eqref{rig-eq: weaklocA1} this implies \eqref{rig-eq: weaklocA1.2}.   We note that in this case for components $\Gamma_j = \partial X_j$ with $X_j \subset U_J$ it suffices to consider a corresponding rectangle $M(\Gamma_j)$ with $M(\Gamma_j) \subset U_J$.   For later we observe that this construction yields
\begin{align}\label{rig-eq: mayblast5}
 V \subset U_J, \ \ \ \ \ \ \ \big| \big(\Omega^{6k} \setminus \bigcup\nolimits M(\Gamma_j)\big) \triangle V\big| = 0.
\end{align}
We now see that \eqref{rig-eq: weaklocA2} follows directly from \eqref{rig-eq: weaklocA21}-\eqref{rig-eq: weaklocA22.3}.
 
It remains to show \eqref{rig-eq: weaklocA9} and \eqref{rig-eq: weaklocA9.2}.  By  \eqref{rig-eq: A2--3*****}, \eqref{rig-eq: R diff2} and \eqref{rig-eq: weaklocA22.2.2} we find $\Vert \bar{R}_{j_1} - \bar{R}_{j_2}\Vert^4_{L^\infty(Q_1 \cap Q_2)} \le   C\lambda^{-2}(1+\vartheta)G + C\lambda^{-2}m^{-1} G$ for sets $Q_1,Q_2 \subset U_J$ as considered above.  Recalling the definition of $G$ we then get 
$$\Vert \bar{R}_{j_1} - \bar{R}_{j_2}\Vert^4_{L^\infty(Q_1 \cap Q_2)}\le C  (1+\vartheta)\lambda^{-2} m^{-1} C_m    \epsilon k \le Cs^{-1}  (1+\vartheta) C^{2}_m \epsilon \le  C  (1+\vartheta)\vartheta $$
Likewise,   we derive $\lambda^{-2}\Vert (\bar{R}_{j_1} \, \cdot + \bar{c}_{j_1})  - (\bar{R}_{j_2} \, \cdot + \bar{c}_{j_2})\Vert^2_{L^\infty(Q_1 \cap Q_2 )} \le C(1+\vartheta)\vartheta$ recalling the definition of $G$ and taking  \eqref{rig-eq: weaklocA22.3.2}, \eqref{rig-eq: weaklocA14}  (for $p=4$)  and the triangle inequality into account. Similarly, by \eqref{rig-eq: weaklocA14} for $p=2$ and the observation that $\delta_2(Q^k_i
(q) \cap \hat{W}) \le Ck^2$ as $\Vert \nabla y \Vert_\infty \le C$ we find  using $\epsilon \le k$
\begin{align*}
\lambda^{-2}\Vert (\bar{R}_{j_1} \, \cdot + \bar{c}_{j_1})  - (\bar{R}_{j_2} &\, \cdot + \bar{c}_{j_2})\Vert^2_{L^\infty(Q_1 \cap Q_2 )} \le C\lambda^{-4}m^{-2}k^2(G + \hat{H}_2)\\
& \le C\lambda^{-2}m^{-4} (m^{-2}G + C_m k^2) \le C\lambda^{-2}C_m^2k^2 \le C C_m^3.
\end{align*}
This finishes the proof of \eqref{rig-eq: weaklocA9}.

Finally, to see \eqref{rig-eq: weaklocA9.2}, we repeat the argument in \eqref{rig-eq: uniform estimate}: Let $x \in Q \cap V \subset \tilde{Q}$ for  $Q = Q^{\lambda}_j(p)$ , $\tilde{Q} = Q^{ 3\lambda}_j(p) \cap U_J$ as considered above and let $\bar{R}\,\cdot + \bar{c}$ be the corresponding rigid motion as given in \eqref{rig-eq: weaklocA14}.  Since $y$ is assumed to be  harmonic in $U^\circ$ the mean value property of harmonic function for $r \le s \hat{m}m$ and Jensen's inequality yield
\begin{align*}
| y(x) - (\bar{R}\,x + \bar{c}) |^2 & \le \Big|\frac{1}{|B_r(x)|} \int_{B_r(x)} (y(t) - (\bar{R}\,t + \bar{c})) \, dt \Big|^2 \\&\le C|B_r(x)|^{-1}  (1+\vartheta) k^2 G 
 \le C  (1+\vartheta) m^{-2}  \hat{m}^{-2} s^{-2} k^2 G  \\ &\le  C(1+\vartheta)C_m m^{ -4}  \hat{m}^{-2} l \epsilon s^{-1} \lambda^2 \le C  (1+\vartheta)\vartheta\lambda^2.
\end{align*}
 Here we used \eqref{rig-eq: weaklocA14} and the fact that $B_r(x) \subset  U^\circ \cap \tilde{Q}$ for all $x \in Q \cap V$.  \eop

\subsection{Local rigidity for an extended function}\label{rig-sec: subsub,  ext}

We now state a version of Lemma \ref{rig-lemma: weaklocA} for an extension of the function $y$.

\begin{corollary}\label{rig-cor: weakA}
Let be given the assumptions of Lemma \ref{rig-lemma: weaklocA2}, Lemma \ref{rig-lemma: weaklocA}  and let $U \in  {\cal V}^{s\hat{m}}_{70k}$, $U^H \in {\cal V}^{\lambda}_{35k}$ be the sets provided by Lemma \ref{rig-lemma: weaklocA2}, Remark \ref{rig-rem: crac bd}, respectively.  Moreover, assume that $\vartheta \le 1$. Then the estimates \eqref{rig-eq: weaklocA2}(iii),(iv) hold on $U^H$ for functions $\bar{R}_j$, $\bar{c}_j$, $j=1,\ldots,4$.   Moreover, we find an extension  $\hat{y} \in SBV^2(U^H,\R^2)$ with $\hat{y} = y$ on $U$ and $\nabla \hat{y} \in SO(2)$ on $U^H \setminus W$ a.e. such that for every $Q = Q^{\lambda}_j(p)$, $p \in I^{\lambda}_j(\Omega^{3k})$, with  $Q \cap U^H \neq \emptyset$ we have
\begin{align}\label{rig-eq: corweak}
\begin{split}
(i)& \ \ \Vert \nabla \hat{y} - \bar{R}_j\Vert^p_{L^p(Q)} \le CC_m^2  \, (\bar{G}(N) + \delta_p(N)), \ p=2,4\\
(ii)& \ \ \Vert \hat{y}  - (\bar{R}_j \, \cdot + \bar{c}_j)\Vert^2_{L^2( Q)} \le C\lambda^2C_m^2  \, \min \lbrace \epsilon k, \bar{G}(N) \rbrace,\\
(iii)& \ \ \Vert \hat{y}  - (\bar{R}_j \, \cdot + \bar{c}_j)\Vert^2_{L^1(\partial Q)} \le C\lambda^2C_m^2  \, \min \lbrace \epsilon k, \bar{G}(N)\rbrace,
\end{split}
\end{align}
where $N=N(Q) = \lbrace x \in W: \dist(x, Q) \le Ck\rbrace$ and  for shorthand $\bar{G}(N) = \gamma(N) + \delta_4(N) + \epsilon {\cal H}^1(N \cap \partial  W)$.  Furthermore, we have 
\begin{align}\label{rig-eq: corweak**}
{\cal H}^1(J_{\hat{y}}) \le C_{u} (\Vert W \Vert_* + C\epsilon^{-1}(\gamma + \delta_4)).
\end{align}
\end{corollary}

\Proof Recall the definition of $U$ in  \eqref{rig-eq: U const} and that $U_J$ and $U^H$ coincide up to a set of measure zero by Remark \ref{rig-rem: crac bd}. In Lemma \ref{rig-lemma: weaklocA2} we have defined sets $(\tilde{U}_j)^4_{j=1}$, $\tilde{U}^*_4 \subset \ldots \subset \tilde{U}^*_1$ (see \eqref{rig-eq: new 513**})  and corresponding extensions $\bar{y}_i|_{U_J \cap S_i}$. Moreover,  in \eqref{rig-eq: new800} have seen that all $\Gamma_t(\tilde{U}_i)$ with $\Gamma_t(\tilde{U}_i) \cap  U_J^\circ \cap S_i \neq \emptyset$ satisfy \eqref{rig-eq: weaklocA3}  for $\bar{y}_i$ and $|\Gamma_t(\tilde{U}_i)|_\infty \le \frac{k}{8}$. By Lemma \ref{rig-lemma: weaklocA} we get that  \eqref{rig-eq: weaklocA2}(iii),(iv) hold.

The goal is to provide one single extension $\hat{y}: U^H \to \R^2$ and to confirm \eqref{rig-eq: corweak}. Define
$$\hat{S}_i := \bigcup\nolimits_{p \in I^k_i(\Omega^{3k})}  \overline{Q_i^{\frac{9}{16}k}(p)}\subset S_i$$ 
and let $D_i =  (\tilde{U}_i \cap U^\circ_J) \cup \bigcup_{\Gamma_t(\tilde{U}_i) \subset \hat{S}_i} X_t(\tilde{U}_i)$, where $X_t(\tilde{U}_i)$ is the component corresponding to $\Gamma_t(\tilde{U}_i)$.  We now show that $U^\circ_J \subset \bigcup\nolimits_{i=1}^4 D_i.$
To see this, it suffices to prove
\begin{align}\label{rig-eq: claimSi}
S_i  \cap U^\circ_J \subset \bigcup\nolimits^4_{n=1} D_n, \ \ i=1,\ldots,4.
\end{align}
Fix $i$ and assume that \eqref{rig-eq: claimSi} has already be established for $j > i$. As $S_i \cap U^\circ_J \subset  \Omega^{5k} \subset  H(\tilde{U}_i) =  \tilde{U}_i \cup \bigcup_{\Gamma_t(\tilde{U}_i)} X_t(\tilde{U}_i)$ by the definition of $U_J$, we find $(S_i \cap U^\circ_J) \setminus D_i \subset (S_i \cap U^\circ_J)  \cap \bigcup_{\Gamma_t(\tilde{U}_i) \not\subset \hat{S}_i} X_t(\tilde{U}_i)$.  To see \eqref{rig-eq: claimSi} for $i$, it now suffices to show that each $\Gamma_t(\tilde{U}_i)$ with $\Gamma_t(\tilde{U}_i) \cap  U_J^\circ \cap S_i \neq \emptyset$ satisfies $U_J^\circ \cap X_t(\tilde{U}_i) \subset \bigcup_{n=1}^4 D_n$. Since $|\Gamma_t(\tilde{U}_i)|_\infty \le \frac{k}{8}$ for all such components, we derive  $ X_t(\tilde{U}_i) \subset \hat{S}_j$ for some $j=1,\ldots,4$.  If $j < i$, by the construction of the sets  $\tilde{U}^*_1 \supset \ldots \supset \tilde{U}^*_4$ we find $(X_{t_s}(\tilde{U}_j))_s$ such that 
$$X_t(\tilde{U}_i) = (\tilde{U}_j \cap X_t(\tilde{U}_i)) \cup \bigcup\nolimits_s X_{t_s}(\tilde{U}_j).$$
As $X_{t_s}(\tilde{U}_j) \subset \hat{S}_j$, this implies  $X_t(\tilde{U}_i) \cap U^\circ_J \subset D_j$. The case $j=i$ is clear. If $j>i$, we obtain $X_t(\tilde{U}_i) \cap U^\circ_J \subset S_j \cap U^\circ_J \subset \bigcup\nolimits^4_{n=1} D_n$ by \eqref{rig-eq: claimSi}. This yields the claim. 

Set $\bar{y} = \bar{y}_4 $ on $D_4 \cap U_J$,  $\bar{y} = \bar{y}_j $ on $(D_{j} \setminus D_{j+1}) \cap U_J$ for $j=3,2,1$.  It is not hard to see that $\bar{y}$ is defined on $U^H$  (as $|U^H \setminus U_J^\circ| =  0$) and $\bar{y} = y$ on $U$. Moreover, by construction  there is a set of components $(X_t)_t$ consisting of components of $(\hat{U}_i)_i$ such that 
$$J_{\bar{y}} \subset \bigcup\nolimits_t \partial X_t \subset \bigcup\nolimits^4_{i=1} \bigcup\nolimits_t \Gamma_t( \hat{U}_i).$$
By  \eqref{rig-eq: construction} we have $\bar{y}(x) =  \bar{y}_{i_t}(x) = R_{ t}\,(\Id + A_t)\, x + R_{ t}\,c_t$ for $x \in X_t$, where $R_t \in SO(2)$, $A_t \in \R^{2 \times 2}_{\rm skew}$, $c_t \in \R^2$ and $1 \le i_t \le 4$ appropriately. Note that due the the definition of the extensions in \eqref{rig-eq: construction} the components $X_t$ are associated to the sets $(\hat{U}_i)_i$, not to $(\tilde{U}_i)_i$. By Remark \ref{rig-rem: crac bd}(ii) this yields \eqref{rig-eq: corweak**} for $\bar{y}$.

Consider $Q = Q^{\lambda}_j(p)$ with  $Q \cap U_J \neq \emptyset$. Let $\tilde{Q} = Q^{{3}\lambda}_j(p) \cap U_J$  and observe $|\tilde{Q} \cap U_J| \sim \lambda^2$.  Let $ {\cal I} \subset \lbrace 1, \ldots, 4\rbrace$ such that for each $\iota \in {\cal I}$ we can select some $Q^k_\iota(q_\iota)$ such that $\tilde{Q} \subset Q^{\frac{5}{8}k}_\iota(q_\iota)$. Note that $\# {\cal I}> 1$ is possible.  It is not hard to see that for all $X_t$ with $X_t \cap Q \neq \emptyset$ we get $i_t \in {\cal I}$. This follows from the construction of the sets $(D_i)_i$ and the fact that $ \tilde{Q} \not\subset S_\iota$ implies $ \tilde{Q} \cap \hat{S}_\iota = \emptyset$ as $\lambda \ll k$. Following the lines of \eqref{rig-eq: A2--3},  \eqref{rig-eq: weaklocA14}-\eqref{rig-eq: weaklocA15}  and using $\hat{H}_4 \le CG$ we find $\bar{R}^\iota \in SO(2)$, $\bar{c}^\iota \in \R^2$ such that
\begin{align}\label{rig-eq: iota0} 
\Vert \bar{y}_\iota - (\bar{R}^\iota\,\cdot +\bar{c}^\iota)\Vert^2_{L^2(\tilde{Q})} \le Ck^2 G, \ \  \ \ \Vert \nabla \bar{y}_\iota- \bar{R}^\iota  \Vert^p_{L^p(\tilde{Q})} \le   C\hat{H}_p
 \end{align}
 for $\iota \in {\cal I}$.  Note that for a special choice of $\iota \in {\cal I}$  (for $\iota = i$ with $i$ as considered in \eqref{rig-eq: CG}ff.) we obtain the rigid motion $\bar{R}_j \, x  + \bar{c}_j$ which we defined in Lemma \ref{rig-lemma: weaklocA}. Then arguing as in   \eqref{rig-eq: weaklocA22.2.2} and \eqref{rig-eq: weaklocA22.3.2}, in particular employing the triangle inequality and using \eqref{rig-eq: iota0}, we derive
 \begin{align}\label{rig-eq: iota} 
 \begin{split}
&\Vert (\bar{R}_j\,\cdot +\bar{c}_j) - (\bar{R}^\iota\,\cdot +\bar{c}^\iota)\Vert^2_{L^2(\tilde{Q})} \le Cm^{-2}k^2  G, \\ & \Vert \bar{R}_j - \bar{R}^\iota  \Vert^p_{L^p(\tilde{Q})} \le  Cm^{-1} \hat{H}_p
\end{split}
 \end{align}
for $\iota \in {\cal I}$.  Likewise we obtain by \eqref{rig-eq: weaklocA3} 
\begin{align}\label{rig-eq: weaklocA3_*}
\int\nolimits_{J_{\bar{y}} \cap \overline{Q}} |[\bar{y}]|^2 \,d{\cal H}^1 \le C\sum_{\iota \in {\cal I}}\int\nolimits_{J_{\bar{y}_\iota} \cap \overline{Q}} |[\bar{y}_\iota]|^2 \,d{\cal H}^1\le C \sum_{\iota \in {\cal I}}  k C_m \epsilon|\partial \hat{U}_\iota \cap Q^k_\iota(q_\iota) |_{\cal H}.
\end{align}
Here we used that all $X_t$ with  $\overline{Q} \cap X_t \neq \emptyset$ satisfy $|\partial X_t|_\infty \le \frac{k}{8}$ and thus $X_t \subset Q^k_\iota(q)$.  Now we obtain 
\begin{align*}
\Vert \nabla \bar{y} - \bar{R}_j\Vert^p_{L^p(Q)} &\le \sum\nolimits_{\iota \in {\cal I}}\Vert \nabla \bar{y}_\iota - \bar{R}_j\Vert^p_{L^p(Q)} \\ & \le C\sum\nolimits_{\iota \in {\cal I}}\big(\Vert \nabla \bar{y}_\iota - \bar{R}^\iota\Vert^p_{L^p(Q)} + \Vert \bar{R}^\iota - \bar{R}_j\Vert^p_{L^p(Q)}\big) 
\end{align*}
for $p=2,4$.  Choosing the constant in the definition of $N$ sufficiently large and recalling the definition of $G$ and  $\hat{H}_p$ (see \eqref{rig-eq: A2--3*****}) we obtain by  \eqref{rig-eq: iota0}  and  \eqref{rig-eq: iota}
\begin{align*}
\Vert \nabla \bar{y} - \bar{R}_j\Vert^p_{L^p(Q)} 
\le   CC^2_m (\gamma(N) + \delta_p(N) + \epsilon |\partial W \cap  N|_{\cal H}).
\end{align*}
Similarly,  recalling  $\lambda = mk$ we derive
\begin{align*}
\Vert \bar{y}  - (\bar{R}_j \, \cdot + \bar{c}_j)\Vert^2_{L^2( Q)} \le C\lambda^2C_m^2 \min\lbrace(\gamma(N) + \delta_4(N) + \epsilon |\partial W \cap N|_{\cal H}), \epsilon k\rbrace.
\end{align*}
Consequently, \eqref{rig-eq: corweak}(i),(ii) hold for $\bar{y}$. 

For later purposes, it is convenient to have an extension satisfying $\nabla \hat{y}(x) \in SO(2)$ for a.e. $x \in  U^H \setminus  W$.  Arguing as in \eqref{rig-eq: R,A diff}  for all components $X_t$ we find $\tilde{R}_t \in SO(2)$ such that $|\tilde{R}_t - ( R_{t} +  R_{t}A_t)|^2 \le C|A_t|^4$. Therefore, by Poincar\'e's inequality we find  for some possibly different $\tilde{c}_t \in \R^2$ 
\begin{align}\label{rig-eq: A2--3.2}
\Vert \tilde{R}_t \, \cdot + \tilde{c}_t - ( R_{t}\,(\Id + A_t)\, \cdot +  R_{t}\,c_t) \Vert^2_{L^2(X_t)} \le C |\partial X_t|_\infty^2 \  |X_t| |A_t|^4
\end{align}
for all $X_t$  and likewise passing to the trace we get
$$
\Vert \tilde{R}_t \, \cdot + \tilde{c}_t - ( R_{t}\,(\Id + A_t)\, \cdot +  R_{t}\,c_t) \Vert^2_{L^2(\partial X_t)} \le C |\partial X_t|_\infty^2 \ |\partial X_t|_{\cal H} |A_t|^4.
$$
In particular, note the the constants above do not depend on the shape of $X_t$ as the involved functions are affine. We set $\hat{y}:U^H \to \R^2$ by  $\hat{y}(x) = \tilde{R}_t \, x + \tilde{c}_t$ for $x \in X_t$ and $\hat{y} =y$ else.  First, we see that \eqref{rig-eq: corweak**} holds since ${\cal H}^1(J_{\bar{y}}) = {\cal H}^1(J_{\hat{y}})$. The definition together with \eqref{rig-eq: weaklocA3_*}  yields
\begin{align*}
\begin{split}
\int\nolimits_{J_{\bar{y}} \cap \overline{Q}} |[\hat{y}]|^2 \,d{\cal H}^1 & \le \int\nolimits_{J_{\bar{y}} \cap \overline{Q}} |[\bar{y}]|^2 \,d{\cal H}^1+ C\sum\nolimits_{X_t \cap \overline{Q} \neq\emptyset} \vert \partial X_t \vert^2_\infty  |\partial X_t|_{\cal H} |A_t|^4 \\ & 
\le \int\nolimits_{J_{\bar{y}} \cap \overline{Q}} |[\bar{y}]|^2 \,d{\cal H}^1+  Ck \sum\nolimits_{X_t \cap \overline{Q} \neq\emptyset} \vert \partial X_t \vert^2_\infty |A_t|^4 \\ & 
\le  C C_m k \sum\nolimits^4_{\iota =1} \epsilon  |\partial \hat{U}_\iota \cap   Q^k_\iota(q_\iota)|_{\cal H}.
\end{split}
\end{align*}
In the second step we used $|\partial X_t|_{\cal H} \le Ck$ which follows  from \eqref{rig-eq: new104} and \eqref{rig-eq: weaklocA6}. In the last step we used Lemma \ref{rig-lemma: A neigh} similarly as in the derivation of \eqref{rig-eq: A2--3***} and employed $s \ge \epsilon$. Using once more that $|J_{\bar{y}} \cap \overline{Q}|_{\cal H}\le \sum\nolimits^4_{\iota =1} |\partial \hat{U}_\iota \cap  Q^k_\iota(q_\iota)|_{\cal H} \le Ck$, H\"older's inequality  and \eqref{rig-eq: new104} yield
\begin{align}\label{rig-eq: weaklocA3_**}
\begin{split}
\Big(\int\nolimits_{J_{\bar{y}} \cap \overline{Q}} |[\hat{y}]|\,d{\cal H}^1 \Big)^2& \le  |J_{\bar{y}} \cap \overline{Q}|_{\cal H} \cdot \int\nolimits_{J_{\bar{y}} \cap \overline{Q}} |[\hat{y}]|^2 \,d{\cal H}^1 \\
&\le  C  C_m k^2 \sum\nolimits^4_{\iota =1} \epsilon  |\partial \hat{U}_\iota \cap  Q^k_\iota(q_\iota)|_{\cal H}.\\
&\le  CC_m^2 \lambda^2 \min \big\{\big(\gamma(N) + \delta_4(N) + \epsilon |\partial W \cap  N|_{\cal H}\big), \epsilon k \big\} .
\end{split}
\end{align}  
 Recalling $|\tilde{R}_t - ( R_{t} +  R_{t}A_t)|^2 \le C|A_t|^4$, $|A_t|\le C$ and again using \eqref{rig-eq: A2--3***}, \eqref{rig-eq: new104} we obtain
\begin{align*}
\Vert \nabla \bar{y} - \nabla \hat{y}\Vert^p_{L^p(Q)} &\le C\sum_{X_t \cap Q \neq \emptyset}  \vert  \partial X_t \vert^2_\infty |A_t|^4 \le  CC_m   ( \gamma(N) + \delta_4(N ) + \epsilon|\partial W \cap  N|_{\cal H})
\end{align*}
for $p=2,4$, and analogously by \eqref{rig-eq: A2--3.2} we get
\begin{align*}
\Vert  \bar{y} -  \hat{y}\Vert^2_{L^2(Q)} \le C\sum_{X_t \cap Q \neq \emptyset} \vert \partial X_t \vert^4_\infty |A_t|^4 \le CC_m^2 \lambda^2 \big (\gamma(N) + \delta_4(N ) + \epsilon| \partial W \cap  N|_{\cal H}\big), 
\end{align*}
where we employed $|\partial X_t|_\infty \le Ck = C\lambda m^{-1}$. Likewise we derive $\Vert  \bar{y} -  \hat{y}\Vert^2_{L^2(Q)} \le C C_m^2\lambda^2 \epsilon k$. Together with the estimates for $\bar{y}$ this shows  \eqref{rig-eq: corweak}(i),(ii). It remains to prove  \eqref{rig-eq: corweak}(iii). By \eqref{rig-eq: corweak}(i) for $p=4$, \eqref{rig-eq: linearization} and the fact that $\nabla \hat{y}(x) \in SO(2)$ for a.e. $x \in U^H \setminus W$ we find $\Vert \bar{e}_{\bar{R}_j}(\nabla \hat{y})\Vert^2_{L^2(Q)} \le  C C^2_m( \gamma(N) +  \delta_4(N) + \epsilon |N \cap \partial W|_{ \cal H})$. This together with \eqref{rig-eq: weaklocA3_**}, $|Q| \le C\lambda^2$ and H\"older's inequality yields 
$$(|E({\bar{R}_j}^T\hat{y}-\id)|(Q))^2 \le CC_m^2 \lambda^2( \gamma(N) +  \delta_4(N) + \epsilon |\partial W \cap  N|_{\cal H}).$$
Then Theorem \ref{rig-th: tracsbv} and a rescaling argument show 
\begin{align*}
\Vert {  \hat{y}}  - (\bar{R}_j \, \cdot + \bar{c}_j)\Vert^2_{L^1(\partial Q)} &\le C\lambda^{-2}\Vert  \hat{y}  - (\bar{R}_j \, \cdot + \bar{c}_j)\Vert^2_{L^1( Q)} + C(|E({\bar{R}_j}^T\hat{y}-\id)|(Q))^2 \\ & \le C\lambda^2C_m^2(\gamma(N) + \delta_4(N) + \epsilon |\partial W \cap  N|_{\cal H}).
\end{align*}
In the last step we have used H\"older's inequality and \eqref{rig-eq: corweak}(ii).  Similarly as before we also derive $\Vert {  \hat{y}}  - (\bar{R}_j \, \cdot + \bar{c}_j)\Vert^2_{L^1(\partial Q)}\le CC_m^2\lambda^2 \epsilon k.$ \eop

\section{Modification of the deformation}\label{rig-sec: sub, local}

The goal of the section is to replace the deformation by an $H^1$-function on $U_J$. In particular, we modify the deformation in such a way that the least crack length is increased. Recall $\nu = sd = \lambda m$.

\begin{lemma}\label{rig-th: global estimate}
Let $k > s ,  \epsilon > 0$ such that $l:=\frac{k}{s} = dm^{-2}$ for $m^{-1},d \in \N$ with  $m^{-1}, d \gg 1$. Then there is a constant $C>0$ such that for all $W \in  {\cal V}^{s m}_{(s,3k)}$ with $W \subset \Omega^{3k}$ and for all $y \in H^1(W)$ with $\Vert \nabla y \Vert_\infty\le C$, $\gamma$ as defined in \eqref{rig-eq: weakloc gamma} and $\delta_2, \delta_4$  as given in \eqref{rig-eq: weaklocA5}  we have the following: \\
There are sets  $U \in {\cal V}^{s\hat{m}^2}_{71k}$ and $U^H \in {\cal V}^\nu_{72k}$ with $ U, U^H \subset \Omega^{6k}$, $|U  \setminus W|=0$,  $|U^H \setminus H^\lambda(U)|=0$, $|(W \setminus U) \cap \Omega^{6k}|  + |U \setminus U^H| \le C_{u}k\Vert  U \Vert_*$ and
\begin{align}\label{rig-eq: glob1}
\Vert U \Vert_*\le (1 + C_{u}m)\Vert W \Vert_* + C\epsilon^{-1}(\gamma + \delta_4)
\end{align}
as well as a function $\tilde{y} \in H^{1}(U^H)$ such that
\begin{align}\label{rig-eq: 5 prop2}
\begin{split}
(i) & \ \ \Vert\dist(\nabla\tilde{y},SO(2) )\Vert^2_{L^2(U^H)} \le C  \min_{p=2,4} (1 + \vartheta_p^3) C_m^2 (\gamma + \delta_{p} + \epsilon \Vert W \Vert_*), \\
(ii) & \ \ \Vert \dist(\nabla \tilde{y},SO(2)) \Vert^2_{L^\infty(U^H)} \le C\bar{\vartheta}(1 + \bar{\vartheta}), \\
(iii) & \ \ \Vert \nabla y-  \nabla \tilde{y}\Vert^2_{L^2(U)} \le C   C_m^2(\gamma + \delta_2 + \epsilon \Vert W \Vert_*),\\
(iv) & \ \ \Vert\tilde{y} -  y \Vert^2_{L^2(U)} \le  C  C_m^2(1+\vartheta) \lambda^2(\gamma + \delta_4 + \epsilon \Vert W \Vert_*),
\end{split}
\end{align}
where $\bar{\vartheta} = \min \lbrace \vartheta (1+\vartheta), C_m^3\rbrace$ and $\vartheta_2 = 1$, $\vartheta_4  = \vartheta$. 
Under the additional assumption that $\Delta y = 0$ in  $W^\circ$ we get 
\begin{align}\label{rig-eq: glob4}
\Vert \nabla y-  \nabla \tilde{y}\Vert^4_{L^4(U)} \le C C_m^2 \delta_4 + C C_m^2 \vartheta (1+ \vartheta)^2 (\gamma + \delta_4 + \epsilon \Vert W \Vert_*).
\end{align}
\end{lemma}

\Proof Apply  Lemma  \ref{rig-lemma: weaklocA} to obtain sets $V \in {\cal V}^{s\hat{m}^2}_{71k}$, $U_J \in {\cal V}^\lambda$ satisfying  \eqref{rig-eq: weaklocA1.2} and \eqref{rig-eq: weaklocA2} for mappings $\bar{R}_j: U_J \to SO(2)$ and $\bar{c}_j: U_J\to  \R^2$, $j=1,\ldots,4$. We first define $U=V$ and see that  the estimate in \eqref{rig-eq: glob1}. Moreover, we recall that $\Omega^{6k} \setminus U$ is the union of rectangular components (see \eqref{rig-eq: mayblast5}). For the components $\Gamma_1(H^\lambda(V)), \ldots, \Gamma_n(H^\lambda(V))$ we let $N(\Gamma_j) \in {\cal U}^\nu$ denote the smallest rectangle with $N(\Gamma_j) \supset X_j$, where as before $X_j$ denotes the component corresponding to $\Gamma_j(H^\lambda(V))$. 

As $\frac{\nu}{\lambda}= m$, we find $|\pi_i \partial N(\Gamma_j)| \le |\pi_i \Gamma_j(H^\lambda(V))| + C_u m|\Gamma_j(H^\lambda(V))|_\infty$ for $i=1,2$.  Arguing similarly as in the construction of  \eqref{rig-eq: new111} we have that $N(\Gamma_{j_1}) \setminus N(\Gamma_{j_2})$ is connected for $1 \le j_1,j_2 \le n$. We apply Lemma \ref{rig-lemma: modifica}(i)  to obtain pairwise disjoint, connected sets $(X'_j)^n_{j=1}$ such that $\bigcup^n_{j=1} \overline{N(\Gamma_j)} = \bigcup^n_{j=1} \overline{X'_j}$ and define 
$$ U^H = \Omega^{6k} \setminus \bigcup\nolimits^n_{j=1} X'_j.$$
It is not hard to see that $U^H \in {\cal V}^\nu_{72k}$. Moreover, we find $U^H \subset H^\lambda(U)$ up to a set of negligible measure and recalling \eqref{rig-eq: mayblast5} we obtain $(U^H)^\circ \subset U_J$.   For later we also observe that  
\begin{align}\label{rig-eq: mayblast 10}
\Vert U^H \Vert_* \le (1+C_u m) \Vert H^\lambda(U) \Vert_*.
\end{align}
This also implies $|U \setminus U^H| \le C_uk \Vert U \Vert_*.$

Let $T_j = \bigcup\nolimits_{p \in I^\lambda_j(\Omega^{3k})} Q_j^{\frac{3}{4}\lambda}(p)$ and define a partition of unity $(\eta_j)^4_{j=1}$ with $\eta_j \in C^\infty(U_J,[0,1])$, $\text{supp}(\eta_j) \subset T_j$ and $\Vert \nabla\eta_j\Vert_\infty \le \frac{C}{\lambda}$. Define $\tilde{y}: U_J \to \R^2$ by
\begin{align*}
\tilde{y}(x) = \sum\nolimits^4_{j=1} \eta_j(x) (\bar{R}_j \,x + \bar{c}_j)                                 
\end{align*}
and observe that $\tilde{y} \in H^1(U_J)$ as the functions $\bar{R}_j, \bar{c}_j$ are constant on each $Q^{\lambda}_j(p)$, $p \in I^{\lambda}_j(U_J)$. The derivative reads as 
\begin{align}\label{rig-eq: derivative}
\begin{split}
\nabla \tilde{y}(x)& = \sum\nolimits^4_{j=1} \big( \eta_j(x) \bar{R}_j   + ( \bar{R}_j\,x + \bar{c}_j) \otimes\nabla \eta_j(x)\big).
\end{split}  
\end{align}
Since $\sum^4_{j=1} \nabla \eta_j = 0$ we find
$$\nabla \tilde{y}(x) = \bar{R}_1 + \sum\nolimits^4_{j=2} \big( \eta_j(x) (\bar{R}_j - \bar{R}_1)   + ( \bar{R}_j\,x + \bar{c}_j - (\bar{R}_1\,x + \bar{c}_1 )) \otimes\nabla \eta_j(x)\big). $$
First,  we compute by  \eqref{rig-eq: weaklocA9} 
\begin{align*}
\Vert\nabla \tilde{y} - \bar{R}_1\Vert^4_{L^4(U_J)} & \le   C\sum^4_{j=2} \Big( \Vert \bar{R}_j - \bar{R}_1 \Vert^4_{L^4(U_J)} + \frac{1}{\lambda^4}\Vert \bar{R}_j\,\cdot + \bar{c}_j - (\bar{R}_1\,\cdot + \bar{c}_1)\Vert^4_{L^4(U_J)} \Big) \\
& \le  C\sum^4_{j=2} \Big(  \Vert \bar{R}_j - \bar{R}_1 \Vert^4_{L^4(U_J)} + \frac{ \bar{\vartheta}}{\lambda^2}\Vert \bar{R}_j\,\cdot + \bar{c}_j - (\bar{R}_1\,\cdot + \bar{c}_1)\Vert^2_{L^2(U_J)} \Big),
\end{align*}
where $\bar{\vartheta} = \min \lbrace \vartheta (1+\vartheta), C_m^3\rbrace$. 
By \eqref{rig-eq: linearization} we find $\bar{e}_{\bar{R}_1}(\bar{R}_j) \le C |\bar{R}_j - \bar{R}_1|^2$ and thus 
\begin{align*}
\Vert \bar{e}_{\bar{R}_1}(\nabla \tilde{y})\Vert^2_{L^2(U_J)}& \le  C \sum^4_{j=2} \Big( \Vert \bar{e}_{\bar{R}_1}(\bar{R}_j)\Vert^2_{L^2(U_J)} + \frac{1}{\lambda^2}\Vert \bar{R}_j\,\cdot + \bar{c}_j - (\bar{R}_1\,\cdot + \bar{c}_1)\Vert^2_{L^2(U_J)} \Big) \\
& \le  C\sum^4_{j=2} \Big( \Vert \bar{R}_j - \bar{R}_1 \Vert^4_{L^4(U_J)} + \frac{1}{\lambda^2}\Vert \bar{R}_j\,\cdot + \bar{c}_j - (\bar{R}_1\,\cdot + \bar{c}_1)\Vert^2_{L^2(U_J)} \Big).  
\end{align*}
 Again using \eqref{rig-eq: linearization} and \eqref{rig-eq: weaklocA2}(iii),(iv) we derive
 $$\Vert \dist(\nabla \tilde{y},SO(2)) \Vert^2_{L^2(U_J)} \le C  (1 +  \vartheta^3) C_m^2(\gamma + \delta_4 + \epsilon \Vert W \Vert_*). $$
Similarly, we get
$$  \Vert\nabla \tilde{y} - \bar{R}_1\Vert^2_{L^2(U_J)} \le C\sum^4_{j=2} \Big( \Vert \bar{R}_j - \bar{R}_1 \Vert^2_{L^2(U_J)} + \frac{1}{\lambda^2}\Vert \bar{R}_j\,\cdot + \bar{c}_j - (\bar{R}_1\,\cdot + \bar{c}_1)\Vert^2_{L^2(U_J)} \Big) $$ 
and thus we find by \eqref{rig-eq: weaklocA2}(iii),(iv)
\begin{align*}  
  \Vert \dist(\nabla \tilde{y},SO(2)) \Vert^2_{L^2(U_J)} \le C  C_m^2 (\gamma + \delta_2 + \epsilon \Vert W \Vert_*), 
  \end{align*}
  where we used that $\delta_4 \le C\delta_2$. This gives \eqref{rig-eq: 5 prop2}(i) as $(U^H)^\circ \subset U_J$. Likewise, we may replace the $L^2, L^4$-norms in the above estimates by the $L^\infty$-norm. Consequently, by  \eqref{rig-eq: weaklocA9} we  obtain $\Vert\nabla \tilde{y} - \bar{R}_1\Vert^4_{L^\infty(U_J)} \le C \bar{\vartheta}(1 +\bar{\vartheta}) $ and  $\Vert \bar{e}_{\bar{R}_1}(\nabla \tilde{y})\Vert^2_{L^\infty(U_J)} \le C\bar{\vartheta}$ which then implies $\Vert \dist(\nabla \tilde{y},SO(2)) \Vert^2_{L^\infty(U_J)} \le C \bar{\vartheta}(1 + \bar{\vartheta})$. It remains to show \eqref{rig-eq: 5 prop2}(iii),(iv) and \eqref{rig-eq: glob4}.  By \eqref{rig-eq: weaklocA2}(i) and the fact that  $U=V$  we obtain
\begin{align*}
\Vert\tilde{y} -  y \Vert^2_{L^2(U)} & \le  \sum\nolimits^4_{j=1} C\Vert y - (\bar{R}_j\,\cdot + \bar{c}_j)\Vert^2_{L^2(U)}  \le  C C_m^2  \lambda^2 (1+ \vartheta) (\gamma + \delta_4 + \epsilon \Vert W \Vert_*).
\end{align*} 
By \eqref{rig-eq: derivative} and the fact that $\sum^4_{j=1}  \eta_j =1$, $\sum^4_{j=1} \nabla \eta_j = 0$ we derive
 $$\nabla y(x) -  \nabla \tilde{y}(x) = \sum\nolimits^4_{j=1} \big( \eta_j(x) (\nabla y(x) - \bar{R}_j)   + ( y(x) - (\bar{R}_j\,x + \bar{c}_j)) \otimes\nabla \eta_j(x)\big).$$
Therefore, by \eqref{rig-eq: weaklocA2}(i)(ii) for $p=2$ we get 
\begin{align*}
\Vert\nabla \tilde{y} - \nabla y \Vert^2_{L^2(U)} & \le   C\sum\nolimits^4_{j=1} \Big( \Vert \nabla y - \bar{R}_j \Vert^2_{L^2(U)} + \frac{1}{\lambda^2}\Vert y - (\bar{R}_j\,\cdot + \bar{c}_j)\Vert^2_{L^2(U)} \Big) \\
& \le  C C_m^2 (\gamma + \delta_2 + \epsilon \Vert W \Vert_*) ,
\end{align*} 
where we used that $\delta_4 \le C\delta_2$. Finally, in the case that $\Delta y = 0$ in $W^\circ$ we obtain by \eqref{rig-eq: weaklocA2}(i)(ii) for $p=4$ and \eqref{rig-eq: weaklocA9.2} 
 \begin{align*}
\Vert\nabla \tilde{y} - \nabla y \Vert^4_{L^4(U)} & \le   C\sum\nolimits^4_{j=1} \Big( \Vert \nabla y - \bar{R}_j \Vert^4_{L^4(U)} + \frac{  \vartheta (1+ \vartheta)}{\lambda^2}\Vert y - (\bar{R}_j\,\cdot + \bar{c}_j)\Vert^2_{L^2(U)} \Big) \\
& \le  C C_m^2 \delta_4 + C C_m^2 \vartheta  (1+ \vartheta)^2 (\gamma + \delta_4 +\epsilon \Vert W \Vert_*).
\end{align*} 
 \eop

\section{SBD-rigidity up to small sets}\label{rig-sec: sub, proof}

In this section we prove a slightly weaker version of the rigidity estimate given in Theorem \ref{rig-th: rigidity} and postpone the proof of the general version to the next section.  Recall definition \eqref{rig-eq: SBVfirstdef}.

\begin{theorem}\label{rig-th: rigidity2}
Let $\Omega \subset \R^2$ open, bounded with Lipschitz boundary. Let $M>0$ and $0 < \eta, \rho, h_* \ll 1$.  Let $q \in \N$ sufficiently large.  Then there are constants $C_1=C_1(\Omega,M,\eta)$, $C_2= C_2(\Omega,M,\eta,\rho,h_*,q)$  and a universal constant $c>0$ such that the following holds  for $\eps > 0$ small enough:\\ For each $y \in SBV_M(\Omega)$ with ${\cal H}^1(J_y) \le M$  and $\int_\Omega \dist^2(\nabla y,SO(2) ) \le M\eps$, there is  a set $\Omega_y \in {\cal V}^{\hat{s}}_{ c\rho^{q-1}}$, $\hat{s}>0$, with $\Omega_y \subset \Omega$, $|\Omega\setminus\Omega_y| \le C_1\rho$,  a modification  $\tilde{y} \in H^1(\Omega_y) \cap SBV_{cM}(\Omega_y)$ with $\Vert y - \tilde{y}\Vert^2_{L^2(\Omega_y)}  + \Vert \nabla y - \nabla \tilde{y}\Vert^2_{L^2(\Omega_y)}\le C_1\eps \rho$, a partition $(P_i)_i$ of $\Omega_y$ and for each $P_i$ a corresponding rigid motion $R_i \, x +c_i$, $R_i \in SO(2)$ and $c_i \in \R^2$, such that the function $u: \Omega \to \R^2$ defined by
\begin{align}\label{rig-eq: u def} 
u(x) := \begin{cases} \tilde{y}(x) - (R_i\,\,x +c_i) & \ \ \text{ for } x \in P_i \\
                      0                      & \ \ \text{ else} \end{cases}
\end{align}
satisfies
\begin{align}\label{rig-eq: main properties}
\begin{split}
(i) & \ \, \Vert \Omega_y\Vert_* \le (1 + C_1h_*){\cal H}^1(J_y) + C_1\rho, \  \ \, (ii) \,  \ \Vert u\Vert^2_{L^2(\Omega_y)} \le C_2\eps, \\
(iii) & \ \, \sum\nolimits_i \Vert  e(R^T_i \nabla u)\Vert^2_{L^2(P_i)} \le C_2\eps, \  \ \ \, \  \ \ \  \  (iv)  \ \, \Vert \nabla u\Vert^2_{L^2(\Omega_y)} \le C_2\eps^{1-\eta}.
\end{split}
\end{align}  
\end{theorem}
We divide the proof into three steps. We begin with a version where the least crack length is almost of macroscopic size. Afterwards, we assume that the jump set consists only of a finite number of cracks of arbitrary size. Finally, we treat the general case applying a suitable approximation argument.

  In what follows, constants indicated by $C_1$ only depend on $M,\eta,\Omega$. Generic constants $C$ may additionally depend on $h_*$. All constants  do not depend on $\rho$ and $q$ unless stated otherwise. As we will eventually let $h_* \sim \rho$ in Section \ref{rig-sec: sub, proof-main}, it is essential that the constant in \eqref{rig-eq: main properties}(i) does not depend on $h_*$.

\subsection{Step 1: Deformations with least crack length}\label{rig-sec: subsub,  step1}

We first treat the case that the least crack length is almost of macroscopic size.

\begin{theorem}\label{rig-thm: V1}
Theorem \ref{rig-th: rigidity2} holds under the additional assumption that there is an $\tilde{\Omega}_y \subset \Omega^{s}$, $\tilde{\Omega}_y \in {\cal V}^{s}_{\rho^{q-1}}$ for some $s \ge  \rho^{q-1}\eps^{\frac{\eta}{8}}$ such that $y \in H^1(\tilde{\Omega}_y)$, $\Vert \tilde{\Omega}_y\Vert_* \le (1 + C_1h_*){\cal H}^1(J_y) + C_1\rho$ and $|\Omega \setminus \tilde{\Omega}_y| \le C_1\rho$ for a constant $C_1 = C_1(\Omega, M,\eta)$.
\end{theorem}

\Proof Let $y \in H^1(\tilde{\Omega}_y)$ be given. Let $\rho$ and define $\varrho= \rho^{q}$ for some  $q \in \N$, $q\ge 2$ large enough to be specified in the proof of Theorem \ref{rig-th: rigidity}  (see Section \ref{rig-sec: sub, proof-main}). Assume without restriction $\rho^{-1} \in \N$ large.   We apply Theorem \ref{rig-thm: harmonic} and consider the harmonic part $w$ of $y$ satisfying 
\begin{align}\label{rig-eq: w energy}
\begin{split}
&\Vert \nabla y - \nabla w\Vert^2_{L^2(\tilde{\Omega}_y)} \le C\Vert \dist(\nabla y,SO(2))\Vert^2_{L^2(\tilde{\Omega}_y)} \le C\eps, \\
&\Vert \nabla y - \nabla w\Vert^4_{L^4(\tilde{\Omega}_y)} \le C\Vert \dist(\nabla y,SO(2))\Vert^4_{L^4(\tilde{\Omega}_y)} \le C\eps.
\end{split}
\end{align}
In the last inequality we used $\Vert \nabla y\Vert_\infty \le M$.  Let $k= \varrho \rho^{-1} = \rho^{q-1}$. Apply Lemma \ref{rig-lemma: weaklocR} on $ \tilde{\Omega}_y \cap \Omega^k$ for the function $w$ and $\epsilon= \hat{c}\rho^{-1}\eps$, $m=\rho$,  where $\hat{c}>0$ is sufficiently large. (Possibly passing to a smaller $s$ we can assume that $k \eps^{\frac{\eta}{8}}\le s \ll k=\rho^{q-1}$.) We find a set $W \subset \Omega^{3k}$, $W \in {\cal V}^{sm}_{(s,3k)}$ such that 
\begin{align}\label{rig-eq: V1.41}
\Vert W \Vert_* \le (1+ C_1\rho)\Vert \tilde{\Omega}_y \Vert_* + C\epsilon^{-1} \eps \le (1+ C_1\rho)\Vert \tilde{\Omega}_y \Vert_* + \rho
\end{align}
by \eqref{rig-eq: weaklocR2} and  $|(\tilde{\Omega}_y \setminus W) \cap \Omega^{3k}| \le C_1k \le C_1\rho$.  (Here and in the following we choose the constant $\hat{c} = \hat{c}(h_*)$  always  larger then the constant $C$.)  Moreover, there are mappings $\hat{R}_i:  W^\circ \to SO(2)$, $i=1,\ldots,4$, which are constant on  the connected components of $Q^{k}_i(p) \cap  W^\circ$, $p \in I^{k}_i(\Omega)$, such that by \eqref{rig-eq: weaklocR1}(i) for $i=1,\ldots,4$
\begin{align}\label{rig-eq: V1.11.2}
\Vert\nabla y - \hat{R}_i\Vert^2_{L^2(W)} \le C\eps + C\Vert\nabla w - \hat{R}_i\Vert^2_{L^2(W)} \le Cl^4\eps \le  C\eps^{1-\eta},
\end{align}
where $l=k s^{-1} \le C\eps^{-\frac{\eta}{8}}$. Moreover, as $\vartheta = l^{9} C_m^{2} s^{-1}\eps \le C(\rho)s^{-10}\eps \le C(\rho)\eps^{1-\frac{5}{4}\eta} \le  1$ for $\eta, \eps$ small enough  (recall \eqref{rig-eq: vartheta def}) we also get
\begin{align}\label{rig-eq: V1.11}
\Vert\nabla y - \hat{R}_i\Vert^4_{L^4(W)}\le C\eps + C\Vert\nabla w - \hat{R}_i\Vert^4_{L^4(W)} \le  C\eps
\end{align}
by \eqref{rig-eq: weaklocR1}(ii).  Now we apply Corollary \ref{rig-cor: weakA} on $W \subset \Omega^{3k}$ for $k=\rho^{q-1}$, $ \lambda= 3\varrho$, $m={3}\rho$  and $\epsilon= \hat{c}\rho^{-1}\eps$.  We obtain a set  $\Omega_y  \in {\cal V}^{s\hat{m}}_{9k}$   with $\Omega_y \subset \Omega^{5k}$,  $|\Omega_y \setminus \tilde{\Omega}_y| = 0$ such that by \eqref{rig-eq: weaklocA1}, \eqref{rig-eq: V1.41}  and \eqref{rig-eq: V1.11} we find 
\begin{align}\label{rig-eq: V1.7}
\Vert \Omega_y \Vert_* &\le (1+ C_1\rho) \Vert W \Vert_* + C\epsilon^{-1} \eps  \le (1+ C_1h_*){\cal H}^1(J_y) + C_1\rho
\end{align}
and  $| (\tilde{\Omega}_y \setminus \Omega_y) \cap \Omega^{5k}| \le C_1k$. This together with the assumption $|\Omega \setminus \tilde{\Omega}_y| \le C_1\rho$ and the fact that $|\Omega \setminus \Omega^{5k}| \le C(\Omega)k$ yields $|\Omega \setminus \Omega_y| \le C_1\rho$. Moreover, there is a set $\Omega^H_y \in {\cal V}^\lambda$ with $H^{\lambda}(\Omega_y) \subset \Omega_y^H$ and mappings $\bar{R}_j: \Omega^H_y \to SO(2)$, $\bar{c}_j: \Omega^H_y \to \R^2$ being constant on $Q^{3\varrho}_j(p)$, $p \in I^{3\varrho}_j(\Omega^{ 3k})$, and an extension  $\hat{y} \in SBV_M(\Omega^H_y,\R^2)$ such that by  \eqref{rig-eq: corweak}(ii) we derive 
\begin{align}\label{rig-eq: V1.10}
\Vert \hat{y} - (\bar{R}_j \, \cdot + \bar{c}_j)  \Vert^2_{L^2(\Omega^H_y)} \le C\varrho^{2}  \rho^{-2} C_\rho^4(\eps + \epsilon \Vert W \Vert_*)  \le C \rho^{2q- 3} C_\rho^4\eps
\end{align}
where $C_\rho = C_{ \frac{m}{3}}$ is the constant defined in \eqref{rig-eq: vartheta def}.  Here we used that each $x \in W$ is contained in at most $\sim \rho^{-2}$ different neighborhoods $N(Q)$ considered in Corollary \ref{rig-cor: weakA}. Moreover, the constant $\hat{c}$ was absorbed in $C$. Similarly,  recalling $\vartheta \le 1$ we get by \eqref{rig-eq: weaklocA2}(iii),(iv), \eqref{rig-eq: corweak}(i)  and \eqref{rig-eq: V1.11.2}, \eqref{rig-eq: V1.11} 
\begin{align}\label{rig-eq: V1.9.2}
\begin{split}
&\Vert \nabla  \hat{y} - \bar{R}_j   \Vert^2_{L^2(\Omega^H_y)} + \Vert \bar{R}_{j_1} - \bar{R}_{j_2}   \Vert^2_{L^2(\Omega^H_y)} \le C\rho^{-3}C_\rho^2\eps^{1-\eta}, \\
& \Vert \nabla  \hat{y} - 
\bar{R}_j  \Vert^4_{L^4(\Omega^H_y)} + \Vert \bar{R}_{j_1} - \bar{R}_{j_2}   \Vert^4_{L^4(\Omega^H_y)}\le C\rho^{-3}C_\rho^2\eps,  \\
& \Vert (\bar{R}_{j_1} \, \cdot + \bar{c}_{j_1})  - (\bar{R}_{j_2} \, \cdot + \bar{c}_{j_2})\Vert^2_{L^2(\Omega^H_y)}\le C\rho^{2q-3} C_\rho^2\eps, 
\end{split}
\end{align}
for $j=1,\ldots,4$ and $1 \le j_1, j_2 \le 4 $. 

Denote  the connected components of  $(\Omega^H_y)^\circ \in {\cal U}^{3\varrho}$ by  $(P^H_i)_i$ and define  $P_i = P_i^H \cap \Omega_y$.  Let $J_i \subset I^{\varrho}(\Omega)$ be the index set such that $Q^{ \varrho}(p) \subset  P^H_i$ for all $p \in J_i$. We now estimate the variation of the rigid motions defined on these squares.  Let $Q_1 = Q^{\varrho} (p_{1})$, $Q_2 = Q^{\varrho} (p_{2})$ for $p_1,p_2 \in J_i$ such that $\overline{Q_1} \cap \overline{Q_2} \neq \emptyset$. Let $R_t = \bar{R}_4|_{Q_t}$ and $c_t = \bar{c}_4|_{Q_t}$ for $t=1,2$. Then  we find some  $j=1,\ldots,4$ such that $\bar{R}_j$ is constant on $Q_1 \cup Q_2$ and thus  $\varrho^2|R_1 - R_2|^p \le C\sum_{t=1,2} \Vert \bar{R}_j - R_{t}\Vert^p_{L^p(Q_1 \cup Q_2)}$ for $p=2,4$. Using the arguments in \eqref{rig-eq: A difference} and \eqref{rig-eq: A,c difference2}  we get 
\begin{align}\label{rig-eq: rig2}
\begin{split}
\varrho^4|R_1 - R_2|^2 + &\Vert (R_1 - R_2)\,\cdot + c_1 - c_2\Vert^2_{L^2( Q_1 \cup Q_2)} \\
&\le C\sum\nolimits_{t=1,2}  \Vert (\bar{R}_{j}\, \cdot + \bar{c}_j) - (R_{t}\, \cdot + c_t)\Vert^2_{L^2(Q_1 \cup Q_2)}.
\end{split}
\end{align}
Consequently, considering chains as in \eqref{rig-eq: path} and  \eqref{rig-eq: A diff}, respectively, following the arguments in the proof of Lemma \ref{rig-lemma: weak rig} and \eqref{rig-eq: A diff} and recalling Remark \ref{rig-rem: 1}(ii), we obtain $R_i \in SO(2)$, $c_i \in \R^2$ such that
\begin{align*}
\Vert \hat{y} - (R_i\, \cdot + c_i)\Vert^2_{L^2(P^H_i)}& \le C\Vert \hat{y} - (\bar{R}_4\, \cdot + \bar{c}_4)\Vert^2_{L^2(P^H_i)}\\
& \ \ \ \ + C\varrho^{-8}\sum\nolimits_{1 \le j_{1},j_2 \le 4}\Vert (\bar{R}_{j_1} \, \cdot + \bar{c}_{j_1})  - (\bar{R}_{j_2} \, \cdot + \bar{c}_{j_2})\Vert^2_{L^2(P^H_i)},\\
\Vert \nabla  \hat{y} - R_i\Vert^p_{L^p(P^H_i)}& \le  C\Vert \nabla \hat{y} - \bar{R}_4\Vert^p_{L^p(P^H_i)} \\
&  \ \ \ \ + C\varrho^{-2p}\sum\nolimits_{1 \le j_{1},j_2 \le 4}\Vert \bar{R}_{j_1} - \bar{R}_{j_2}   \Vert^p_{L^p(P^H_i)}, \ \ \ p=2,4.
\end{align*}
In the first estimate we used H\"older's inequality (cf. \eqref{rig-eq: A diff}). Summing over all connected components, \eqref{rig-eq: V1.10} and \eqref{rig-eq: V1.9.2} implies
\begin{align}\label{rig-eq: compare}
\begin{split}
&\sum\nolimits_i\Vert  \hat{y} - (R_i\, \cdot + c_i)\Vert^2_{L^2(P^H_i)} \le C(\rho,q)\eps,\\
& \sum\nolimits_i \Vert \nabla  \hat{y} - R_i\Vert^4_{L^4(P^H_i)} \le C(\rho,q)\eps, \ \, \sum\nolimits_j \Vert \nabla  \hat{y} - R_i\Vert^2_{L^2(P^H_i)} \le C(\rho,q)\eps^{1- \eta}
\end{split}
\end{align}
for $C(\rho,q)$ large enough. Defining $u$ as in \eqref{rig-eq: u def} (for $\tilde{y} = y$)  and taking also  \eqref{rig-eq: V1.7} into account, we immediately get \eqref{rig-eq: main properties}(i)(ii),(iv). Finally, \eqref{rig-eq: main properties}(iii) is a consequence of the linearization formula \eqref{rig-eq: linearization2} and \eqref{rig-eq: compare}. \eop

\subsection{Step 2: Deformations with a finite number of cracks}\label{rig-sec: subsub,  step2}

We now prove a version where the crack set consists of a finite number of components. We first assume that each crack is at least of atomistic size. The strategy will be to establish an estimate of the form \eqref{rig-eq: V1.11.2} and \eqref{rig-eq: V1.11} by iterative modification of $y$ according to Lemma \ref{rig-th: global estimate}. 

First, we introduce some notation and derive preliminary estimates. Let $\rho>0$, set $\varrho= \rho^{q}$ and assume without restriction $\rho^{-1} \in \N$ large.  As before we assume $\Vert \dist(\nabla y, SO(2))\Vert^2_{L^2(\Omega)}\le C\eps$. Choose $t^{-1}\in \N$ such that $t\le \rho$ and set $t_j = t^{j+1}$. By Remark \ref{rig-rem: z}(i) we can assume that $T:= t^{z+18} \le C^{-2}_{t} t^{18}$ for $z \in \N$ sufficiently large (recall \eqref{rig-eq: vartheta def} for the definition of $C_t$).  Moreover, set $T_j = T^{j+1}$. Let $\tilde{\Omega}_y \subset \Omega^s$ for some $s>0$ be given. Let 
\begin{align}\label{rig-eq: B def}
B_j =  \big(\Vert  \tilde{\Omega}_y \Vert_* + C_* \rho\big)\ \cdot \text{\scriptsize $\sum\nolimits$}^{j-1}_{i=0} t^{i} \ \cdot \Pi^{j-1}_{i=0} (1+ C_* t^{i+1})
\end{align}
and $B = \lim_{j\to\infty} B_j$ for a constant $C_*  = C_*(M,\eta,\Omega) \ge 1$ to be specified below. Furthermore, let $P= \hat{c}^2 (1  + \rho^{-1} B)$  for $\hat{c}=\hat{c}(h_*)$ sufficiently large. Set $s_0 = \kappa \eps$ for  $\kappa$ sufficiently large, let  $\epsilon_0 =   \hat{c}^2\rho^{-1} \eps$ and  subsequently define $\epsilon_{j+1} = PT^{-1}_j\epsilon_j$. We set $r = \frac{1}{ 18}$, $\omega = \frac{\eta}{36}$ for notational convenience and for $j \ge 0$ we define 
\begin{align}\label{rig-eq: l def}
d_j = \Big\lfloor\min \Big\{\Big(\frac{s_j}{\epsilon_j}\Big)^r , \eps^{-\omega}\Big\}\Big\rfloor, 
\end{align}
where $s_j = s_0 \Pi^{j-1}_{i=0} d_i$. In accordance with Sections \ref{rig-sec: sub, first-weak}, \ref{rig-sec: sub, local} we also define \begin{align}\label{rig-eq: lds}
l_j = d_j t^{-2}_j, \ \ \ \ \ \lambda_j=s_j d_j t^{-1}_j,  \ \ \ \ \ k_j = s_j l_j.
\end{align}
As noted before, $d_j$ describes the increase of the minimal distance of different cracks and $PT^{-1}_j$ will be the factor of energy increase. Below we will show that indeed  $d_j \gg 1$ for all $0 \le j \le J^*$, where
$$J^* = \lceil\text{log}_{1+r}(\text{log}_T \eps^{\omega})) +\tfrac{1}{\omega}\rceil.$$
One of the main reasons why the iterative application of Lemma \ref{rig-th: global estimate} works is the fact that $d_j$ increases much faster than $PT^{-1}_j$. We define the quotient $q_j := \frac{d_j  }{P T_j^{-1}}$ and observe $q_0 = \frac{d_0 T_0}{P} = TP^{-1} (s_0 \epsilon_0^{-1})^r$  for $\eps$ sufficiently small. Recalling \eqref{rig-eq: B def} and the definition $s_0 = \kappa \eps$,  $\epsilon_0 =   \hat{c}^2 \rho^{-1} \eps$  we can first choose $T = T(\rho,h_*)$ so small and then $\kappa  = \kappa(T,\rho, h_*, \bar{z})  $ so large that 
\begin{align}\label{rig-eq: V3.3}
q_0 T^{1/r} \ge T^{- \bar{z}}  \ge T^{-1} \ge  \hat{c}^4 P^2 > 1
\end{align}  
 for $\bar{z} \in \N$ to be specified below. For the third inequality we used the fact that $P\le C$ for some $C=C(C_*,\rho, h_*,  M)$ independent of $T$. We find
\begin{align}\label{rig-eq: V3.2}
q_j =  T^{-1/r} (q_0 T^{1/r})^{(1+r)^j}
\end{align}
for $j \le \hat{J}$, where  $\hat{J} \in \N$ is the largest index such that $\frac{s_j}{\epsilon_j} \le \eps^{-\frac{\eta}{2}}$ for all $j \le \hat{J}$. Indeed,  we first note that the formula is trivial for $j=0$. Assume \eqref{rig-eq: V3.2} holds for $j  \le \hat{J}-1$, then we compute 
\begin{align*}
q_{j+1} = \frac{T_{j+1}}{P} \Big(\frac{s_{j+1}}{\epsilon_{j+1}}\Big)^r = \frac{T_{j+1}}{P} \Big(\frac{s_{j}d_{j}}{PT^{-1}_{j}\epsilon_{j}}\Big)^r = \frac{q^r_{j} T_{j+1}  }{P  }\Big(\frac{s_{j}}{\epsilon_{j}}\Big)^r= \frac{q^r_{j} d_j { T_{j}} T}{P   } = T q^{1+r}_{j} 
\end{align*}  
which gives \eqref{rig-eq: V3.2} for $j+1$, as desired.  In particular, taking \eqref{rig-eq: V3.3} into account, \eqref{rig-eq: V3.2} implies $q_j > 1$ and thus $d_j  = q_j PT^{-1}_j \gg 1$ for all $j \le \hat{J}$. For $\hat{J} < j \le J^*$ we get $d_j = \eps^{-\omega}$. In fact, using \eqref{rig-eq: V3.3} and $ \epsilon_0 \le  \hat{c}^2 t^{-1}\eps$  we observe for $C$ sufficiently large
\begin{align}\label{rig-eq: V3.4_*}
\begin{split}
\epsilon_j & = \epsilon_0 \Pi^{j-1}_{i=0} (PT^{-1}_i) \le  \hat{c}^{-2} \epsilon_0  \Pi^{j-1}_{i=0} (T^{-(i+1)}T^{-\frac{1}{2}} )\le  \hat{c}^{-2} \epsilon_0 T^{-\frac{1}{2}(j+1)^2}  \\ 
& \le \eps  T^{-C- [\text{log}_{1+r}(\text{log}_T\eps^{\omega})]^{2}}= \eps  o\big(T^{- \text{log}_{T^{-1}}\eps^{-\omega}}\big) = \eps \cdot o(\eps^{-\omega})
\end{split}
\end{align}
for $\eps \to 0$ for all $1 \le j \le J^*$. Consequently, if $\frac{s_j}{\epsilon_j} \ge \eps^{-\frac{\omega}{r}} = \eps^{-\frac{\eta}{2}}$, then $d_j = \eps^{-\omega}$, $PT^{-1}_j = o(\eps^{-\omega})$  (see \eqref{rig-eq: V3.4_*}) and thus $\frac{s_{j+1}}{\epsilon_{j+1}}  = \frac{d_j s_j}{ PT^{-1}_j \epsilon_{j}} \ge \eps^{-\frac{\omega}{r}}$. This then implies $d_{j} = \eps^{-\omega}$ for all $\hat{J} < j \le J^*$.

We introduce   $\vartheta_j = s_j^{-1} \epsilon_j l_j^{9} C^{2}_{t_j}$ (recall definition \eqref{rig-eq: vartheta def} and $l_j = d_jt^{-2}_j$) and close the preparations by showing that 
\begin{align}\label{rig-eq: vartheta char}
\vartheta_{j}  \le \frac{\epsilon_0}{ \hat{c}^2\epsilon_{j+1}} T_{j} \ \ \text{for} \ \ 0 \le j \le J^*. 
\end{align}
 This particularly implies $\vartheta_j \le 1$ for all $j$ as $\epsilon_j \ge \epsilon_0$ for all $j$. By \eqref{rig-eq: l def}-\eqref{rig-eq: V3.2}  we  obtain
\begin{align}\label{rig-eq: V3.4}
\begin{split}
s_j& \ge \epsilon_j \eps^{-\frac{\eta}{2}} \ \ \text{ or } \ \ s_{j} = \epsilon_j d^{1/r}_{j} \ge  \epsilon_j q^{1/r}_{j} \ge \epsilon_jT^{ -\frac{\bar{z}}{r} (1+r)^{j}} \ge \epsilon_j T^{-9 (j+1)^2}.
\end{split}
\end{align}
for all $0 \le j \le J^*$. The last step holds for  $\bar{z} \in \N$ sufficiently large as $\lim_{j \to \infty} \frac{1}{r}(1+r)^{j} (9 (j+1)^2)^{-1} = \infty$.  Similarly as in \eqref{rig-eq: V3.4_*} we see that $T^{-9 (j+1)^2} = o(\eps^{-\omega})$ for $j \le J^*$ as $\eps \to 0$.  Since $\eps^{-\omega} =o(\eps^{-\frac{\eta}{2}})$, we find $s_j \ge \epsilon_j T^{-9 (j+1)^2}$ for all $0 \le j \le J^*$. Therefore,  we derive by \eqref{rig-eq: l def}, \eqref{rig-eq: V3.3},  the first line of \eqref{rig-eq: V3.4_*}  and $r = \frac{1}{18}$
\begin{align*}
\vartheta_{j} \epsilon_{j+1} &= s_j^{-1} \epsilon_j d_j^{9} t_j^{-18} C_{t_j}^2 \ PT^{-1}_j \epsilon_j \le s_j^{-\frac{1}{2}}  \epsilon_j^{\frac{3}{2}}  \hat{c}^{-2} T^{-3}_j \le   \hat{c}^{-2}\epsilon_0 T^{4(j+1)^2} T^{-3}_j \le  \hat{c}^{-2}\epsilon_0 T_{j}   
\end{align*}
for all $0 \le j \le J^*$, as desired.  In the second step we used $ C^2_{t_j} t_j^{-18} \le T^{-1}_j$ and $P \le T^{-1}_j$.  Recall the definition of $\kappa$ and $k_0$ above (see \eqref{rig-eq: lds} and \eqref{rig-eq: V3.3}).

\begin{theorem}\label{rig-thm: V2}
Theorem \ref{rig-th: rigidity2} holds under the additional assumption that there is an $\tilde{\Omega}_y \subset \Omega^{ s}$, $\tilde{\Omega}_y \in {\cal V}^{s}_{ k_0}$ for some $s \ge \kappa \eps$, such that $y \in H^1(\tilde{\Omega}_y)$, $\Vert \tilde{\Omega}_y\Vert_* \le (1 + C_1h_*){\cal H}^1(J_y) + C_1\rho$ and $|\Omega \setminus \tilde{\Omega}_y| \le C_1\rho$ for a constant $C_1 = C_1(\Omega, M,\eta)$.
\end{theorem}

\Proof Let $y \in H^1(\tilde{\Omega}_y)$ be given.  If $s \ge \eps^{\frac{\eta}{8}}$ we can apply Theorem \ref{rig-thm: V1}, so it suffices to consider $s \le \eps^{\frac{\eta}{8}}$. Recall $s_0 = \kappa \eps$ for some $\kappa  = \kappa(T,\rho, h_*, \bar{z}) \gg 1$  and assume $s \ge s_0$. The strategy is to apply Lemma \ref{rig-th: global estimate} iteratively. Set $W_0 = W^H_{-1}  = W^H_0 = \tilde{\Omega}_y \in {\cal V}^s_{k_0}$ and  $y_0 = y$. Recall $\epsilon_0 = \hat{c}^2\rho^{-1} \eps$ and define 
$$\gamma_0 := \Vert \dist(\nabla y_0, SO(2))\Vert^2_{L^2(\tilde{\Omega}_y)} \le C \frac{\rho\epsilon_0}{ \hat{c}^2}, \  \alpha_0 := \Vert \dist(\nabla y_0, SO(2))\Vert^4_{L^4(\tilde{\Omega}_y)} \le C\frac{\rho\epsilon_0}{ \hat{c}^2}.$$
 In the last inequality we used $\Vert \nabla y \Vert_\infty \le M$. Recall \eqref{rig-eq: lds}.  Set  $\hat{s}_j  = s_{j}\hat{t}^2_j$ for $j\ge 0$ and $\hat{s}_{-1} = s$, where $\hat{t}_j = C_2(t_j,h_*)$ (see \eqref{rig-eq: vartheta def}). Assume $W_{j}\in {\cal V}^{\hat{s}_{j-1}}_{k_{j}}$, $W^H_{j}  \in {\cal V}^{s_{j}}_{ k_{j}}$  are given with $ W_{j},  W^H_j \subset  \Omega^{6k_{j-1}}$,  $|W_{j} \setminus W^H_{j-1}| = 0$ and $|\tilde{\Omega}_y \setminus W_j| \le C_1 \sum^{j-1}_{i=0} k_{i}$, where we set $k_{-1} = s$.  Recall that $|W_j \setminus W^H_j|\le C_1 k_{j-1}$ and $|W^H_j \setminus H^{\lambda_{j-1}}(W_{j})|=0$, where  $\lambda_{-1} = 0$. Set 
$\beta_j = \Vert H^{\lambda_{j-1}}(W_{j})\Vert_*$ and $\beta^d_j = \Vert W_j\Vert_* - \Vert H^{\lambda_{j-1}}(W_{j})\Vert_{*}.$ Moreover, suppose there is a function  $y_{j} \in H^1(W^H_j )$ with
$$\gamma_j := \Vert \dist(\nabla y_j, SO(2))\Vert^2_{L^2(W^H_j)}, \ \ \alpha_j := \Vert \dist(\nabla y_j, SO(2))\Vert^4_{L^4(W^H_j)}$$
 such that for $j\ge 1$
\begin{align}\label{rig-eq: V3.1} 
\begin{split}
(i) &\ \ \beta_j + \beta^d_j \le (1 +  C_1 t_{j-1})\beta_{j-1} + C\epsilon^{-1}_{j-1}\gamma_{j-1} \le B_j,\\
(ii) &\ \ \gamma_{j} \le CT^{-1}_{j-1} t_{j-1}(\gamma_{j-1} + \epsilon_{j-1} \beta_{j-1}) \le  \hat{c}^{-1} t_{j-1}\rho \epsilon_{j},\\ 
(iii) &\ \ \alpha_{j}  \le  C\vartheta_{j-1}\gamma_j \le C \eps T_{j-1}, \\
(iv) & \ \ \Vert \dist (\nabla y_j,SO(2))\Vert^2_{L^\infty(W^H_j)} \le C \vartheta_{j-1},\\
(v)& \ \ \Vert \nabla y_j - \nabla y_{j-1} \Vert^4_{L^4(W_j)} \le  C\eps T_{j-1}, \\
(vi)& \ \ \Vert \nabla y_j - \nabla y_{j-1} \Vert^2_{L^2(W_j)}\le CT^{-1}_{j-1} (l^{4}_{j-1} \gamma_{j-1} + \epsilon_{j-1} \beta_{j-1})  \le  Cl^{4}_{j-1} \epsilon_{j}.
\end{split}
\end{align}
Setting $\vartheta_{-1} = 1$ and $t_{-1} = 1$, we note that, provided $\hat{c}$ is sufficiently large,   in the case $j=0$ (iii),(iv) are clearly satisfied  for $y_0 = y$ and (i),(ii) hold neglecting the second terms. We now construct $y_{j+1}$, $W_{j+1}\in {\cal V}^{\hat{s}_j}_{ k_{j+1}}$  with  $W_{j+1} \subset \Omega^{ 6k_{j}}$, $| W_{j+1} \setminus   W_j^H| = 0$  and $|\tilde{\Omega}_y  \setminus W_{j+1}| \le C_1\sum^j_{i=0}k_i$ as well as $W^H_{j+1} \in {\cal V}^{s_{j+1}}_{ k_{j+1}}$. 

First we apply Theorem \ref{rig-thm: harmonic} and let $w_j   \in H^1(W_j^H)$ be the harmonic part of $y_j$ such that similarly as in \eqref{rig-eq: w energy} 
\begin{align}\label{rig-eq: w energy3}
\Vert \nabla y_j - \nabla w_j \Vert^2_{L^2(W^H_j)} \le C\gamma_j, \ \ \  \Vert \nabla y_j - \nabla w_j \Vert^4_{L^4(W^H_j)} \le C\alpha_j
\end{align}
and so in particular $\Vert \dist(\nabla w_j, SO(2))\Vert^2_{L^2(W^H_j)} \le C\gamma_j$.  Recall $W_j^H \in {\cal V}^{s_j}_{k_j}$, $W_j \subset \Omega^{6{k_{j-1}}}$ and note $\Omega^{k_j} \subset \Omega^{6k_{j-1}}$. Then apply Lemma \ref{rig-lemma: weaklocR} with $s=s_{j}$, $k=   k_j =  s_j l_j$, $m =  t_j= t^{j+1}$,  $\epsilon = \epsilon_{j}$, $ U=W^H_j  \cap \Omega^{k_j}$, $y = w_j$ and obtain a set $\tilde{W}^H_j \in {\cal V}^{s_jt_{j}}_{ (s_j,3k_j)}$ such that 
$$\delta_4 := \sum\nolimits^4_{i=1}\Vert \nabla w_j  - \hat{R}_i\Vert^4_{L^4(\tilde{W}^H_j)} \le C\vartheta_j\gamma_j, \  \delta_2 := \sum\nolimits^4_{i=1}\Vert \nabla w_j - \hat{R}_i\Vert^2_{L^2(\tilde{W}^H_j)} \le Cl_j^4\gamma_j$$
for mappings $\hat{R}_i: (\tilde{W}^H_j)^{\circ} \to SO(2)$, $i=1,\ldots,4$, which are constant on the connected components of $Q^{k}_i(p) \cap (\tilde{W}^H_j)^\circ$, $p \in I^{k}_i(\Omega^k)$. We now use Lemma \ref{rig-th: global estimate} with $m =  t_j$, $s=s_{j}$, $\epsilon = \epsilon_{j}$, $d=d_{j}$, $W=\tilde{W}^H_j$,  $y = w_j$ and show \eqref{rig-eq: V3.1} for $j+1$. First, we obtain $W_{j+1} \in {\cal V}^{\hat{s}_{j}}_{71k_j} \subset {\cal V}^{\hat{s}_{j}}_{k_{j+1}}$, with  $W_{j+1} \subset \Omega^{6k_{j}}$, $|W_{j+1}  \setminus W^H_j| = 0$, $|(W_j^H \setminus W_{j+1})\cap \Omega^{6k_j}| \le Ck_j \Vert W_{j+1}\Vert_*$ and $W^H_{j+1}  \in {\cal V}^{s_{j+1}}_{72k_j} \subset {\cal V}^{s_{j+1}}_{k_{j+1}}$  with $|W^H_{j+1} \setminus  H^{\lambda_{j}}( W_{j+1})| = 0$ and $|W_{j+1} \setminus W^H_{j+1}| \le C_1 k_j$. Recall $\Vert W^H_j \Vert_* \le (1+C_1 t_j) \beta_j$ by \eqref{rig-eq: mayblast 10}. Thus, we have
\begin{align}\label{rig-eq: V3.9}
\Vert W_{j+1} \Vert_* \le (1 +  C_1 t_j)\Vert W^H_j \Vert_* + C\epsilon_j^{-1} (\gamma_j + \vartheta_j\gamma_j) \le (1 +  C_1 t_j)\beta_j + C\epsilon_j^{-1} \gamma_j
\end{align}
by \eqref{rig-eq: weaklocR2}, \eqref{rig-eq: glob1} and the fact that  $\vartheta_j \le 1$ (see \eqref{rig-eq: vartheta char}).  Moreover, we get a function $y_{j+1} \in H^1(W^H_{j+1})$ with  (see \eqref{rig-eq: 5 prop2}, \eqref{rig-eq: glob4})
\begin{align}\label{rig-eq: V3.13}
\begin{split}
(i) & \ \ \Vert\dist(\nabla y_{j+1},SO(2) )\Vert^2_{L^2(W_{j+1}^H)} \le CC^2_{t_j}(\gamma_j + \epsilon_j \beta_j), \\
(ii) & \ \ \Vert \nabla w_j-  \nabla y_{j+1}\Vert^2_{L^2(W_{j+1})} \le  CC^2_{t_j}(\gamma_j + l_j^{4} \gamma_j  + \epsilon_j \beta_j),\\
(iii) &\ \ \Vert \nabla w_j-  \nabla y_{j+1}\Vert^4_{L^4(W_{j+1})} \le  CC^2_{t_j}\vartheta_j (\gamma_j + \epsilon_j \beta_j),\\
(iv) & \ \ \Vert \dist (\nabla y_{j+1},SO(2))\Vert^2_{L^\infty(W^H_{j+1})} \le C\vartheta_j,
\end{split}
\end{align}
where we again used that $\vartheta_j \le 1$. The first inequality in \eqref{rig-eq: V3.1}(ii) follows directly noting that $T^{-1}_{j}t_j \ge C^2_{t_j}$ and for the second inequality we use \eqref{rig-eq: V3.1}(i),(ii) for iteration step $j$ as well as \eqref{rig-eq: B def} to see 
\begin{align}\label{rig-eq: Tgammabeta}
CT^{-1}_j  (\gamma_j + \epsilon_j \beta_j) \le CT^{-1}_j  \rho\epsilon_j (1 + \rho^{-1} B_j) \le \rho  \hat{c}^{-1} PT^{-1}_j  \epsilon_j =  \hat{c}^{-1}  \rho \epsilon_{j+1},
\end{align}
where we choose  $\hat{c}$ sufficiently large. Likewise, \eqref{rig-eq: V3.1}(i) follows by \eqref{rig-eq: V3.9}, the fact that $\Vert W_{j+1} \Vert_* = \beta_{j+1} + \beta^d_{j+1}$  and 
\begin{align*}
\beta_{j+1} + \beta^d_{j+1} &\le (1 +  C_1 t_j)B_j +   \rho t_{j-1} \\
&\le  \big( \Vert \tilde{\Omega}_y \Vert_* + C_* \rho\big) \cdot \text{\scriptsize $\sum\nolimits$}^{j-1}_{i=0}  t^{i} \cdot \pi^{j}_{t=0} (1+ C_* t^{i+1}) 
+   \rho t^{j}\\
&\le  \big( \Vert \tilde{\Omega}_y \Vert_* + C_* \rho\big) \cdot \text{\scriptsize $\sum\nolimits$}^{j}_{i=0}  t^{i} \cdot \pi^{j}_{t=0} (1+ C_* t^{i+1}) =B_{j+1} .
\end{align*}
Here we have again chosen $\hat{c}$ and $C_*$ large enough  (with respect to $C$ and $C_1$, respectively). This also implies $|(W_j \setminus W_{j+1}) \cap \Omega^{6k_j}| \le C k_{j}$ by \eqref{rig-eq: V3.1}(i) and thus $|(\tilde{\Omega}_y \setminus W_{j+1})| \le C \sum^j_{i=0}k_{i} + |\Omega \setminus \Omega^{6k_j}| \le C\sum^j_{i=0}k_{i}$. 

Estimate \eqref{rig-eq: V3.1}(iv) follows from \eqref{rig-eq: V3.13}(iv). The first inequality in \eqref{rig-eq: V3.1}(iii) is a consequence of \eqref{rig-eq: V3.1}(iv), the second inequality is implied by the fact that $\eps =  \hat{c}^{-2} \rho \epsilon_0$, \eqref{rig-eq: V3.1}(ii) and \eqref{rig-eq: vartheta char}. Moreover, \eqref{rig-eq: V3.1}(v) follows from  \eqref{rig-eq: V3.1}(iii),  \eqref{rig-eq: w energy3}, \eqref{rig-eq: V3.13}(iii) and the fact that $\vartheta_j C^2_{t_j}(\gamma_j +\epsilon_j \beta_j)  \le \vartheta_j  \rho \epsilon_{j+1} \le C\eps T_j$  by \eqref{rig-eq: vartheta char} and \eqref{rig-eq: Tgammabeta}. Similarly, \eqref{rig-eq: V3.1}(vi) follows from \eqref{rig-eq: V3.13}(ii), \eqref{rig-eq: w energy3} and \eqref{rig-eq: Tgammabeta}.

We now choose $j^* \in \N$ such that 
\begin{align}\label{rig-eq: final en}
\eps^{3\omega} \ge s_{j^*} \ge \eps^{4\omega}, \ \ \ \ \ \ \ \epsilon_{j^*} \le C\eps^{1 - \omega} T^{2}_{j^*}
\end{align}
holds  for $\eps$ sufficiently small. The first inequality is possible by \eqref{rig-eq: l def} and we obtain $j^* \le J^*= \lceil\text{log}_{1+r}(\text{log}_T \eps^{\omega}))  + \frac{1}{\omega}\rceil$.  Indeed, by  \eqref{rig-eq: V3.4}  and the fact that $\bar{z} \ge 1$ we get $s_j \ge \eps^{-\frac{\omega}{r}} \epsilon_j = \eps^{-\frac{\eta}{2}} \epsilon_j$ for $j > \lceil\text{log}_{1+r}(\text{log}_T \eps^{\omega}))\rceil$ and therefore $\hat{J} \le \lceil\text{log}_{1+r}(\text{log}_T \eps^{\omega}))\rceil$. The second inequality can be derived arguing as in \eqref{rig-eq: V3.4_*}. Similarly, proceeding as in \eqref{rig-eq: V3.4_*} we have $t^{-2}_{j_*}  = o(\eps^{-\omega})$ for $\eps \to 0$ and thus $k_{j^*} = s_{j^*} d_{j^*} t^{-2}_{j^{*}} = o(\eps^{\omega})$. This implies $\Omega^{6k_{j^*}} \supset \Omega^{\varrho}$ for $\eps$ small enough. We let 
$$y_* = y_{j^*}, \ \  \ W^H_* = W^H_{j^*} \cap \Omega^{\varrho}, \ \ \   W_* = \bigcap\nolimits^{j^*}_{i=0} W_i \cap \Omega^{\varrho}.$$
It is not hard to see that $|\tilde{\Omega}_y \setminus W_*| \le C_1 \sum^{j^*}_{i=0} k_{i}  \le C \varrho$.   As $\hat{s}_j = s_{j} \hat{t}^2_j$ is increasing in $j$  (note that $d_j \ge \hat{t}^{-2}_j$ for all $j$, see e.g. \eqref{rig-eq: V3.4}), we find $W_* \in {\cal V}^{\hat{s}_0}$. 

 The strategy is now to establish an estimate of the form \eqref{rig-eq: V1.11.2} and \eqref{rig-eq: V1.11}. Observe that $s_{j^*} \ge \eps^{\frac{\eta}{8}}$, i.e. for the function $y_* \in H^1(W_*^H)$ we may proceed as in Theorem \ref{rig-thm: V1}  (replacing $s$ by $s_{j^*}$). Similarly as in \eqref{rig-eq: w energy}, we apply Theorem \ref{rig-thm: harmonic} and let $w_*$ be the harmonic part of $y_*$ with 
\begin{align}\label{rig-eq: harmonic estimate}
\Vert \nabla w_* - \nabla y_*\Vert^2_{L^2(W^H_*)} \le C \eps^{1-\frac{\eta}{2}}, \ \ \ \Vert \nabla w_* - \nabla y_*\Vert^4_{L^4(W^H_*)} \le C\eps T^{j^*}.
\end{align}
by \eqref{rig-eq: V3.1}, \eqref{rig-eq: final en} and $\omega \le \frac{\eta}{2}$. Apply Lemma \ref{rig-lemma: weaklocR} on $ W^H_* \subset \Omega^\varrho$ for the function $w_*$ and $k=\rho^{q-1} = \varrho \rho^{-1}$,  $s=\eps^{4\omega}$,  $\epsilon= \hat{c}\rho^{-1}\eps^{1-\frac{\eta}{2} }$, $m=\rho$.  (Without restriction we can assume $s^{-1} \in \N$.) We find a set $W^H \subset \Omega^{3k}$, $W^H \in {\cal V}^{s_{j^*} m}_{3k}$ such that  
\begin{align}\label{rig-eq: V1.41_*}
\Vert W^H \Vert_* \le (1+ C_1\rho)\Vert W_*^H \Vert_* + C \hat{c}^{-1}\rho\eps^{\frac{\eta}{2} - 1} \eps^{1-\frac{\eta}{2}} \le \Vert W_*^H \Vert_* + C_1\rho
\end{align}
by \eqref{rig-eq: weaklocR2} as well as $|W^H_* \setminus W^H| \le |(W^H_* \setminus W^H) \cap \Omega^{3k}| + C_1k \le C_1k \le C_1\rho$. Moreover, there are mappings $\hat{R}_i: (W^H)^{\circ} \to SO(2)$, $i=1,\ldots,4$, which are constant on the connected components of $Q^{k}_i(p) \cap (W^H)^{\circ}$, $p \in I^{k}_i(\Omega)$, such that by \eqref{rig-eq: weaklocR1}(i) and \eqref{rig-eq: harmonic estimate}
\begin{align*}
\Vert\nabla y_* - \hat{R}_i\Vert^4_{L^4(W^H)} \le C\Vert\nabla w_* - \hat{R}_i\Vert^4_{L^4(W^H)} + C\eps T^{j^*} \le    C\vartheta \eps^{1-\frac{\eta}{2}} + C\eps\le C\eps,
\end{align*}
where similarly as before equation \eqref{rig-eq: V1.11} we compute (recall \eqref{rig-eq: final en} and $\omega = \frac{\eta}{36}$) $\vartheta \le C(\rho,q)  s^{-10} \epsilon \le C(\rho,q) \eps^{-40\omega}\eps^{1-\omega} = C(\rho,q)\eps^{1 - \frac{41}{36} \eta} \le \eps^{\frac{\eta}{2}}$ {for $\eps,\eta$ small enough}. Likewise, we derive
$$\Vert\nabla y_* - \hat{R}_{i}\Vert^2_{L^2(W^H)} \le C\Vert\nabla w_* - \hat{R}_i\Vert^2_{L^2(W^H)} + C\eps^{1-\frac{\eta}{2}} \le C(1+l^4) \eps^{1-\frac{\eta}{2}} \le C\eps^{1-\eta}$$
as $l = \frac{k}{s} \le C\eps^{-4\omega} \le \eps^{-\frac{\eta}{8}}$.

 We now will construct a set $W \in {\cal V}^{\hat{s}_0}_{143k}$ which is contained in $W^H \cap W_* \cap \Omega^{3k} \in {\cal V}^{\hat{s}_0}$, where the two sets coincide up to a set of measure smaller than $C_1 \rho$. (Similarly as before the difference of the sets is related to the definition of the boundary components.) Before we give the exact definition of $W$ and establish an estimate of the form \eqref{rig-eq: V1.41}, we first observe $| \tilde{\Omega}_y \setminus W| \le C_1\rho$ arguing as before  and derive estimates similar to \eqref{rig-eq: V1.11.2} and \eqref{rig-eq: V1.11}.

We iteratively apply \eqref{rig-eq: V3.1}(v) and derive for $i=1,\ldots,4$
\begin{align}\label{rig-eq: R4}
\Vert\nabla y - \hat{R}_{ i}\Vert^4_{L^4(W)} \le  C \Big(\sum\nolimits^{j^*}_{\iota=1} (\eps T_{\iota-1})^{\frac{1}{4}}\Big)^4   + C\Vert\nabla y_* - \hat{R}_i\Vert^4_{L^4(W)}  \le C \eps.
\end{align}
Likewise, observe that by \eqref{rig-eq: l def}, \eqref{rig-eq: lds} and \eqref{rig-eq: final en} we have  $l_{j-1}^{4} \epsilon_j \le l_{j}^{4} \epsilon_j= d_j^{4} t^{-8(j+1)} \epsilon_j\le \eps^{-4\omega}\eps^{1-\omega} T_j \le \eps^{1-\eta} T_j$. We derive  by \eqref{rig-eq: V3.1}(vi)
$$\Vert\nabla y - \hat{R}_{i}\Vert^2_{L^2(W)} \le C\eps^{1-\eta} \Big(\sum\nolimits^{j^*}_{\iota=1} T^{\frac{1}{2}}_{\iota}\Big)^2  + C\Vert\nabla y_* - \hat{R}_i\Vert^2_{L^2(W)}  \le C\eps^{1-\eta}$$
for $i=1,\ldots,4$.

It remains to give the exact definition of $W \in {\cal V}^{\hat{s}_0}_{143k}$ and to establish $\Vert  W\Vert_* \le (1 + Ch_*){\cal H}^1(J_y) + C\rho$. Recall  $W_0 = \tilde{\Omega}_y$ and define $W_{j^*+1} := W^H$ for notational convenience. We now define $W$ inductively.

Let $Y_0 = Y_0' = Y_0'' = W_0$. Assume $Y_j \in {\cal V}^{\hat{s}_0}$ and $Y_j' \in {\cal V}^{\hat{s}_0}_{k_j}$, $Y_j'' \in {\cal V}^{\hat{s}_0}$ are given with $|Y'_j \setminus Y_j| + |Y'_j \triangle Y_j''| = 0$, $|Y_j \setminus Y'_j| \le C_1 k_{j-1}$ and
$$\max \lbrace\Vert Y'_{j} \Vert_*, \Vert Y''_{j} \Vert_*\rbrace \le \Vert Y_j \Vert_* \le \Vert W_j \Vert_* + \sum\nolimits^{j-1}_{i=1} \beta^d_i,$$
where $Y_j''$ has the property that all components not intersecting $\partial H^{\lambda_{j-1}}(W_j)$ coincide with components of $Y_j'$ and the set $\big(X_t(H^{\lambda_{j-1}}(W_j))\big)_t$ of components of $H^{\lambda_{j-1}}(W_j)$ is a subset of the components of $Y_j''$. Moreover, suppose that $|Y'_j \setminus \bigcap^j_{i=0} W_i| = 0$ and $| \bigcap^j_{i=0} W_i\setminus Y'_j|\le \sum^{j-1}_{i=0} k_i$.

We now pass to step $j+1$. Let $X_1(W_{j+1}), \ldots,X_{n_{j+1}}(W_{j+1})$ be the components of $W_{j+1}$ and define 
$$Y_{j+1} = \big(Y_j'' \setminus \bigcup\nolimits^{n_{j+1}}_{t=1} X_t(W_{j+1})\big) \cup \bigcup\nolimits^{n_{j+1}}_{t=1} \partial X_t(W_{j+1}) \in {\cal V}^{\hat{s}_0}.$$
 First observe that $Y_{j+1}$ satisfies $|Y_{j+1} \setminus \bigcap^{j+1}_{i=0} W_i| = 0$ and $| \bigcap^{j+1}_{i=0} W_i\setminus Y_{j+1}|\le \sum^{j-1}_{i=0} k_i$. As $|W_{j+1} \setminus W_j^H| = 0$, we obtain  $\bigcup^{n_{j+1}}_{t=1} \overline{X_t(W_{j+1})} \supset \bigcup_{t} \overline{X_t(W^H_{j})}$ and then by the fact that $|W_j^H \setminus H^{\lambda_{j-1}}(W_j)| = 0$ we get $\bigcup^{n_{j+1}}_{t=1} \overline{X_t(W_{j+1})}  \supset \bigcup_{t} \overline{X_t(H^{\lambda_{j-1}}(W_j))}$. As by hypothesis  the components of $H^{\lambda_{j-1}}(W_j)$ are also components of $Y_j''$, we  derive recalling $\beta^d_i = \Vert W_{i} \Vert_* - \Vert H^{\lambda_{j-1}}(W_j)\Vert_* $ and $\beta^d_0 = 0$
\begin{align*}
\Vert Y_{j+1} \Vert_* &\le \Vert Y''_j \Vert_* +\Vert W_{j+1}\Vert_* - \Vert H^{\lambda_{j-1}}(W_j) \Vert_*  = \Vert W_{j+1}\Vert_* + \sum\nolimits^{j}_{i=1} \beta^d_i.
\end{align*}
 Observe that possibly $Y_{j+1} \notin {\cal V}^{\hat{s}_0}_{\rm con}$. However, by Lemma \ref{rig-lemma: modifica}(ii) we find a set $Y_{j+1}' \in {\cal V}^{\hat{s}_0}$ with $|Y_{j+1} \setminus Y'_{j+1}| \le C_1 k_j$ and $\Vert Y'_{j+1} \Vert_* \le \Vert Y_{j+1} \Vert_*$.  Here we essentially used the rectangular shape of the boundary components given by \eqref{rig-eq: mayblast5} and \eqref{rig-eq: new111}, respectively. Then it is elementary to see that $Y_{j+1}' \in {\cal V}^{\hat{s}_0}_{143k_j} \subset {\cal V}^{\hat{s}_0}_{k_{j+1}}$ and $| \bigcap^{j+1}_{i=0} W_i\setminus Y'_{j+1}|\le \sum^{j}_{i=0} k_i$. Moreover, if $j+1 \le j^*$,  we let $Y_{j+1}'' = (Y_{j+1}'\cap H^{\lambda_j}(W_{j+1})) \cup \partial H^{\lambda_j}(W_{j+1})$ and observe that $Y_{j+1}''$ has the desired properties. In fact, $\Vert Y_{j+1}'' \Vert_* \le \Vert Y_{j+1}\Vert_*$ follows as before. Components not intersecting $\partial H^{\lambda_j}(W_{j+1})$ are clearly components of $Y_{j+1}'$. Finally, by definition components of $H^{\lambda_j}(W_{j+1})$ are also components of $Y_{j+1}''$.

We finally define $W = Y'_{j^*+1} \cap \Omega^{3k} \in {\cal V}^{\hat{s}_0}_{143k}$. By \eqref{rig-eq: B def}  and  \eqref{rig-eq: V3.1}(i),(ii) we have
$$\beta^d_i \le \beta_{i-1} - \beta_i + C_1t^{i}\beta_{i-1}+ C\epsilon^{-1}_{i-1}\gamma_{i-1} \le \beta_{i-1} - \beta_i + C_1t^{i}B+ \rho t^{i - 1}$$
for $i=1,\ldots, j^*$. Recalling  $\beta_ 0= \Vert \tilde{\Omega}_y \Vert_*$, $  \Vert W^H_* \Vert_* \le (1+C_1 t_{j^*})\beta_{j^*} $ and using \eqref{rig-eq: B def},  \eqref{rig-eq: V1.41_*}  as well as $t \le \rho$  we conclude 
\begin{align*}
\begin{split}
\Vert W\Vert_* & \le \Vert Y'_{j^*+1}\Vert_* \le \Vert W^H\Vert_*    + \sum\nolimits^{j^*}_{i=1} (\beta_{i-1} - \beta_i + C_1t^{i}B+ \rho t^{i-1}) \\
& \le \Vert W^H\Vert_* - \beta_{j^*} + \beta_0  + C_1 \rho B  +  C_1\rho \le   C_1\rho  + \Vert \tilde{\Omega}_y \Vert_*  +  C_1 \rho B  +  C_1\rho\\
& \le (1+ C_1\rho)\Vert \tilde{\Omega}_y \Vert_* +  C_1\rho \le  (1 + C_1h_*){\cal H}^1(J_y) + C_1\rho,
\end{split}
\end{align*}
as derided.

We now proceed as in the proof of Theorem \ref{rig-thm: V1} after equation \eqref{rig-eq: V1.11} with the only difference that we take $\hat{s}_{0}$ instead of $s \sim \eps^{\frac{\eta}{8}}$ in the application of Corollary \ref{rig-cor: weakA}. However, this does not change the analysis. This leads to a set $\Omega_y \in {\cal V}^{\hat{s}_0 \hat{m}}_{ck}$ with $\Omega_y \subset \Omega^{5k}$  and  $|\Omega \setminus \Omega_y |\le C_1\rho$ for $k= \rho^{q-1}, m = 3\rho$ for which \eqref{rig-eq: main properties} can be established. \eop

We now additionally treat the subatomistic regime by dropping the assumption $s \ge \kappa \eps$.

\begin{theorem}\label{rig-thm: V2.5}
Theorem \ref{rig-th: rigidity2} holds under the additional assumption that there is an $\tilde{\Omega}_y \subset \Omega^{s}$, $\tilde{\Omega}_y \in {\cal V}^{s}_{\eps}$ for some  $0 < s \ll \eps$ such that $y \in H^1(\tilde{\Omega}_y)$, $\Vert \tilde{\Omega}_y\Vert_* \le (1 + C_1h_*){\cal H}^1(J_y) + C_1\rho$ and $|\Omega \setminus \tilde{\Omega}_y| \le C_1\rho$ for a constant $C_1 = C_1(\Omega, M,\eta)$.
\end{theorem}

\Proof
Let again $\rho^{-1}\in \N$,  $s_0 = \kappa \eps$  and recall $\Vert \dist(\nabla y, SO(2))\Vert^2_{L^2(\Omega)} \le C\eps$.  As $\kappa \gg 1$ was chosen in dependence of $T$ and $T=T(\rho, h_*)$ (see \eqref{rig-eq: V3.3}), we can suppose $\kappa = \kappa(\rho, h_*)$. Applying Lemma \ref{rig-lemma: local estimate2} for $s$, $k=\rho^{-2}  \kappa \eps$, $m=  \rho $ and $\epsilon=  \rho^{-2} \kappa\eps$,  $U = \tilde{\Omega}_y \cap \Omega^k$ there is a set $W \subset \Omega^{3k}$ with $W \in {\cal V}^s_{k}$, $|\tilde{\Omega}_y \setminus W| \le C_1 k\le C_1\rho$ for $\eps$ small enough and 
$$\Vert W\Vert_* \le \Vert  \tilde{\Omega}_y \Vert_* + C\epsilon^{-1} \eps \le \Vert \tilde{\Omega}_y \Vert_* + \rho.$$
The last inequality holds by choosing $\kappa$ larger than $C$. Moreover, there are mappings $\hat{R}_i: \Omega^{3k} \to SO(2)$, $i=1,\ldots,4$, which are constant on $Q^{k}_i(q) \cap W$, $q \in I^{k}_i(\Omega^{k})$, such that
\begin{align*}
\Vert \nabla y  - \hat{R}_i\Vert^2_{L^2(W)} \le C \eps + C\eps \rho^{-2}\kappa \Vert \tilde{\Omega}_y \Vert_* \le C  \rho^{-2}\kappa \eps.
\end{align*}
Clearly, we also get $\Vert \nabla y  - \hat{R}_i\Vert^4_{L^4(W)}  \le C \rho^{-2}\kappa \eps$ as $\Vert \nabla y\Vert_\infty \le M$. Then we apply Lemma \ref{rig-th: global estimate} for $k=\rho^{-2}s_0$, $\nu=s_0$, $m= \rho$  and $\epsilon=\hat{c}\rho^{-3}\kappa\eps$ to get sets $U \in {\cal V}^{s\hat{m}^2}_{71k}$ and $U^H  \in {\cal V}^{\nu}_{ 72k}$ with $U,  U^H \subset  \Omega^{6k}$,  $|U \setminus W| =0$, $|U^H \setminus H^{\frac{\nu}{m}}(U)|=0$ and 
\begin{align*}
\Vert  U \Vert_* \le (1+ C_1\rho)\Vert W\Vert_* + C\epsilon^{-1}  \rho^{-2}\kappa\eps \le \Vert \tilde{\Omega}_y \Vert_* + C_1\rho
\end{align*}
as well as $|W \setminus U| \le C_1 k \le C_1\rho$ for $\eps$ small enough. Moreover, we find a function $\hat{y} \in H^{1}(U^H)$  such that by \eqref{rig-eq: 5 prop2} 
\begin{align*}
(i) & \ \ \Vert\dist(\nabla\hat{y},SO(2) )\Vert^2_{L^2(U^H)} \le CC^{2}_\rho (\rho^{-2}\kappa \eps +  \rho^{-3}\kappa \eps \Vert W \Vert_*) \le C  C^{2}_\rho \rho^{-3} \kappa  \eps, \\
(ii) & \ \ \Vert \dist (\nabla \hat{y},SO(2))\Vert^2_{L^\infty(U^H)} \le C C^6_\rho,\\
(iii) & \ \ \Vert \nabla y-  \nabla \hat{y}\Vert^2_{L^2(U')} \le C  C^2_\rho \rho^{-3} \kappa  \eps, \ \ \ \Vert \nabla y-  \nabla \hat{y}\Vert^4_{L^4(U')} \le C  C^{8}_\rho \rho^{-3} \kappa  \eps,
\end{align*}
 where the second part of (iii) follows from (ii). Note that this also implies $\Vert\dist(\nabla\hat{y},SO(2) )\Vert^4_{L^4(U^H)} \le C  C^{8}_\rho \rho^{-3} \kappa  \eps$.  Setting $W_1 = U$, $W_1^H = U^H$, $y_1 = \hat{y}$ we can now follow the proof of Theorem \ref{rig-thm: V2} beginning with \eqref{rig-eq: V3.1}  with the essential difference that we have to replace $\eps$ by $  C C^{8}_\rho \rho^{-3} \kappa  \eps$. We then obtain the desired result for a possibly larger constant $C_2$ in \eqref{rig-eq: main properties}. \eop

\subsection{Step 3: General case}\label{rig-sec: subsub,  step3}

We are now in a position to prove the general version of Theorem \ref{rig-th: rigidity2}.

\noindent {\em Proof of Theorem \ref{rig-th: rigidity2}.}  Let $y \in SBV_{M}(\Omega) \cap L^2(\Omega)$  be given and let $\rho>0$. It suffices to find a set $\tilde{\Omega} \in {\cal V}^s_{\eps}$, $s>0$, and a function $\tilde{y} \in H^1(\tilde{\Omega})$ with  $\Vert \tilde{y} \Vert_{L^\infty(\tilde{\Omega})}  + \Vert \nabla \tilde{y} \Vert_{L^\infty(\tilde{\Omega})}\le cM$ for a universal  constant $c>0$ such that 
\begin{align}\label{rig-eq: V3.20}
\begin{split}
&|\Omega \setminus  \tilde{\Omega}| \le C_1\rho, \ \  \Vert\tilde{\Omega}\Vert_*  \le (1 + C_1h_*){\cal H}^1(J_y) + C_1\rho, \\ &\Vert y - \tilde{y}\Vert^2_{L^2(\tilde{\Omega})}  + \Vert \nabla y - \nabla \tilde{y}\Vert^2_{L^2(\tilde{\Omega})}\le C_1\eps\rho.
\end{split}
\end{align} 
Then the result follows from Theorem \ref{rig-thm: V2.5}  applied on the function $\tilde{y}$. (Accordingly, replace $M$ by $cM$ in all estimates.) Note that we cannot just apply density results for SBV functions (see \cite{Cortesani:1997})  since in general such approximations do not preserve an $L^\infty$ bound for the derivative. The problem may be bypassed  by  construction of a different approximation (see \cite{Chambolle:2004} and \cite{Friedrich:15-1}) at the cost of a non exact approximation of the jump set which, however, suffices for our purposes.  

Let $\mu = \eps \rho $. Recall that $J_{y}$ is rectifiable (see \cite[Section 2.9]{Ambrosio-Fusco-Pallara:2000} ), i.e. there is a countable union of $C^1$ curves $(\Gamma_i)_{i \in \N}$ such that ${\cal H}^1(J_{y} \setminus \bigcup_i \Gamma_i) = 0$. Covering $J_{y}$ with small balls and applying Besicovitch's covering theorem (see \cite[Corollary 1, p. 35]{EvansGariepy92}) we find finitely many closed, pairwise disjoint balls $\overline{B_j} = \overline{B_{r_j}(x_j)}$, $j=1,\ldots,n$ with $r_j \le \mu$ such that ${\cal H}^1(J_{y}\setminus \bigcup^n_{j=1} B_j) \le \mu$. Moreover, we get ${\cal H}^1(J_{y} \cap \overline{B_j}) \ge 2(1-\mu)r_j$ and for each $B_j$ we find a $C^1$ curve $\Gamma_{i_j}$ such that  $\Gamma_{i_j} \cap \overline{B_j}$ is connected and ${\cal H}^1((\Gamma_{i_j} \triangle J_{y})  \cap \overline{B_j}) \le 2 \mu r_j \le \frac{\mu}{1 - \mu}  {\cal H}^1(J_{y} \cap \overline{B_j})$. For a detailed proof we refer to \cite[Theorem 2]{Chambolle:2004}.

We choose rectangles $R_j$ with $|\partial R_j|_\infty \le 2\sqrt{2}r_j$ such that  ${\cal H}^1(\Gamma_{i_j} \cap (B_j \setminus R_j)) = 0$ and $|\partial R_j|_\infty \le {\cal H}^1(\Gamma_{i_j} \cap \overline{B_j})$. We then obtain
\begin{align*}
\sum\nolimits_j |\partial R_j|_\infty &  \le  \sum\nolimits_j {\cal H}^1(\Gamma_{i_j} \cap \overline {B_j}) \\
& \le \Big(1 + \frac{\mu}{1-\mu}\Big) \sum\nolimits_j {\cal H}^1(J_y \cap \overline {B_j}) \le (1+C_1\mu){\cal H}^1(J_{y})
\end{align*}
and likewise $\sum_j |\partial R_j|_{\cal H} \le C_1{\cal H}^1(J_{y})$. Choose rectangles $Q_j$ with $R_j \subset \subset Q_j$ such that $|\partial Q_j|_* \le (1+\mu)|\partial R_j|_*$ and
\begin{align}\label{rig-eq: no boundary}
{\cal H}^1\Big(\bigcup\nolimits_j \partial Q_j \cap J_y \Big) = 0.
\end{align}
As before it is not hard to see that $R_{j_1} \setminus R_{j_2}$ is connected for $1 \le j_1, j_2 \le n$. The rectangles $(Q_j)_j$ can be chosen in a way such that they fulfill the same property. Possibly replacing the rectangles by infinitesimally larger rectangles we can assume that there is some $s>0$ such that $R_j,Q_j \in {\cal U}^s$ for $j=1,\ldots,n$. Then by Lemma \ref{rig-lemma: modifica}(i) we find sets $W,V \in {\cal V}^s_{\eps}$ with $|V \triangle  (\Omega^\rho \setminus \bigcup_j R_j)| = 0$ and  $|W \triangle  (\Omega^\rho \setminus \bigcup_j Q_j)| = 0$. Note that $ W^\circ  \subset \subset V^\circ$ and $|\Omega \setminus W| \le C_1\rho$. It is not restrictive to assume that  corners of $R_j,Q_j$ do not coincide and thus $ W^\circ,V^\circ$ are Lipschitz domains. We get (recall Lemma \ref{rig-lemma: infty})
\begin{align}\label{rig-eq: v4.3}
\Vert W\Vert_* \le (1+\mu)\sum\nolimits_j |\partial R_j|_*  \le (1+ C_1\rho + C_1h_*){\cal H}^1(J_{y}).
\end{align}
Moreover,  as ${\cal H}^1(J_{y}\setminus \bigcup^n_{j=1} B_j) \le \mu$ we get
\begin{align}\label{rig-eq: V4.1}
{\cal H}^1(J_{y}\setminus \bigcup\nolimits^n_{j=1} R_j)&\le \mu + {\cal H}^1\Big(\bigcup\nolimits^n_{j=1} J_{y} \cap  (B_j \setminus R_j)\Big) \notag\\ 
& \le \mu + {\cal H}^1\Big(\bigcup\nolimits^n_{j=1} \Gamma_{i_j} \cap (B_j \setminus R_j)\Big) + {\cal H}^1\Big(\bigcup\nolimits^n_{j=1} (\Gamma_{i_j} \triangle J_y) \cap \overline{B_j}\Big)\notag \\
& \le \mu + \frac{\mu}{1-\mu} {\cal H}^1(J_y) \le C_1\mu,
\end{align} 
where in the last step we have used ${\cal H}^1(\Gamma_{i_j} \cap (B_j \setminus R_j)) = 0$. We now show that there is a function $\hat{y} \in SBV( W^\circ)$ with  $\Vert y - \hat{y}\Vert^2_{L^2(W)}  + \Vert \nabla y - \nabla \hat{y}\Vert^2_{L^2(W)}\le C_1\eps\rho$ such that $\Vert \nabla \hat{y} \Vert_\infty \le cM$ and $J_{\hat{y}}$ is a finite union of closed segments satisfying ${\cal H}^1(J_{\hat{y}}) \le C_1\mu \le C_1\rho$. We apply a result by Chambolle obtained in \cite{Chambolle:2004} in an SBD-setting and rather cite the result as repeating the arguments. Therefore, we first obtain a control only  over the symmetric part of the gradient. To derive the desired result we repeat the arguments for the function $v= (y^2,y^1)$ instead of $y = (y^1,y^2)$ to control also the skew part.

We define
$$E(y,W^\circ) = \int_{W^\circ} V(e(\nabla y)) + {\cal H}^1(J_y \cap W^\circ)$$
and $E_{c}(y,W^\circ) = E(y,W^\circ) + c{\cal H}^1(J_y \cap W^\circ)$, where $V(A) := \frac{1}{2\pi} \int_{S^1} (\xi^T A\xi)^2 \, d\xi$ for $A \in \R^{2 \times 2}$.  As $y \in SBV_{M}( W^\circ) \cap L^2(W^\circ)$ with $E(y, W^\circ) < + \infty$ and $ W^\circ$ has Lipschitz boundary, by \cite[Theorem 1]{Chambolle:2004} we find a sequence $y_n \in SBD(W^\circ) \cap L^2(W^\circ)$ with $\Vert y_n - y \Vert_{L^2(W^\circ)} \to 0$ such that $\overline{J_{y_n}}$ is  a finite union of closed segments and
\begin{align}\label{rig-eq: V4.2}
\begin{split}
\limsup_{n \to \infty} E(y_n,W^\circ) &\le E_c(y, W^\circ) \le E(y,W^\circ) + C_1\mu \\ & \le \int_{W^\circ} V(e(\nabla y)) + C_1\mu.
\end{split}
\end{align}
In the second  and third step we used \eqref{rig-eq: V4.1}.  The proof is based on a discretization argument. Consequently, as a preparation an extension $y'$  to some set $W' \supset \supset W^\circ$ with $E(y',W') \le E(y,W^\circ) + \delta$ for arbitrary $\delta>0$ had to be constructed (see \cite[Lemma 3.2]{Chambolle:2004}). In our framework we can choose $y' = y$ due to $W^\circ \subset \subset V^\circ$ and \eqref{rig-eq: no boundary}. 
Moreover, $\Vert y_n \Vert_{\infty}  \le \Vert  y \Vert_{\infty}$ holds. Although not stated explicitly in the theorem, the approximations satisfy $\Vert \nabla y_n \Vert_{L^\infty(W^\circ)} \le c\Vert \nabla y' \Vert_{L^\infty(W')} \le c\Vert \nabla y \Vert_{L^\infty(V)} \le  cM$. (For a precise argument see the proof of \cite[Theorem 3.1]{Chambolle-Giacomini-Ponsiglione:2007}, where a similar construction is used.) Strictly speaking, the theorem only states that $J_{y_n}$ is essentially closed and contained in a finite union of closed segments. However, the proof shows that up to an infinitesimal perturbation of $y_n$ (do not set $y_n=0$ on a `jump square', but $y_n= \tilde{c}$ for $\tilde{c} \approx 0$) the desired property can be achieved. 

 By \cite[Lemma 5.1]{Chambolle:2004} we obtain  weak convergence $e(\nabla y_n) \rightharpoonup e(\nabla y)$ in $L^2(W^\circ)$ up to a not relabeled subsequence. Together with the lower semicontinuity results $\int_{W^\circ} V(e(\nabla y)) \le \liminf_{n \to \infty} \int_{W^\circ}  V(e(\nabla y_n))$ and  ${\cal H}^1(J_y) \le \liminf_{n \to \infty} {\cal H}^1(J_{y_n})$ (see \cite[Lemma 5.1]{Chambolle:2004}) we find  by \eqref{rig-eq: V4.2}
$$\int_{W^\circ}  V(e(\nabla y)) \le \limsup_{n \to \infty} \int_{W^\circ}  V(e(\nabla y_n)) \le \int_{W^\circ}  V(e(\nabla y)) + C_1\mu.$$
Consequently, by weak convergence we obtain 
\begin{align*}
\limsup_{n \to \infty} \Vert e(\nabla y_n) - e(\nabla y)\Vert^2_{L^2({W^\circ})} &\le c\limsup_{n \to \infty} \int_{W^\circ} V(e(\nabla y_n - \nabla y))\\
& \le c \limsup_{n \to \infty} \Big(\int_{W^\circ} V(e(\nabla y_n)) - \int_{W^\circ} V(e(\nabla y))\Big)  \\ & \le C_1\mu = C_1\eps\rho.
\end{align*}
Then by \eqref{rig-eq: V4.2} we also get $\limsup_{n \to \infty} {\cal H}^1(J_{y_n}) \le C_1\mu \le C_1\rho$. 
We now repeat the argument for $v = (y^2,y^1)$ instead of $y$ and observe that by construction the approximations can be chosen as $v_n = (y_n^2,y_n^1)$. We find that $y_n \in SBV(W^\circ)$ and $\limsup_{n \to \infty} \Vert \nabla y_n - \nabla y\Vert^2_{L^2(W^\circ)} \le C_1\eps\rho$. Now choose $n$ large enough such that $\hat{y} := y_n$ satisfies
$$ \Vert y - \hat{y}\Vert^2_{L^2(W^\circ)} + \Vert \nabla y - \nabla \hat{y}\Vert^2_{L^2(W^\circ)} \le C_1\eps\rho, \ \ \ {\cal H}^1(J_{\hat{y}}) \le C_1\rho$$
for $C_1>0$ large enough. Choose a finite number of closed segments $(S_i)^m_i$ such that $\overline{J_{\hat{y}}} \cap  W^\circ \subset \bigcup_i S_i$ and ${\cal H}^1(\bigcup_i S_i) \le C_1\rho$. For $s>0$ small choose $T_i \in {\cal U}^s$ as the smallest  rectangle with $S_i \subset T_i$. Then by Lemma \ref{rig-lemma: modifica}(i) we obtain a set $\tilde{\Omega} \in {\cal V}^s_{\eps}$ with 
$$ \big|\tilde{\Omega}  \triangle \big(W \setminus \bigcup\nolimits^m_{j=1} T_m\big)\big|=0.$$
Observe that for $s$ sufficiently small we obtain $\Vert \tilde{\Omega}\Vert_* \le \Vert W \Vert_* + C_1\rho$ and $|W \setminus \tilde{\Omega}| \le C_1 \rho$. This together with \eqref{rig-eq: v4.3} gives the two first parts of \eqref{rig-eq: V3.20}. Finally, define the function $\tilde{y} \in H^1(\tilde{\Omega})$ by $\tilde{y} = \hat{y}|_{\tilde{\Omega}}$ and observe that $\tilde{y}$ satisfies \eqref{rig-eq: V3.20}. \eop

\section{Proof of the main SBD-rigidity result}\label{rig-sec: sub, proof-main}

This last section is devoted to the proof of the main SBD-rigidity result.  We start with some preparations and then split up the proof into two steps concerning a suitable construction of the jump set and the definition of an extension.   As before constants indicated by $C_1$ only depend on $M,\eta,\Omega$ and all constants  do not depend on $\rho$ and $q$ unless stated otherwise. 

Let $y \in SBV_M(\Omega) \cap L^2(\Omega)$ be given and let $\rho > 0$, $\varrho = \rho^{q}$ for $q \in \N$ to be specified below. Set $k = \rho^{q-1}$ and $m=\rho$. Recall the definition $\Omega_\rho = \lbrace x \in \Omega: \dist(x,\partial \Omega) > C\rho \rbrace$. We apply Theorem \ref{rig-th: rigidity2} and obtain a set $\Omega_y \subset \Omega_{\rho}$ with $\Omega_y \in {\cal V}^{s}_{ck}$ for $s$ sufficiently small  and  $|\Omega \setminus \Omega_y| \le C_1\rho$ such that  \eqref{rig-eq: main properties} holds  for a modification $\tilde{y} \in H^1(\Omega_y) \cap SBV_{cM}(\Omega_y)$ with $\Vert y - \tilde{y} \Vert^2_{L^2(\Omega_y)} + \Vert \nabla y - \nabla \tilde{y} \Vert^2_{L^2(\Omega_y)}\le C_1\eps\rho$.  Recall from the proof of Theorem \ref{rig-thm: V1} and Corollary \ref{rig-cor: weakA} that there is a set $\Omega_y^H \in {\cal V}^{3\varrho}_{ck}$ with $ \Omega^\circ_y \subset\Omega_y^H$ and  an extension $\hat{y}: \Omega_y^H \to \R^2$ of $\tilde{y}$ satisfying  \eqref{rig-eq: corweak**} and estimates of the form \eqref{rig-eq: corweak}.

We first construct a modification of $\Omega_y^H$ and appropriate Jordan curves which separate the connected components. For a (closed) Jordan curve $\gamma$ we denote by ${\rm int}(\gamma)$ the interior of the curve.  As connected components may be not simply connected we further introduce a generalization: We say a curve $ \gamma = \gamma_0 \cup{\bigcup}^m_{j=1} \gamma_j$ is a \emph{generalized Jordan curve} if it consists of pairwise disjoint Jordan curves $\gamma_0, \ldots, \gamma_m$ with $\gamma_j \in {\rm int}(\gamma_0)$ for $j=1,\ldots,m$. We define the interior of $\gamma$ by ${\rm int}(\gamma) = {\rm int}(\gamma_0) \setminus \bigcup^m_{j=1} {\rm int}(\gamma_j)$.

\begin{lemma}\label{rig-lemma: jordan}
Let $\rho > 0$, $M>0$ and $q \in \N$. There is a constant $C_1=C_1(M)>0$  such that for all $\Omega_y^H \in {\cal V}^{3\varrho}_{ck}$ as given above we find  $\hat{\Omega} \subset \Omega_\rho$ with ${\cal H}^1(\partial \hat{\Omega})\le  C_1$, $|\Omega_y^H \setminus \hat{\Omega}| \le C_1\rho$
and a set $S \subset \Omega_\rho \setminus \hat{\Omega}$ such that
\begin{itemize}
\item[(i)]  ${\cal H}^1(S) \le \Vert \Omega^H_y\Vert_* + C_1\rho$,
\item[(ii)] $\text{for all } \hat{P}_i \text{ there is a  generalized Jordan curve }  \gamma  \text{ in }  S \cup \partial \Omega_\rho \text{ such that } \\ {\rm int}(\gamma) \cap \hat{\Omega} = \hat{P}_i, \text{where } (\hat{P}_i)_i \text{ denote the connected components of } \hat{\Omega}$,
\item[(iii)] ${\rm int}(\gamma) \cap \hat{\Omega} \neq \emptyset \text{ for all Jordan curves } \gamma  \text{ in }  S \cup \partial \Omega_\rho,$
\item[(iv)] $\dist(x,S) \le C_1\rho^{q-2} \text{ for all } x \in \Omega_\rho \setminus \hat{\Omega}$,
\item[(v)] $(S\cup \partial \Omega_\rho) \cap X_t(\hat{\Omega}) \text{ is connected for all components } X_t(\hat{\Omega}) \text{ of } \Omega_\rho \setminus \hat{\Omega}.$
\end{itemize}
\end{lemma}

\Proof  In contrast to the previous sections, where it was essential to  avoid the combination of different cracks, we now combine boundary components: Choose a set $\hat{\Omega}_y^H \in {\cal V}^{3\varrho}$ satisfying $\hat{\Omega}_y^H  \subset \Omega_y^H$, $|\Omega_y^H \setminus \hat{\Omega}_y^H| = 0$ and $|\Gamma_j(\hat{\Omega}_y^H ) \cap \Gamma_l(\hat{\Omega}_y^H )|_{\cal H} =  0$ for $j \neq l$.  Clearly, by  \eqref{rig-eq: V1.7} and \eqref{rig-eq: abs crac} we have ${\cal H}^1(\hat{\Omega}_y^H) \le {\cal H}^1(\Omega_y^H)\le C_1$.

 Letting $Y_1, \ldots,Y_m$ be the connected components of $\hat{\Omega}_y^H$ satisfying $|\partial Y_j|_\infty \le  \rho^{q-2} $ for $j=1,\ldots,m$ we define $\tilde{\Omega} = \hat{\Omega}_y^H  \setminus \bigcup\nolimits^m_{j=1} Y_j$.  As $|\partial Y_j|_\infty \le \rho^{q-2}$ for $j=1,\ldots,m$, the isoperimetric inequality implies $|\bigcup_{j=1}^m Y_j| \le C_1\rho^{q-2}\Vert \hat{\Omega}_y^H\Vert_{\cal H} \le C_1\rho$  and thus  $|\Omega_y^H \setminus \tilde{\Omega}| \le C_1\rho$. 

 Let $Z \subset \Omega_\rho \setminus  \tilde{\Omega}$ be the largest set in ${\cal U}^{\rho^{q-2}}$ such that $\dist_\infty(x,\partial  \tilde{\Omega} \setminus \partial \Omega_\rho) \ge  \rho^{q-2} $ for all $x \in Z$ and define $\hat{\Omega} = \tilde{\Omega} \cup  \overline{Z}.$ (Observe that $Z$ is typically not connected.)  It is not hard to see that
\begin{align}\label{rig-eq: Z,Y} 
\dist(x,\partial \hat{\Omega} \setminus \partial \Omega_\rho) \le C_1 \rho^{q-2} \ \ \  \text{for all}  \ \ x \in \Omega_\rho \setminus \hat{\Omega}.
\end{align} 
Moreover, we get $|\Omega_y^H \setminus \hat{\Omega}| \le C_1\rho$ and ${\cal H}^1(\hat{\Omega}) \le C_1$. In fact, for each connected component $Z^i$ of $\overline{Z}$ we find boundary components $(X^i_j=X^i_j(\Omega_y^H))_j$ and $(Y_j^i)_j$ such that $\partial Z^i \subset \bigcup_j \overline{X^i_j} \cup \bigcup_j \overline{Y_j^i} $ and thus by $|\partial X^i_j|_\infty \le c\rho^{q-1}$, $|\partial Y^i_j|_\infty \le \rho^{q-2}$ we obtain $|\partial Z^i|_{\cal H} \le C_1 (\sum_j |\partial X^i_j|_{\cal H} +\sum_j |\partial Y^i_j|_{\cal H})$. We recall ${\cal H}^1(\Omega_y^H)\le C_1$ and observe that for different components $Z^{i_1}, Z^{i_2}$ one has $(\bigcup_j \overline{X^{i_1}_j} \cup \bigcup_j \overline{Y^{i_1}_j})\cap (\bigcup_j \overline{X^{i_2}_j} \cup \bigcup_j \overline{Y^{i_2}_j}) = \emptyset$.

Let $\hat{P}_1,\ldots, \hat{P}_n$ be the connected components of $\hat{\Omega}$ and define ${\cal F}(\hat{P}_i) = \lbrace X_j = X_j(\Omega_y^H): \overline{X_j} \cap \overline{\hat{P}_i} \neq \emptyset\rbrace$. (Here it is essential that we take the components of $\Omega_y^H$.)  By $Z_j \in {\cal U}^{3\varrho}$ we denote the smallest rectangle containing $X_j$. 

(I)  As a preparation we  consider the special case that there is only one connected component $\hat{P}_1$.  Moreover, we first suppose that $\Omega_\rho \setminus \hat{\Omega}$ consists of one connected component only. We can choose a set $S$ in $\bigcup_{Z_j \in {\cal F}(\hat{P}_1)} \overline{Z_j}$ consisting of segments such that  $S\cup (\partial \Omega_\rho \setminus \hat{\Omega})$ is connected, 
\begin{align}\label{rig-eq: S construction2}
{\cal H}^1(S ) \le  (1 + C_1\rho)\sum\nolimits_{X_j \in {\cal F}(\hat{P}_1)} |\Gamma_j|_\infty \le  (1 + C_1\rho)\sum\nolimits_{X_j \in {\cal F}(\hat{P}_1)} |\Gamma_j|_*
\end{align}
and $\dist(x,S) \le C_1 \rho^{q-2} $ for all  $x \in  \partial \hat{P}_1  \setminus \partial \Omega_\rho$ for a sufficiently large constant. Indeed, a set with the desired properties can be constructed in the following way.  By the definition of $|\cdot|_\infty$ we first see that we can choose a piecewise affine  Jordan curve $\gamma$ in $\bigcup_{X_j \in {\cal F}(\hat{P}_1)} \overline{Z_j} \cup \partial \Omega_\rho$ such that $\hat{P}_1 \subset {\rm int}(\gamma)$ and $S_0 := \gamma \cap \Omega^\circ_\rho$ satisfies ${\cal H}^1(S_0 ) \le \sum\nolimits_{X_j \in {\cal G}(S_0)} |\Gamma_j|_{ \infty}$, where ${\cal G}(S_0) = \lbrace X_j : X_j \cap S_0 \neq \emptyset\rbrace$. (If $\gamma \cap \Omega_\rho^\circ  = \emptyset$, we let $S_0  =  \lbrace p_0 \rbrace$ for some point $p_0 \in \Omega_\rho \setminus \hat{\Omega}$.) Assume a connected set $S_l$ consisting of segments has been constructed such that 
\begin{align}\label{rig-eq: S construction}
{\cal H}^1(S_l) \le \sum\nolimits_{X_j \in {\cal G}(S_l)} |\Gamma_j|_{\infty} + C_1l \rho^{q-1}.
\end{align}
 If $\dist(x,S_l) \le C_1 \rho^{q-2}$ for all  $x \in \partial \hat{P}_1   \setminus \partial \Omega_\rho $, we stop. Otherwise,  there is some $y \in  \partial \hat{P}_1  \setminus \partial \Omega_\rho$ such that $\dist(y,S_l) > C_1\rho^{q-2}$. By the definition of $|\cdot|_\infty$ it is elementary to see that we can find a piecewise affine, continuous curve $T_{l+1}$ with $T_{l+1} \cap S_{l} \neq \emptyset$, $y \in T_{l+1}$, $\# ({\cal G}(T_{l+1}) \cap {\cal G}(S_l)) = 1$ such that  ${\cal H}^1(T_{l+1}) \le \sum_{X_j \in {\cal G}(T_{l+1})} |\Gamma_j|_{\infty}$. Then using that  $|\Gamma(\Omega_y^H)|_\infty \le 2\sqrt{2} \cdot ck \le C_1\rho^{q-1}$ and  $\# ({\cal G}(T_{l+1}) \cap {\cal G}(S_l)) = 1$ we find that \eqref{rig-eq: S construction} is  satisfied for $S_{l+1} := S_l \cup T_{l+1}$. 

 After a finite number of iterations $n \in \N$ we find that $\dist(y,S_n) \le C_1\rho^{q-2}$ for all $y \in \partial \hat{P}_1  \setminus \partial \Omega_\rho$ and set $ S_*= S_n$. Indeed, this follows from the fact that in each iteration ${\cal G}(S_l)$ increases and the assertion clearly holds if $S_l$ intersects all boundary components since $\max_j|\Gamma_j(\Omega_y^H)|_\infty \le  C_1\rho^{q-1}$.  As $ {\cal H}^1(T_l) > C_1\rho^{q-2}$, it is not hard to see that $n \le C_1\rho^{2-q}\sum\nolimits_{X_j \in {\cal F}(\hat{P}_1)} |\Gamma_j|_\infty$ and thus \eqref{rig-eq: S construction2} holds  replacing $S$ by $S_*$. 
 
 Observe that possibly $S_* \cup (\partial \Omega_\rho \setminus \hat{\Omega})$ is not connected. Therefore, we choose some point $y$ in each connected component of $\partial \Omega_\rho \setminus \hat{\Omega}$ (which may be several if $\Omega_\rho$ is not simply connected) and repeat the construction below \eqref{rig-eq: S construction} for each $y$. We obtain a set $S$ such that \eqref{rig-eq: S construction2} still holds and $S \cup (\partial \Omega_\rho \setminus \hat{\Omega})$ is connected.
 
 If $\Omega_\rho \setminus \hat{\Omega}$ consists of several connected components $X_t(\hat{\Omega})$, we repeat the arguments on each component separately possibly starting with $S_0 = \lbrace p_0 \rbrace$ for some $p_0 \in X_t(\hat{\Omega})$.

 We see that (i),(v) are satisfied, (ii)  holds with $\gamma$  and (iii) follows from the fact that in the construction of the sets $T_l$ above we do not obtain additional `loops'. Moreover, (iv) follows from the fact that each $x \in \Omega_\rho \setminus \hat{\Omega}$ satisfies $\dist(x, \partial \hat{P}_1  \setminus \partial \Omega_\rho) \le C_1 \rho^{q-2}$ by  \eqref{rig-eq: Z,Y}.

(II)  We now consider  an arbitrary number of connected components. Choose Jordan curves $\gamma^i$ in $\bigcup_{X_j \in {\cal F}(\hat{P}_i)} \overline{Z_j} \cup \partial \Omega_\rho$ such that $\hat{P}_i   \subset {\rm int}(\gamma^i) \cap \hat{\Omega}$ and  ${\cal H}^1(\gamma^i\cap  \Omega^\circ_\rho ) \le \sum\nolimits_{X_j \in {\cal G}(\gamma^i )} |\Gamma_j|_\infty$.  We first assume that $\hat{P}_i  = {\rm int}(\gamma^i) \cap \hat{\Omega}$, i.e. ${\rm int}(\gamma^i)$ does not contain other components of $\hat{\Omega}$, and treat the general case in (III). As the sets $({\cal F}(\hat{P}_i))_{i=1}^n$ might be not disjoint, we have to combine the different curves in a suitable way.  Define $G_i= \bigcup_{X_j \in {\cal G}(\gamma^i)}  \overline{Z_j}$. It is not restrictive to assume that $\bigcup_{1 \le i \le n}  G_i$ is  connected  as otherwise we apply the following arguments on each component separately.  For $B \subset \R^2$ we define $${\rm Int}(B) = \lbrace x \in \R^2: \exists \ {\rm Jordan \  curve } \ \gamma^i \ {\rm in } \ B: x \in {\rm int}(\gamma^i) \rbrace.$$
  It is not hard to see that we can order the sets $(\hat{P}_i)_i$ in a way such that for all $1 \le l \le n$ we have  $\bigcup_{1 \le i \le l}  G_i$ is connected and  ${\rm Int}(\bigcup_{1 \le i \le l}  G_i) \cap \hat{P}_j = \emptyset$ for all $j > l$. In fact, to see the second property, assume the first $l$ sets $G_1, \ldots, G_l$ have already been chosen. Select some other component $\hat{P}_k$, $k>l$, with corresponding $G_k$. If the desired property is satisfied, we reorder and set $G_{l+1} = G_k$, otherwise we find some $\hat{P}_{k'}$, $k'>l, k' \neq k$, with corresponding $G_{k'}$ such that $\hat{P}_{k'} \subset {\rm Int}(\bigcup_{1 \le i \le l} G_i \cup G_k)$. Possibly repeating this procedure we finally find a set $G_{l+1}$ such that ${\rm Int}(\bigcup_{1 \le i \le l+1} G_i) \cap \hat{P}_j = \emptyset$ for all $j > l+1$. 
  
  We now proceed iteratively. Set $ S_0 = \emptyset$ and assume a connected set $S_l$ has been constructed with
\begin{align}\label{rig-eq: mayblast 13}
(a) & \ \ {\cal H}^1(S_l \cap \Omega_\rho) \le  (1+ C_1\rho)\sum\nolimits_{X_j \in \bigcup_{1\le i \le l} {\cal 
G}({\rm int}(\gamma^i))} |\Gamma_j|_* + C_1( l-1)  \rho^{q-1},\notag\\
(b)& \ \ \text{for all } 1 \le i \le l \text{ there is a  Jordan curve }  \gamma  \text{ in }  S_l  \text{ such that }  {\rm int}(\gamma) \cap \hat{\Omega} = \hat{P}_i,\notag\\
(c)& \ \ \dist(x,S_l) \le C_1 \rho^{q-2} \text{ for all } x \in \bigcup\nolimits^l_{i=1} \partial \hat{P}_i  \setminus \partial \Omega_\rho.
\end{align} 
  Let $T_{l+1}$ be the (unique)  connected component of $\gamma^{l+1} \setminus \bigcup_{1 \le i \le l}  G_i$ such that  $\hat{P}_{l+1} \subset {\rm Int}(\bigcup_{1 \le i \le l}  G_i \cup T_{l+1})$. Now choose two segments $T^j_{l+1}$, $j=1,2$, with ${\cal H}^1(T^j_{l+1}) \le C_1 \rho^{q-1}$, $T^j_{l+1} \cap S_l \neq \emptyset$,   $T^j_{l+1} \cap T_{l+1} \neq \emptyset$ such that $\hat{S}_{l+1} := S_l \cup T_{l+1} \cup \bigcup_{j=1,2} T^j_{l+1}$ satisfies $\hat{P}_{l+1} \subset {\rm Int}(\hat{S}_{l+1})$ and
$$ {\cal H}^1(\hat{S}_{l+1} \cap \Omega_\rho) \le   (1+ C_1\rho)\sum\nolimits_{X_j \in \bigcup_{1\le i \le l} {\cal 
G}({\rm int}(\gamma^i)) \cup {\cal G}(T_{l+1})} |\Gamma_j|_* + C_1 l \rho^{q-1}.$$
 By the order of the sets $(\hat{P}_i)_i$ it is not hard to see that there is a Jordan curve $\gamma$ in $ \hat{S}_{l+1}$ with ${\rm int}(\gamma) \cap \hat{\Omega} = \hat{P}_{l+1}$.  Observe that  $\dist(x,\gamma^{l+1}) \le C_1 \rho^{q-2}$ for all $x \in \partial \hat{P}_{l+1}  \setminus \partial \Omega_\rho$ might not hold. Therefore, following the lines of (I) we choose a (possibly not connected) set $R_{l+1} \subset \text{int}(\gamma^{l+1})$ such that  such that $S_{l+1} := \hat{S}_{l+1} \cup R_{l+1}$ is connected  in each component of $\Omega_\rho \setminus \hat{\Omega}$, $\dist(x,S_{l+1}) \le C_1\rho^{q-2}$ for all $x \in \partial \hat{P}_{l+1}  \setminus \partial \Omega_\rho$ and
$${\cal H}^1(R_{l+1}) \le (1+ C_1\rho)\sum\nolimits_{X_j \in  {\cal G}({\rm int}(\gamma^{l+1})) \setminus {\cal G}(\hat{S}_{l+1}) } |\Gamma_j|_*.$$  
Now it is not hard to see that (a)-(c) are satisfied for $S_{l+1}$.

After the last iteration step we define $ S_* = S_n \cap \Omega_\rho$. Observe that  by construction (see before \eqref{rig-eq: Z,Y}) each $\hat{P}_i$ satisfies $|\partial \hat{P}_i|_\infty\ge \rho^{q-2} $. Thus $n \le C_1  \rho^{2-q} $ and then we obtain ${\cal H}^1(S_*) \le \Vert \Omega_y^H\Vert_* + C_1 \rho$ since $n  \rho^{q-1}  \le C_1\rho$.  Similarly as before, $S_* \cup \partial \Omega_\rho$ might not be connected in the components of $\Omega_\rho \setminus \hat{\Omega}$. Consequently, we proceed as in (I) (see construction below \eqref{rig-eq: S construction}) to find a set $S \supset S_*$ such that (i) still holds and $S \cup \partial \Omega_\rho$ is connected in the components of $\Omega_\rho \setminus \hat{\Omega}$. This gives (v). Moreover, (b) implies (ii)  and similarly as in (I) also (iii) holds.  (Here we do not have to consider generalized Jordan curves.)  Finally, to see (iv) we use (c) and the fact that each $x \in \Omega_\rho \setminus \hat{\Omega}$ satisfies $\dist(x, \partial \hat{\Omega} \setminus \partial \Omega_\rho) \le C_1 \rho^{q-2} $ by  \eqref{rig-eq: Z,Y}. 

(III) We now finally treat the case that the components $(\hat{P}_i)^n_{i=0}$ may also contain other components of $\hat{\Omega}$. To simplify the exposition we assume that there is exactly one component, say $\hat{P}_0$, such that 
$\hat{P}_0  \neq {\rm int}(\gamma^0) \cap \hat{\Omega}$. The general case follows by inductive application of the following arguments. 

We proceed as in (II) (assuming we had $\hat{P}_0  = {\rm int}(\gamma^0) \cap \hat{\Omega}$) and construct a set $S'$ particularly satisfying (i),(iii),(v). We have to verify (ii) for $\hat{P}_0$ and find a set $S \supset S'$ such that (iv) is satisfied and (i),(iii),(v) still hold. By $(\hat{P}_{i_j})_j$ we denote the components with $\hat{P}_{i_j} \subset {\rm int}(\gamma_0)$. As (ii) holds for these components we find pairwise disjoint Jordan curves $\gamma_1,\ldots, \gamma_m$ with $\bigcup_j \hat{P}_{i_j} \subset \bigcup_{j=1}^m {\rm int}(\gamma_j) \subset {\rm int}(\gamma_0)$. Consequently, defining the generalized Jordan curve $\gamma = \bigcup^m_{j=0} \gamma_j$ we find $\hat{P}_0  = {\rm int}(\gamma) \cap \hat{\Omega}$ which gives (ii).

 Let $(Y_j)_j$ be the components of $\Omega_\rho \setminus \hat{\Omega}$ which are completely contained in ${\rm int}(\gamma_0)$. We observe that (iv) may be violated for $x \in Y^*:=\bigcup_j Y_j \setminus \bigcup_{j=1}^m {\rm int}(\gamma_j)$. We now proceed similarly as in (I) to obtain a set $R \subset Y^* $ such that $S : = S' \cup R$ is connected in the connected components of $\Omega_\rho \setminus \hat{\Omega}$ and $\dist(x,S) \le C_1 \rho^{q-2}$ for all $x \in \partial \hat{P}_1 \cap Y^*$. This implies (iii),(v) and together with \eqref{rig-eq: Z,Y}  also (iv). Arguing similarly as in (II) we find that (i) is still satisfied since the sum in \eqref{rig-eq: mayblast 13}(a) does not run over the components contained in $Y^*$.  \eop

We finally can give the proof of Theorem \ref{rig-th: rigidity} by constructing an extension $\hat{y}$ of the function $\tilde{y}$. We briefly note that the function $\hat{y}$ has to be defined as an extension of the approximation and not of the original deformation $y$ as only in this case we obtain the correct surface energy due to the higher regularity of the jump set of $\tilde{y}$  and the available trace estimates. Recall the definition of $E_ \eps^\rho(y,U)$ in \eqref{rig-eq: Griffith en2}, in particular $f_\eps^\rho(x) = \min\lbrace\frac{x}{\sqrt{\eps}\rho} ,1 \rbrace$.

\noindent {\em Proof of Theorem \ref{rig-th: rigidity}.}  Let  $\Omega_y \subset \Omega_{\rho}$ with $\Omega_y \in {\cal V}^{s}$ and $\Omega_y^H \in {\cal V}^{3\varrho}$ with $ \Omega^{ \circ}_y \subset\Omega_y^H $ be given.  Recall that $|\Omega \setminus \Omega_y| \le C_1 \rho$.  Let $\tilde{y} \in H^1(\Omega_y) $ be the approximation of $y \in SBV_M(\Omega)$ with $\Vert y - \tilde{y} \Vert^2_{L^2(\Omega_y)} + \Vert \nabla y - \nabla \tilde{y} \Vert^2_{L^2(\Omega_y)} \le C_1\eps\rho$ and let $\hat{y} \in  SBV_{cM}(\Omega^H_y)\cap L^2(\Omega^H_y)$ be the extension of $\tilde{y}$ given by Corollary \ref{rig-cor: weakA}. Let $\hat{\Omega}$ be the set constructed in Lemma \ref{rig-lemma: jordan}.  We first consider the jumps  of $\hat{y}$ in $(\Omega_y^H \cap \hat{\Omega})^\circ$. By \eqref{rig-eq: weaklocA3_**} and H\"older's inequality we find
\begin{align*}
\Big(\int_{J_{\hat{y}} \cap (\Omega_y^H)^\circ} |[\hat{y}]|\,d{\cal H}^1\Big)^2 &\le \Big(\sum\nolimits_{Q_t \subset \Omega_y^H} \int_{J_{\hat{y}} \cap \overline{ Q_t} } |[\hat{y}]|\,d{\cal H}^1\Big)^2 \\
& \le \sum\nolimits_{Q_t } | J_{\hat{y}} \cap  \overline{Q_t}|_{\cal H} \cdot \sum\nolimits_{Q_t} | J_{\hat{y}}\cap  \overline{Q_t}|_{\cal H}^{-1} \Big( \int_{J_{\hat{y}} \cap \overline{Q_t} } |[\hat{y}]|\,d{\cal H}^1\Big)^2 \\
& \le  C{\cal H}^1(J_{\hat{y}}) \cdot \sum\nolimits_{Q_t } CC_\rho^2 \varrho^2 ( \gamma(N_t) + \delta_4(N_t ) + \epsilon|\partial  W \cap  N_t|_{\cal H}),
\end{align*}
where  $W$ as defined in \eqref{rig-eq: V1.41}, $N_t := N(Q_t)= \lbrace x \in W: \dist(x, Q_t) \le C\rho^{q-1}\rbrace$  and  $\gamma(N_t) = \Vert \nabla  \dist(\nabla \hat{y},SO(2))\Vert^2_{L^2(W)}$, $\delta_4(N_t) = \sum_{i=1}^4 \Vert \nabla  \hat{y} - \hat{R}_i  \Vert^4_{L^4(W)}$ (recall \eqref{rig-eq: V1.11}). As each $x \in \Omega$ is contained in at most $\sim \rho^{-2}$ different $N_t$ we find by \eqref{rig-eq: w energy}, \eqref{rig-eq: V1.41}, \eqref{rig-eq: V1.11},  \eqref{rig-eq: corweak**} and the fact that $\epsilon =  \hat{c} \rho^{-1} \eps$
\begin{align*}
\Big(\int_{J_{\hat{y}} \cap (\Omega_y^H \cap \hat{\Omega})^\circ} |[\hat{y}]|\,d{\cal H}^1\Big)^2 \le  C\rho^{-2} \cdot CC_\rho^2\varrho^2 \epsilon \le C \varrho^2\rho^{-3}C_\rho^2 \eps.
\end{align*}
(Note that in the general case the set $W$ and the rigid motions $\hat{R}_i$ were defined differently (see e.g. \eqref{rig-eq: R4}), but here and in the following we prefer to refer to the proof of Theorem \ref{rig-thm: V1} for the sake of simplicity.)
By Remark \ref{rig-rem: z}(i) we get for $q = q(h_*)$ sufficiently large
$$\int_{J_{\hat{y}} \cap (\Omega_y^H \cap \hat{\Omega})^\circ} |[\hat{y}]|\,d{\cal H}^1 \le CC_{\rho} \rho^{ q- \frac{3}{2}}\sqrt{\eps} =  C\rho^{q-(\frac{3}{2} + z)}\sqrt{\eps} \le \rho^{2} \sqrt{\eps}.$$
Recalling that $f_\eps^{\rho}(x) \le \rho^{-1} \frac{x}{\sqrt{\eps}}$ for $x \ge 0$ we get 
\begin{align}\label{rig-eq: rig9}
\int_{J_{\hat{y}} \cap (\Omega_y^H \cap \hat{\Omega})^\circ}  f_\eps^{\rho}(|[\hat{y}]|)\,d{\cal H}^1  \le \eps^{-1/2} \rho^{-1}\int_{J_{\hat{y}} \cap (\Omega_y^H \cap \hat{\Omega})^\circ}  |[\hat{y}]|\,d{\cal H}^1  \le  \rho.
\end{align}
We now concern ourselves with the components of $\partial \hat{\Omega}$. Let $Y_t$ be a connected component of $\Omega_\rho \setminus (\hat{\Omega} \cup S)$, where $S$ is the set constructed in Lemma \ref{rig-lemma: jordan}. Set $S_t = S \cap \overline{Y_t}$ and $\Gamma_t = \overline{Y_t} \cap \partial \hat{\Omega}$.  We observe  that  by Lemma \ref{rig-lemma: jordan}(ii),(iii) $\Gamma_t$  is a Jordan curve if $\overline{Y_t} \cap \partial \Omega_\rho = \emptyset$. 

 Define $J  = I^{\varrho}(\hat{\Omega})$ and for $\Gamma_t$ we choose   $J(\Gamma_t) \subset J$ such that $\overline{Q^{{ \varrho} }(p)} \cap \Gamma_t \neq \emptyset$ for all $p \in J(\Gamma_t)$. We set $M(\Gamma_t) = \bigcup_{p \in J(\Gamma_t)} \overline{Q^{ \varrho}(p)}$. For later purpose, for components with $|\Gamma_t|_{\infty} >   2\rho^{q-2}$ we introduce a finer partition of $M(\Gamma_t)$: Define $J(\Gamma_t)= I_1 \dot{\cup} \ldots \dot{\cup} I_n$ and the connected sets $B_i =  \bigcup_{p \in I_i} \overline{Q^{\varrho}(p)} $ such that $\rho^{-2} \le \# I_i \le C\rho^{-2}$,   $i=1,\ldots,n$, for a constant $C \gg 1$. For  $|\Gamma_t|_{\infty} \le  2\rho^{q-2}$ we let $I_1 = J(\Gamma_t)$.  It is elementary to see that $n \le \max\lbrace C|\Gamma_t|_{\cal H}\rho^{2-q} , 1\rbrace \le C|\Gamma_t|_{\cal H}\rho^{-q}$, where we used $|\Gamma_t|_{\cal H} \ge C\rho^q$.

Consider  $\bar{R}_j: \Omega^H_y \to SO(2)$ and $\bar{c}_j: \Omega_y^H  \to \R^2$, $j=1,\ldots,4$, as given in \eqref{rig-eq: V1.9.2}. Recall the definition $\tilde{\Omega} = \hat{\Omega} \setminus \overline{Z} \subset \Omega_y^H$ before \eqref{rig-eq: Z,Y}. We extend the function $\hat{y}$ to $\hat{\Omega}$ by setting $\hat{y} = \id$ on $\hat{\Omega} \setminus  \tilde{\Omega}$ and likewise let $\bar{R}_j = \Id$, $\bar{c}_j = 0$ on $\hat{\Omega} \setminus  \tilde{\Omega}$.  (If $\overline{Z} \cap \Omega_y^H \neq \emptyset$, we redefine the function on this set.)  Applying Corollary \ref{rig-cor: weakA} on each $Q^{ 3\varrho}_j(p) \subset \hat{\Omega}$ with $Q^{3\varrho}_j(p) \cap M(\Gamma_t) \neq \emptyset$, we get
\begin{align}\label{rig-eq: rig3}
\begin{split}
\Vert \hat{y} - (\bar{R}_j\, \cdot + \bar{c}_j) \Vert^2_{L^2(B_i)} &\le C\varrho^2 C_\rho^2 \cdot \rho^{-2}  \rho^{q- 1}  \epsilon \cdot \#I_i = C\rho^{3q-6}C_\rho^2 \eps, \\
\Vert \hat{y} - (\bar{R}_j\, \cdot + \bar{c}_j) \Vert^2_{L^1( \partial B_i)} &\le  C\rho^{3q-6}C_\rho^2 \eps, 
\end{split} 
\end{align}
for $j=1,\ldots,4$ and $i=1,\ldots,n$. Here we used   $k=\rho^{q-1}$, $\epsilon =  \hat{c}\eps\rho^{-1}$ and the fact that each $N(Q^{ 3\varrho}_{j}(p))$ contains  $\sim m^{-2} = \rho^{-2}$ different $Q^{ 3\varrho}(p)  \subset \Omega_y^H$.  The triangle inequality then yields
\begin{align*}
\Vert (\bar{R}_{j_1}\, \cdot + \bar{c}_{j_1})- (\bar{R}_{j_2}\, \cdot + \bar{c}_{j_2}) \Vert^2_{L^2(B_i)} \le C\rho^{3q-6}C_\rho^2 \eps
\end{align*}
for $1 \le j_1,j_2 \le 4$ and $i=1,\ldots,n$.  The strategy will be to cover $Y_t$ with $n$ different rigid motions.  We argue as in \eqref{rig-eq: rig2}f. and  \eqref{rig-eq: A diff}  to get $\hat{R}_i \in SO(2)$, $\hat{c}_i \in \R^2$ such that
\begin{align*}
\Vert \hat{y} - (\hat{R}_i\, \cdot + \hat{c}_i) \Vert^2_{L^2(B_i)} \le C (\# I_i)^4 \rho^{3q-6}C_\rho^2 \eps   \le C \rho^{3q-{ 14}}C_\rho^2 \eps.
\end{align*}
Here we used H\"older's inequality (cf. \eqref{rig-eq: A diff}).
A similar argument shows that we even find
\begin{align}\label{rig-eq: rig4}
 \sum\nolimits _{j=-1,0,1}\Vert \hat{y} - (\hat{R}_{i+j}\, \cdot + \hat{c}_{i+j}) \Vert^2_{L^2(B_i)} \le C \rho^{3q-14}C_\rho^2 \eps 
 \end{align}
 for $i=1,\ldots,n$,  where (in the case that $\Gamma_t$ is a Jordan curve) we set $\hat{R}_{n+1} = \hat{R}_1$, $\hat{c}_{n+1} = \hat{c}_1$ and $\hat{R}_{0} = \hat{R}_n$, $\hat{c}_{0} = \hat{c}_n$. Without restriction  recalling Remark \ref{rig-rem: 1}(iii) we can assume that $\hat{R}_i \in \text{im}_{\bar{R}_4}(M(\Gamma_t)) \subset SO(2)$, where $\text{im}_{\bar{R}_4}$ denotes the image of the function $\bar{R}_4$.   For shorthand let $\bar{R} = \bar{R}_4$ and $\bar{c} = \bar{c}_4$.  By \eqref{rig-eq: rig3} and \eqref{rig-eq: rig4} we get
\begin{align}\label{rig-eq: ri2}
 \sum\nolimits _{j=-1,0,1}\Vert (\hat{R}_{i+j}\, \cdot + \hat{c}_{i+j})- (\bar{R}\, \cdot + \bar{c}) \Vert^2_{L^2(B_i)} \le  C \rho^{3q-14}C_\rho^2 \eps.
\end{align}
Using H\"older's inequality and passing to the trace on each $Q^{3\varrho}(p)$ we obtain  for all $i=1, \ldots, n$
\begin{align*}
\begin{split}
 \sum\nolimits _{j=-1,0,1}\Vert (\hat{R}_{i+j}\, \cdot &+ \hat{c}_{i+j})- (\bar{R}\, \cdot + \bar{c}) \Vert^2_{L^1(B_i \cap \Gamma_t)} \\& \le C \sum\nolimits _{j=-1,0,1}\vert B_i \cap \Gamma_t \vert_{\cal H}\Vert (\hat{R}_{i+j}\, \cdot + \hat{c}_{i+j})- (\bar{R}\, \cdot + \bar{c}) \Vert^2_{L^2({ B_i \cap } \Gamma_t )} \\ & \le  C\varrho\rho^{ -2} \cdot\varrho^{-1}   \rho^{3q-14}C_\rho^2 \eps  \le C\rho^{3q-{ 16}}C_\rho^2 \eps.
\end{split}
\end{align*}
Together with  \eqref{rig-eq: rig3} this implies
\begin{align*}
\begin{split}
 \sum\nolimits _{j=-1,0,1}\Vert \hat{y} - (\hat{R}_{i+j}\, \cdot + \hat{c}_{i+j}) \Vert^2_{L^1(B_i \cap \Gamma_t)} &\le C\rho^{3q-16}C_\rho^2 \eps.
\end{split}
\end{align*}
This and the fact that  $n \le C|\Gamma_t|_{\cal H}\rho^{-q}$ yield
\begin{align}\label{rig-eq: rig7} 
\begin{split}
H_1 := \sum\nolimits_i  \sum\nolimits _{j=-1,0,1}  \Vert  \hat{y} - (\hat{R}_{i+j}\, &\cdot + \hat{c}_{i+j}) \Vert_{L^1(B_i \cap\Gamma_t)} \le C|\Gamma_t|_{\cal H}\rho^{ \frac{q}{2} - 8 }C_\rho \sqrt{\eps}.
\end{split}
\end{align}
For the difference of the rigid motions we get by the triangle inequality and \eqref{rig-eq: rig4}
$$   \sum\nolimits _{j_1,j_2=-1,0,1}\Vert (\hat{R}_{i+j_1}\, \cdot + \hat{c}_{i+j_1}) - (\hat{R}_{i+j_2}\, \cdot + \hat{c}_{i+j_2}) \Vert^2_{L^2(B_i)} \le  C \rho^{3q-14}C_\rho^2 \eps.$$
Let $\tilde{B}_i = \lbrace x \in \Omega: \dist(x,B_i) \le  \bar{C} \rho^{q-2}\rbrace$. Arguing similarly as in \eqref{rig-eq: A,c difference2} it is not hard to see that 
\begin{align}\label{rig-eq: ri1}
\begin{split}
   \sum\nolimits _{j_1,j_2=-1,0,1}\Vert (\hat{R}_{i+j_1}\, &\cdot + \hat{c}_{i+j_1}) - (\hat{R}_{i+j_2}\, \cdot + \hat{c}_{i+j_2}) \Vert^2_{L^2(\tilde{B}_i)} \\
& \le   C (\rho^{ -2})^2 \cdot \rho^{ -4} \cdot \rho^{3q-14}C_\rho^2 \eps \le C  \rho^{3q-22}C_\rho^2 \eps
\end{split}
\end{align}
as $\frac{|\tilde{B}_i|}{|B_i|} \le C\rho^{ -4}$ and $\frac{|\partial\tilde{B}_i|_\infty}{|\partial B_i|_\infty} \le C\rho^{ -2}$. Define  $\tilde{I}_i = I^{\varrho}(\tilde{B}_i )$. Again using H\"older's inequality, passing  from the traces to a bulk integral and recalling  $n \le C|\Gamma_t|_{\cal H}\rho^{ -q}$, $\# \tilde{I}_{i} \le C\rho^{ -4}$ we derive (let $\cdot =  (\hat{R}_{i+j_1}\, \cdot + \hat{c}_{i+j_1}) - (\hat{R}_{i+j_2}\, \cdot + \hat{c}_{i+j_2})$ for shorthand) 
\begin{align}\label{rig-eq: rig8}
\begin{split}
 H_2 & := \sum\nolimits_{i} \sum\nolimits_{p \in \tilde{I}_i}  \sum\nolimits_{j_1,j_2=-1,0,1} \Vert \cdot \Vert_{L^1(\partial Q^{ \varrho}(p))}  \\
 & \le C\sum\nolimits_{i} \sum\nolimits_{p \in \tilde{I}_i}  \sum\nolimits_{j_1,j_2=-1,0,1} \varrho^{1/2} \Vert \cdot \Vert_{L^2(\partial Q^{ \varrho}(p))}\\
 &\le  C \sum\nolimits_{i}  (\# \tilde{I}_i)^{\frac{1}{2}} \Big(\sum\nolimits_{p \in \tilde{I}_i}  \sum\nolimits_{j_1,j_2=-1,0,1} \varrho \Vert \cdot \Vert^2_{L^2(\partial Q^{ \varrho}(p))}\Big)^{1/2}\\
 & \le C \sum\nolimits_{i} \rho^{ -2} \Big(\sum\nolimits_{j_1,j_2=-1,0,1}\Vert \cdot \Vert^2_{L^2(\tilde{B}_i)}\Big)^{1/2}\le C|\Gamma_t|_{\cal H}\rho^{\frac{q}{2}-{ 13}}C_\rho \sqrt{\eps}.
\end{split}
 \end{align}
By $(T_j)_j$ we denote the connected components of $ Q^{\varrho}(p) \setminus  (\hat{\Omega} \cup S)$ for all $Q^{\varrho}(p)$ with  $Q^{\varrho}(p) \cap Y_t \neq \emptyset$.  We now choose suitable rigid motions:  Observe that $\dist(\Gamma_t  \cup \partial \Omega_\rho,x) \le C_1\rho^{ q-2}$ for all $x \in Y_t$ by Lemma \ref{rig-lemma: jordan}(iv) and the fact that $Y_t$ is a connected component of $\Omega_\rho \setminus (\hat{\Omega} \cup S)$. Therefore, for every $T_j$  with $\dist(T_j,\partial \Omega_\rho) \gg \rho^{ q-2}$ we find some (possibly non unique) $B_{i_j}$ with $\dist(T_j,B_{i_j})\le C\rho^{ q-2}$.  In particular, we get $T_j \subset \tilde{B}_{i_j}$  choosing $\bar{C}$ in the definition of $\tilde{B}_i$ large enough.  We define
\begin{align}\label{rig-eq: rig10_*}
 \hat{y}(x) = \hat{R}_{i_j}\, x + \hat{c}_{i_j} \ \ \ \text{for } x \in T_j  \cap Y_t \cap   \Omega_{2\rho}
 \end{align}
for all $j$ and note that we have found an extension $\hat{y}$ to $Y_t  \cap \Omega_{2\rho}$. (If $Y_t \cap \Omega_y^H \neq \emptyset$, we redefine the function on this set.)  Taking Lemma \ref{rig-lemma: jordan}(v) into account the choice of $B_{i_j}$ can be done in a way that for neighboring sets $T_1,T_2$ with $\overline{T}_1 \cap \overline{T}_2 \neq \emptyset$ one has $i_1 - i_2 \in \lbrace -1,0,1\rbrace$ and that ${\cal H}^1(J_{\hat{y}} \cap Y_t) \le C_1{\cal H}^1(\Gamma_t)$. Now by \eqref{rig-eq: rig7} and \eqref{rig-eq: rig8} it is not hard to see that
$$\int\nolimits_{(J_{\hat{y}} \cap \overline{Y_t}) \setminus S} |[\hat{y}]| \,d{\cal H}^1 \le CH_1 + CH_2  \le C|\Gamma_t|_{\cal H}\rho^{\frac{q}{2}- 13}C_\rho \sqrt{\eps}.$$
Repeating the arguments for all components $Y_t$ we obtain a  configuration $\hat{y} \in  SBV_{cM}(\Omega_{\rho})$  with $\hat{y} =\tilde{y}$ on $\Omega^*_y := \Omega_y \cap  \tilde{\Omega}$, where by Lemma \ref{rig-lemma: jordan} we have $|\Omega \setminus \Omega_y^*| \le C_1 \rho$. (With a slight abuse of notation we replace $\Omega_y^*$ by $\Omega_y$ in the assertion of Theorem \ref{rig-th: rigidity}.)  Summing over all $Y_t$ and recalling that ${\cal H}^1(\partial \hat{\Omega}) \le C_1$  by Lemma \ref{rig-lemma: jordan} we get 
\begin{align*}
\sum\nolimits_{t}  \int_{(J_{\hat{y}} \cap \overline{Y_t}) \setminus S } f_\eps^\rho(|[\hat{y}]|)\,d{\cal H}^1 &\le C \rho^{\frac{q}{2}-13}C_\rho  \le \rho
\end{align*}
for  $q = q(h_*)$ sufficiently large. Together with \eqref{rig-eq: rig9}, Lemma \ref{rig-lemma: jordan}(i) and \eqref{rig-eq: main properties}(i) this implies
$$\int_{J_{\hat{y}}} f_\eps^\rho(|[\hat{y}]|) \,d{\cal H}^1 \le \int_{J_{\hat{y}}\setminus S} f_\eps^\rho(|[\hat{y}]|) \,d{\cal H}^1 + {\cal H}^1(S)\le (1+ C_1h_*){\cal H}^1(J_y) + C_1\rho.$$
Choosing $h_* = \rho$ we finally get
\begin{align}\label{rig-eq: part + crack2}
\int_{J_{\hat{y}}} f_\eps^\rho(|[\hat{y}]|) \,d{\cal H}^1 \le {\cal H}^1(J_y) + C_1\rho.
\end{align}
 We observe $\nabla \hat{y} \in SO(2)$ on $\Omega_{\rho} \setminus \Omega_y$ (see construction in Corollary \ref{rig-cor: weakA}, \eqref{rig-eq: rig10_*} and recall $\hat{y} = \id$ in $\hat{\Omega} \setminus  \tilde{\Omega}$). As $\Vert \tilde{y} - y \Vert^2_{L^2(\Omega_y)} + \Vert \nabla \tilde{y} - \nabla y \Vert^2_{L^2(\Omega_y)} \le C_1\eps\rho$ we obtain $E_\eps^\rho(\hat{y},\Omega_{\rho}) \le E_\eps(y) + C_1\rho$ which gives \eqref{rig-eq: energy le}.  Here we used $\Vert \nabla  \tilde{y}\Vert_\infty + \Vert \nabla y\Vert_\infty \le cM$ and the regularity of the stored energy density $W$.

Let $(P_j)_j$ be the connected components of $\Omega_{\rho} \setminus S$. By Lemma \ref{rig-lemma: jordan}(ii),(iii) it is not hard to see that for every index $j$ there is a  (unique) connected component $\hat{P}_j$ of $\hat{\Omega}$ such that $\hat{P}_j \subset P_j$. Then there is either a connected component $P^H_j$ of $\Omega^H_y$ such that $\hat{P}_j = P^H_j$ (see proof of Theorem \ref{rig-thm: V1}) or $\hat{y} = \id$ on $\hat{P}_j$ (see construction before \eqref{rig-eq: rig3}). We now  define \eqref{rig-eq: u def2} by $u(x) = \hat{y}(x) - (R_j \, x +c_j)$ for $x \in P_j$, where $R_j \,x+c_j$ is either the rigid motion on $P_j^H$ given in Theorem \ref{rig-th: rigidity2} or $R_j = \Id$, $c_j = 0$, respectively. For later purpose, we note that for \eqref{rig-eq: part + crack2} we can also write 
\begin{align}\label{rig-eq: part + crack}
\sum\nolimits_j \tfrac{1}{2} P(P_j,\Omega_\rho) + \int_{J_{\hat{y}} \setminus \partial P} f_\eps^\rho(|[\hat{y}]|) \,d{\cal H}^1 \le {\cal H}^1(J_y) + C_1\rho,
\end{align}
where $\partial P = \bigcup_j \partial P_j$ and $ P(P_j,\Omega)$ denotes the perimeter of $P_j$ in $\Omega_\rho$.

It remains to confirm \eqref{rig-eq: main properties2}. First, (i) follows by  ${\cal H}^1(J_{\hat{y}} \cap (\Omega_y^H)^\circ)\le C_1$  (see \eqref{rig-eq: corweak**} and \eqref{rig-eq: V1.7}), ${\cal H}^1(\partial \hat{\Omega}) \le C_1$ (see Lemma \ref{rig-lemma: jordan}) and the fact that the ${\cal H}^1$-measure of the jump set added in the construction of $\hat{y}$ (see \eqref{rig-eq: rig10_*}) is controlled by ${\cal H}^1(\partial \hat{\Omega})$  and ${\cal H}^1(S)$. In view of \eqref{rig-eq: main properties}(ii)-(iv) (see also \eqref{rig-eq: compare}) the properties (ii)-(iv) already hold on the set $\hat{\Omega}$ for a sufficiently  large constant $C(\rho,q) = C(\rho)$. (Recall $q=q(h_*)$ and the definition $h_* = \rho$. See also Remark \ref{rig-rem: h_*}.)

Recall that $\Omega_{\rho} \setminus  \hat{\Omega} \subset \bigcup_{t}  \overline{Y_t}$. Repeating the arguments leading to \eqref{rig-eq: compare} we find by \eqref{rig-eq: ri2},  \eqref{rig-eq: ri1} and \eqref{rig-eq: rig10_*}
$$\sum\nolimits_j \Vert \hat{y} - (R_j\, \cdot +c_j)\Vert^2_{L^2(P_j\setminus \hat{\Omega})} \le C(\rho)\eps.$$
This gives  (ii). Moreover, as on each $Q^{ \varrho}(p) \subset P_j \setminus \hat{\Omega}$ we have $\nabla \hat{y} = R$ for some $R \in \text{im}_{\bar{R}_4}(\hat{\Omega})$  (see construction before \eqref{rig-eq: ri2}) we get 
\begin{align*}
\Vert \nabla \hat{y} - R_j\Vert^p_{L^p(P_j\setminus \hat{\Omega})}& \le   C(\rho)\Vert   \bar{R}_4 - R_j \Vert^p_{L^p(P_j\cap \hat{\Omega})}     \\& 
\le C(\rho)\Big(\Vert \nabla \hat{y} - \bar{R}_4\Vert^p_{L^p(P_j\cap \hat{\Omega})} + \Vert \nabla \hat{y} - R_j\Vert^p_{L^p(P_j\cap \hat{\Omega})}\Big)
\end{align*}
for $p=2,4$.  By \eqref{rig-eq: V1.9.2} and \eqref{rig-eq: compare} this yields
\begin{align*}
\sum\nolimits_j \Vert \nabla \hat{y} - R_j\Vert^4_{L^4(P_j\setminus \hat{\Omega})}  \le C(\rho)\eps, \ \ \ 
\sum\nolimits_j \Vert \nabla \hat{y} - R_j\Vert^2_{L^2(P_j\setminus \hat{\Omega})}  \le C(\rho)\eps^{1-\eta}.
\end{align*}
This together with \eqref{rig-eq: linearization} gives  (iii),(iv). \eop

Having completed the main rigidity result, we can now prove the linearized version. We may essentially follow the proof of Theorem \ref{rig-th: rigidity} with some minor changes. The proof, however, is considerably simpler as a lot of estimates and lemmas can be skipped.

\noindent {\em Proof of Theorem \ref{rig-th: rigidity_lin}.} We only give a short sketch of the proof. Define $y = \id +u$. As the approximation argument presented in the proof of Theorem \ref{rig-th: rigidity2} also holds in the SBD-setting, it again suffices to prove the result under the assumption that there is some $\tilde{\Omega}_u\in {\cal V}^s_\eps$ such that $u|_{\tilde{\Omega}_u} \in H^1(\tilde{\Omega}_u)$.  We skip Section \ref{rig-sec: subsub,  est-deriv} and always set $\hat{R}_i = \Id$ for $i=1,\ldots,4$. Similarly as in Lemma \ref{rig-lemma: weaklocA}  we find sets $\Omega_u$, $\Omega^H_u \in {\cal V}^{3\varrho}_{9k}$ for $k=\rho^{q-1}$, $\varrho= \rho^q$, as well as  mappings $\bar{A}_j: \Omega^H_u\to \R^{2 \times 2}_{\rm skew}$ and $\bar{c}_j: \Omega^H_u\to  \R^2$,  which are constant on $Q^{3\varrho}_j(p)$, $p \in I^{3\varrho}_j(\Omega^{3k})$, such that
\begin{align*}
\begin{split}
(i) & \ \ \Vert u  - (\bar{A}_j \, \cdot + \bar{c}_j)\Vert^2_{L^2(\Omega_u)} \le C C_\rho^2\varrho^2 (\alpha + \epsilon \Vert W \Vert_*),\\
(iii) & \ \ \Vert (\bar{A}_{j_1} \, \cdot + \bar{c}_{j_1})  - (\bar{A}_{j_2} \, \cdot + \bar{c}_{j_2})\Vert^2_{L^2(\Omega_u^H )}\le CC_\rho^2 \varrho^2 (\alpha + \epsilon \Vert W \Vert_*) 
\end{split}
\end{align*}
for $j_1,j_2 = 1,\ldots,4$, $j=1,\ldots,4$, where $\alpha = \Vert e(\nabla u) \Vert^2_{L^2(\tilde{\Omega}_u)}$  and $\epsilon = \hat{c}\rho^{-1}\eps$. This can be established following the lines of the proof of Lemma \ref{rig-lemma: weaklocA} with the difference that in \eqref{rig-eq: R,A diff} we do not replace $\Id + A$ by a different rigid motion $\bar{R}$, but proceed with $\Id + A$. Analogously, we find an extension $\Omega_u^H$ as constructed in Corollary \ref{rig-cor: weakA} and then we obtain the result up to a small set following the lines of Theorem \ref{rig-thm: V1}. Finally, the jump set and the extension to $\Omega_\rho$ may be constructed as in Section \ref{rig-sec: sub, proof-main}.  \eop


 \typeout{References}


\begin{thebibliography}{10}



\bibitem{Ambrosio-Coscia-Dal Maso:1997} 
{\sc L~ Ambrosio, A.~Coscia, G.~Dal Maso}.
\newblock {\em Fine properties of functions with bounded deformation}. 
\newblock Arch.\ Ration.\ Mech.\ Anal.\
\newblock {\bf 139} (1997), 201--238.






\bibitem{Ambrosio-Fusco-Pallara:2000} 
{\sc L.~Ambrosio, N.~Fusco, D.~Pallara}.
\newblock {\em Functions of bounded variation and free discontinuity problems}. 
\newblock Oxford University Press, Oxford 2000. 

\bibitem{Bellettini-Coscia-DalMaso:98}
{\sc G.~Bellettini, A.~Coscia, G.~Dal Maso}.
\newblock {\em  Compactness and lower semicontinuity properties in $SBD(\Omega)$}.
\newblock  Math.\ Zl.\
\newblock {\bf 228} (1998), 337--351.

\bibitem{Bourdin-Francfort-Marigo:2008}
{\sc B.~Bourdin, G.~A.~Francfort, J.~J.~Marigo}. 
\newblock {\em The variational approach to fracture}.
\newblock J.\ Elasticity\ 
\newblock {\bf 91} (2008), 5--148. 

\bibitem{Braides-Gelli:2002-2} 
{\sc A.~Braides, M.~S.~Gelli}. 
\newblock {\em Limits of discrete systems with long-range interactions}. 
\newblock J.\ Convex Anal.\
\newblock {\bf 9} (2002), 363--399. 







\bibitem{Chambolle:2003}
{\sc A.~Chambolle}. 
\newblock {\em A density result in two-dimensional linearized elasticity, and applications}.
\newblock Arch.\ Rat.\ Mech.\ Anal.\
\newblock {\bf 167} (2003), 167--211.



\bibitem{Chambolle:2004}
{\sc A.~Chambolle}. 
\newblock {\em An approximation result for special functions with bounded deformation}.
\newblock J.\ Math.\ Pures\ Appl. 
\newblock {\bf 83} (2004), 929--954. 



\bibitem{Chambolle-Giacomini-Ponsiglione:2007}
{\sc A.~Chambolle, A.~Giacomini, M.~Ponsiglione}. 
\newblock {\em Piecewise rigidity}.
\newblock J.\ Funct.\ Anal.\ Solids 
\newblock {\bf 244} (2007), 134--153. 


\bibitem{Conti-Dolzmann-Muller:14}
{\sc S.~Conti, G.~Dolzmann, and S.~M\"uller}.
\newblock {\em Korn's second inequality and geometric rigidity with mixed growth conditions}. 
\newblock Calc.\ Var.\
Partial\ Differential\ Equations\
\newblock {\bf 50} (2014), 437--454.


\bibitem{ContiFaracoMaggi:2005}
{\sc S.~Conti, D.~Faraco, F.~Maggi}. 
\newblock {\em A new approach to counterexamples to $L^1$ estimates: Korn's inequality, geometric rigidity, and regularity for gradients of separately convex functions}.
\newblock Arch.\ Rat.\ Mech.\ Anal.\
\newblock {\bf 175} (2005), 287--300. 



\bibitem{ContiSchweizer:06}
{\sc S.~Conti, B.~Schweizer}.
\newblock {\em Rigidity and Gamma convergence for solid-solid phase transitions with $SO(2)$-invariance}.
\newblock Comm.\ Pure\ Appl.\ Math.\
\newblock {\bf 59} (2006), 830--868.




\bibitem{Cortesani:1997}
{\sc G.~Cortesani}. 
\newblock {\em Strong approximation of GSBV functions by piecewise smooth functions}.
\newblock Ann.\ Univ.\ Ferrara\ Sez.\
\newblock {\bf 43} (1997), 27--49. 


\bibitem{Cortesani-Toader:1999}
{\sc G.~Cortesani, R.~Toader}. 
\newblock {\em A density result in SBV with respect to non-isotropic energies}.
\newblock Nonlinear Analysis
\newblock {\bf 38} (1999), 585--604. 




\bibitem{DalMaso:13}
{\sc G.~Dal Maso}.
\newblock {\em Generalized functions of bounded deformation}.
\newblock J.\ Eur.\ Math.\ Soc.\
\newblock {\bf 15} (2013), 1943--1997.

\bibitem{DalMasoNegriPercivale:02}
{\sc G.~Dal Maso, M.~Negri, D.~Percivale}.
\newblock {\em  Linearized elasticity as $\Gamma$-limit
of finite elasticity}.
\newblock  Set-valued\ Anal.\
\newblock {\bf 10} (2002), 165--183.

\bibitem{DalMaso-Toader:02}
{\sc G.~Dal Maso, R.~Toader}.
\newblock {\em A model for the quasistatic growth of brittle fractures: existence and approximation
results}. 
\newblock Arch.\ Rational\ Mech.\ Anal.,\ 
\newblock {\bf 162} (2002), 101--135.



\bibitem{DeGiorgi-Ambrosio:1988}
{\sc E.~De Giorgi, L.~Ambrosio}. 
\newblock {\em Un nuovo funzionale del calcolo delle variazioni}. 
\newblock Acc.\ Naz.\ Lincei, Rend.\ Cl.\ Sci.\ Fis.\ Mat.\ Natur.\ 
\newblock {\bf 82} (1988), 199--210. 




\bibitem{EvansGariepy92}
{\sc L.~C~Evans, R.~F.~Gariepy}. 
\newblock {\em Measure theory and fine properties of
functions}.
\newblock CRC Press, Boca Raton $\cdot$ London $\cdot$ New York $\cdot$ Washington,
D.C. 1992.

\bibitem{Focardi-Iurlano:13}
{\sc M.~Focardi, F.~Iurlano}.
\newblock {\em Asymptotic analysis of Ambrosio-
Tortorelli energies in linearized elasticity}. 
\newblock SIAM\ J.\ Math.\ Anal.\,
\newblock {\bf 46} (2014), 2936--2955.





\bibitem{Francfort-Larsen:2003}
{\sc G.~A.~Francfort, C,~J.~Larsen}. 
\newblock {\em Existence and convergence for quasi-static evolution in brittle fracture}.
\newblock Comm.\ Pure Appl.\ Math.\ 
\newblock {\bf 56} (2003), 1465--1500. 


\bibitem{Francfort-Marigo:1998}
{\sc G.~A.~Francfort, J.~J.~Marigo}. 
\newblock {\em Revisiting brittle fracture as an energy minimization problem}.
\newblock J.\ Mech.\ Phys.\ Solids 
\newblock {\bf 46} (1998), 1319--1342. 



\bibitem{Friedrich:15-1}
{\sc M.~Friedrich}.
\newblock {\em A Korn-Poincar\'e-type inequality for special functions of bounded deformation}. 
\newblock Preprint,\ 2015. 


\bibitem{Friedrich:15-2}
{\sc M.~Friedrich}.
\newblock {\em A derivation of linearized Griffith energies from nonlinear models}. 
\newblock Preprint,\ 2015. 



\bibitem{FriedrichSchmidt:2014.1}
{\sc M.~Friedrich, B.~Schmidt}.
\newblock {\em An analysis of crystal cleavage in the passage from atomistic models to continuum theory}. 
\newblock Arch.\ Rational\ Mech.\ Anal.,\  
\newblock published online 2014, doi:
10.1007/s00205-014-0833-y.


\bibitem{FriedrichSchmidt:2014.2}
{\sc M.~Friedrich, B.~Schmidt}.
\newblock {\em On a discrete-to-continuum convergence result for a two dimensional brittle material in the
small displacement regime}. 
\newblock Netw.\ Heterog.\ Media. In press. 


\bibitem{FrieseckeJamesMueller:02}
{\sc G.~Friesecke, R.~D.~James, S.~M{\"u}ller}.
\newblock {\em A theorem on geometric rigidity and the derivation of nonlinear plate theory from three-dimensional elasticity}. 
\newblock Comm.\ Pure Appl.\ Math.\ 
\newblock {\bf 55} (2002), 1461--1506. 

\bibitem{Griffith:1921} 
{\sc A.~A.~Griffith}. 
\newblock {\em The phenomena of rupture and flow in solids}.
\newblock Philos.\ Trans.\ R.\ Soc.\ London 
\newblock {\bf 221} (1921), 163--198. 


\bibitem{Iurlano:13}
{\sc F.~Iurlano}.
\newblock {\em A density result for GSBD and its application
to the approximation of brittle fracture energies}. 
\newblock Calc.\ Var. 
\newblock {\bf 51} (2014), 315--342. 


\bibitem{John:1961} 
{\sc F.~John}.
\newblock {\em Rotation and strain}. 
\newblock Comm.\ Pure.\ Appl.\ Math.\
\newblock {\bf 14} (1961), 391--413. 











\bibitem{Kohn:82}
{\sc R.~V.~Kohn}. 
\newblock {\em New integral estimates for deformations in terms of
their nonlinear strains}. 
\newblock Arch.\ Ration.\ Mech.\ Anal.\
\newblock {\bf 78} (1982), 131--172. 



\bibitem{Kristensen:1999}
{\sc J.~Kristensen}.
\newblock {\em Lower semicontinuity in spaces of weakly differentiable functions}. 
\newblock Math.\ Ann.\  
\newblock {\bf 313} (1999), 653--710.


\bibitem{Muller-Scardia-Zeppieri:14}
{\sc S.~M\"uller, L.~Scardia, and C.~I.~Zeppieri}.
\newblock {\em Geometric rigidity for incompatible fields and an application to strain-gradient plasticity}. 
\newblock Indiana\ Univ.\ Math.\ J.\ 
\newblock {\bf 63} (2014), 1365--1396.





\bibitem{NegriToader:2013}
{\sc M.~Negri, R.~Toader}.
\newblock {\em Scaling in fracture mechanics by Ba$\check z$ant's law: from finite to linearized elasticity}. 
\newblock Preprint\ SISSA,\ 
\newblock Trieste, 2013.



\bibitem{Reshetnyak:1961} 
{\sc Y.~G.~Reshetnyak}.
\newblock {\em  Liouville's theory on conformal mappings under minimial regularity assumptions}. 
\newblock Sibirskii\ Math.\ J.\
\newblock {\bf 8} (1967), 69--85. 


\bibitem{Schmidt:08}
{\sc B.~Schmidt}. 
\newblock {\em Linear $\Gamma$-limits of multiwell energies in nonlinear elasticity theory}. 
\newblock Continuum Mech.\ Thermodyn.\  
\newblock {\bf 20} (2008) 375--396. 




\bibitem{SchmidtFraternaliOrtiz:2009}
{\sc B.~Schmidt, F.~Fraternali, M.~Ortiz}. 
\newblock {\em Eigenfracture: an eigendeformation approach to variational fracture}. 
\newblock SIAM\ Mult.\ Model.\ Simul.\ 
\newblock {\bf 7} (2009), 1237--1266. 

\bibitem{Temam:85}
{\sc R.~Temam}.
\newblock {\em Mathematical Problems in Plasticity}.
\newblock Bordas, Paris 1985. 


\bibitem{Zhang:2004}
{\sc K.~Zhang}. 
\newblock {\em An approximation theorem for sequences of linear strains and its applications}. 
\newblock ESAIM\ Control\ Optim.\ Calc.\ Var.\
\newblock {\bf 10} (2004), 224--242. 





\end{thebibliography}
\end{document}